\documentclass[leqno,a4paper]{amsart}
\usepackage{amsmath}
\usepackage{epsfig}
\usepackage{amssymb,latexsym}
\usepackage{texdraw}
\usepackage{mathrsfs}

\newtheorem{thm}{\bf Theorem}[section]
\newtheorem{df}[thm]{\bf Definition}
\newtheorem{prop}[thm]{\bf Proposition}
\newtheorem{cor}[thm]{\bf Corollary}
\newtheorem{lem}[thm]{\bf Lemma}
\newtheorem{rem}[thm]{\bf Remark}
\newtheorem{ex}[thm]{\bf Example}

\newcommand{\wt}{{\rm wt}}
\newcommand{\Z}{\mathcal{Z}}

\newcommand{\cP}{\mathscr{P}}

\newcommand{\pf}{{\bfseries Proof. }}
\numberwithin{equation}{section}

\begin{document}
\title[Crystal bases and string functions]{Crystal bases of the Fock space representations and string functions}
\author[S.-J. Kang and J.-H. Kwon]{Seok-Jin Kang$^{*}$ and Jae-Hoon Kwon$^{\dagger}$}
\address{$^{*}$School of Mathematics\\ Korea Institute for Advanced Study \\
         207-43 Cheongryangri-dong Dongdaemun-gu \\
         Seoul 130-012, Korea }
\address{
         $^{\dagger}$Department of Mathematics \\ University of Seoul \\
         90 Cheonnong-dong, Dongdaemun-gu \\ Seoul,130-743 Korea}
\thanks{$^{*}$This research was supported by KRF Grant 2003-070-C00001.}
\thanks{$^{\dagger}$This research was supported by KOSEF Grant
R01-2003-000-10012-0.} \email{sjkang@kias.re.kr, jhkwon@uos.ac.kr}

\maketitle

\begin{abstract} Let $U_q(\frak{g})$ a the quantum affine algebra
of type $A_n^{(1)}$, $A_{2n-1}^{(2)}$, $A_{2n}^{(2)}$,
$B_n^{(1)}$, $D_n^{(1)}$ and $D_{n+1}^{(2)}$, and let
$\mathcal{F}(\Lambda)$ be the Fock space representation for a
level 1 dominant integral weight $\Lambda$. Using the crystal
basis of $\mathcal{F}(\Lambda)$ and its characterization in terms
of abacus, we construct an explicit bijection between the set of
weight vectors in $\mathcal{F}(\Lambda)_{\lambda-m\delta}$ ($m\geq
0$) for a maximal weight $\lambda$ and the set of certain ordered
sequences of partitions. As a corollary, we obtain the string
function of the basic representation $V(\Lambda)$.

\end{abstract}

\section{Introduction}
Let $\frak{g}$ be a classical affine Kac-Moody algebra of type
$A_n^{(1)}$, $A_{2n-1}^{(2)}$, $A_{2n}^{(2)}$, $B_n^{(1)}$,
$D_n^{(1)}$ and $D_{n+1}^{(2)}$, and let $\Lambda$ be a dominant
integral weight of level 1. The weight multiplicities of the basic
representation $V(\Lambda)$ can be explained in terms of the
string functions. By using, for example, the representations of
Virasoro algebras or modular forms, they are given as well-known
functions which arise naturally in combinatorics and number theory
(cf.\cite{Kac90}).

The purpose of this paper is to understand the combinatorics which
lies behind the string functions, that is, to interpret them in a
combinatorial way. To this end, we will use the representations of
the corresponding quantum affine algebra $U_q(\frak{g})$ and their
crystal bases. Also, instead of the basic representation, we will
use the Fock space representation of $U_q(\frak{g})$ with a nice
combinatorial realization of its crystal basis.

When $\frak{g}=A_n^{(1)}$, it can be explained very nicely by the
Misra and Miwa's Fock space representation \cite{MM}; They
introduced a $U_q(\frak{g})$-module $\mathcal{F}(\Lambda)$ which
is spanned by the set of all partitions. The submodule generated
by the empty partition is isomorphic to $V(\Lambda)$ whose crystal
is given by the set of $n$-reduced (or $n$-restricted) partitions.
We observe that a partition has a weight $\Lambda-m\delta$ ($m\geq
0$), as a crystal element of $\mathcal{F}(\Lambda)$, if and only
if it has empty $n$-core with $n$-weight $m$. Therefore, a weight
vector of $\mathcal{F}(\Lambda)_{\Lambda-m\delta}$ is uniquely
determined by its $n$-quotient, an $n$-tuple of partitions whose
sum is $m$. From this correspondence and the decomposition of
$\mathcal{F}(\Lambda)$ into irreducible highest weight modules, we
obtain the associated string function of $V(\Lambda)$ immediately.

Generalizing the notion of partitions in case of $A_n^{(1)}$, Kang
introduced an abstract crystal $\Z(\Lambda)$, the set of {\it
proper Young walls} which are collections of finite number of
blocks added on the ground state wall \cite{Ka2003}. Then he
showed that the connected component $\mathcal{Y}(\Lambda)$ of the
ground state wall is isomorphic to the crystal $B(\Lambda)$ of the
basic representation. Motivated by the work of Misra and Miwa, in
\cite{KK1,KK2}, we constructed a $U_q(\frak{g})$-module
$\mathcal{F}(\Lambda)$ having $\Z(\Lambda)$ as its crystal, which
decomposes as follows:
\begin{equation}
\mathcal{F}(\Lambda)\simeq \bigoplus_{m\geq 0}V(\Lambda-\epsilon
m\delta)^{\oplus p(m)},
\end{equation}
where $p(m)$ is the number of partitions of $m$ and $\epsilon=2$
if $\frak{g}=D_{n+1}^{(2)}$ and $\epsilon=1$ otherwise. We call
the $\mathcal{F}(\Lambda)$ the {\it Fock space representation} of
$U_q(\frak{g})$. This can be seen as a combinatorial realization
of the Fock space representation by Kashiwara, Miwa, Petersen and
Yung for level 1 case \cite{KMPY}, and it has a natural analogue
of the Lascoux-Leclerc-Thibon's algorithm for computing global
bases element of $V(\Lambda)$ \cite{LLT}. In this paper, we will
not define the module structure on $\mathcal{F}(\Lambda)$ since we
need only the crystal graph of it. But we will give another proof
for the decomposition of the crystal $\Z(\Lambda)$.

Therefore, we reduce our problem to characterizing the proper
Young walls in $\Z(\Lambda)_{\lambda-m\delta}$ ($m\geq 0$) for a
maximal weight $\lambda$. The main result in this paper is the
construction of a bijection between
$\Z(\Lambda)_{\lambda-m\delta}$ and a set of certain ordered
sequences of partitions whose generating function allows us to
recover the associated string function of $V(\Lambda)$. Also, the
bijection is given more explicitly when we take a particular
maximal weight (for example, $\Lambda=\lambda$), and it is
obtained by modifying the method of abacus which were used when
$\frak{g}=A_n^{(1)}$. We remark that it might be possible to give
a similar characterization of
$\mathcal{Y}(\Lambda)_{\lambda-m\delta}$. In fact, for
$A_n^{(1)}$, $A_{2n}^{(2)}$ and $D_{n+1}^{(2)}$, there exists a
bijection between $\mathcal{Y}(\Lambda)_{\lambda-m\delta}$ and a
set of ordered sequences of partitions (see \cite{LLT} for
$A_n^{(1)}$). Then, however, the bijections should include Weyl
group actions even in the case of $\lambda=\Lambda$, and hence
become more complicated than those for
$\Z(\Lambda)_{\lambda-m\delta}$. This is one of the reason we
prefer the Fock space representation rather than the basic
representation.

This paper is organized as follows: in Section 2, we recall the
notion of abstract crystals and proper Young walls. We refer the
reader to \cite{Kas91,Kas93} for a general exposition on crystal
bases and abstract crystals, and \cite{Ka2003} for a detailed
description and more examples of proper Young walls. In Section 3,
we review the results for $A_n^{(1)}$ which we mentioned before.
Then in the following sections, we define the abacus for each type
of $\frak{g}$, and then characterize the proper Young walls of
weight $\Lambda-m\delta$ from their bead configurations in the
abacus to obtain a bijection (Section 4 for $A_{2n}^{(2)}$,
$D_{n+1}^{(2)}$, Section 5 for $A_{2n-1}^{(2)}$, $D_{n+1}^{(1)}$
and Section 6 for $B_n^{(1)}$).\vskip 3mm

{\bf Acknowledgement} The result in this paper was announced at
KIAS International conference on Lie algebras and related topics
(Seoul, Oct.2003). We would like to thank S.-Y. Kang who kindly
taught us how to prove Lemma \ref{q(m)} during the conference.

\section{Affine crystals and Young walls}
Let $I=\{\,0,1,\cdots,n\,\}$ be an index set and let $(A,
P^{\vee}, P, \Pi^{\vee}, \Pi)$ be an affine Cartan datum where

(1) $A=(a_{ij})_{i,j\in I}$ is a generalized Cartan matrix of
affine type,

(2) $P^{\vee}=\mathbb{Z}h_0 \oplus \cdots \oplus \mathbb{Z}h_n
\oplus \mathbb{Z}d$ is the dual weight lattice,

(3) $\frak{h}=\mathbb{Q}\otimes_{\mathbb{Z}}P^{\vee}$ is the
Cartan subalgebra

(4) $P=\{\,\lambda\in \frak{h}^*\,|\,\lambda(P^{\vee})\subset
\mathbb{Z}\,\}$ is the weight lattice,

(5) $\Pi^{\vee}=\{\,h_i\,|\,i\in I\,\}$ is the set of simple
coroots,

(6) $\Pi=\{\,\alpha_i\,|\,i\in I\,\}$ is the set of simple roots.
\vskip 3mm

Let $\frak{g}$ be the corresponding affine Kac-Moody algebra. We
denote by $Q=\bigoplus_{i\in I}\mathbb{Z}\alpha_i$ the root
lattice, and set $Q_+=\bigoplus_{i\in I}\mathbb{Z}_{\geq
0}\alpha_i$, $Q_-=-Q_+$.
\begin{df}
{\rm An {\it {\rm (}affine{\rm )} crystal associated with} $(A,
P^{\vee}, P, \Pi^{\vee}, \Pi)$ is a set $B$ together with the maps
${\rm wt} : B \rightarrow P$, $\varepsilon_i : B \rightarrow
\mathbb{Z}\cup\{-\infty\}$, $\varphi_i : B \rightarrow
\mathbb{Z}\cup\{-\infty\}$, $\tilde{e}_i : B \rightarrow
B\cup\{0\}$, and $\tilde{f}_i : B \rightarrow B\cup\{0\}$
satisfying the following conditions\,:
\begin{itemize}

\item[(i)] for all $i\in I$, $b\in B$, we have
\begin{equation*}
\begin{aligned}\text{}
& \varphi_i(b) = \varepsilon_i(b) + {\rm wt}(b)(h_i), \\
& {\rm wt}(\tilde{e}_i b)={\rm wt}(b) + \alpha_i, \\
& {\rm wt}(\tilde{f}_i b)={\rm wt}(b) - \alpha_i,
\end{aligned}
\end{equation*}

\item[(ii)] if $\tilde{e}_i b \in B$, then
$$\varepsilon_i(\tilde{e}_i b) = \varepsilon_i(b) - 1, \quad
\varphi_i(\tilde{e}_i b) = \varphi_i(b) + 1,$$

\item[(iii)] if $\tilde{f}_i b \in B$, then
$$\varepsilon_i(\tilde{f}_i b) = \varepsilon_i(b) + 1, \quad
\varphi_i(\tilde{f}_i b) = \varphi_i(b) - 1,$$

\item[(iv)] $\tilde{f}_i b = b'$ if and only if $b = \tilde{e}_i
b'$ for all $i\in I$, $b, b' \in B$,

\item[(v)] if $\varepsilon_i(b) = - \infty$, then $\tilde{e}_i b =
\tilde{f}_i b = 0$.
\end{itemize}
}
\end{df}

From now on, we assume that $\frak{g}$ is of type $A_n^{(1)}$,
$A_{2n-1}^{(2)}$, $B_n^{(1)}$, $D_n^{(1)}$, $A_{2n}^{(2)}$ or
$D_{n+1}^{(2)}$. Suppose that we are given the following three
kinds of blocks;\vskip 3mm

\begin{center}
\begin{tabular}{c|c|c|c}
 shape & width & thickness & height  \\
\hline  \raisebox{-0.4\height}{
\begin{texdraw}
\drawdim em \setunitscale 0.1 \linewd 0.5 \move(-10 0)\lvec(0
0)\lvec(0 10)\lvec(-10 10)\lvec(-10 0) \move(0 0)\lvec(5 5)\lvec(5
15)\lvec(-5 15)\lvec(-10 10) \move(0 10)\lvec(5 15)
\end{texdraw}
} $=$ \raisebox{-0.4\height}{
\begin{texdraw}
\drawdim em \setunitscale 0.1 \linewd 0.5 \move(-10 0)\lvec(0
0)\lvec(0 10)\lvec(-10 10)\lvec(-10 0)
\end{texdraw}}
 & 1 & 1 & 1  \\
\raisebox{-0.4\height}{
\begin{texdraw}
\drawdim em \setunitscale 0.1 \linewd 0.5 \move(0 0)\lvec(10
0)\lvec(10 5)\lvec(0 5)\lvec(0 0) \move(10 0)\lvec(15 5)\lvec(15
10)\lvec(5 10)\lvec(0 5) \move(10 5)\lvec(15 10)
\end{texdraw}
} $=$ \raisebox{-0.4\height}{
\begin{texdraw}
\drawdim em \setunitscale 0.1 \linewd 0.5 \textref h:C v:C \move(0
0)\lvec(10 0)\lvec(10 5)\lvec(0 5)\lvec(0 0)
\end{texdraw}
}
& 1 & 1 & $\frac{1}{2}$  \\

\raisebox{-0.4\height}{
\begin{texdraw}
\drawdim em \setunitscale 0.1 \linewd 0.5 \move(-10 0)\lvec(0
0)\lvec(0 10)\lvec(-10 10)\lvec(-10 0) \move(0 0)\lvec(2.5
2.5)\lvec(2.5 12.5)\lvec(-7.5 12.5)\lvec(-10 10) \move(0
10)\lvec(2.5 12.5) \lpatt(0.3 1) \move(0 0)\lvec(-2.5
-2.5)\lvec(-12.5 -2.5)\lvec(-10 0)
\end{texdraw}
} $=$ \raisebox{-0.4\height}{
\begin{texdraw}
\drawdim em \setunitscale 0.1 \linewd 0.5 \move(-10 0)\lvec(0
0)\lvec(0 10)\lvec(-10 0) \lpatt(0.3 1)\move(-10 0)\lvec(-10
10)\lvec(0 10)
\end{texdraw}}
, \hskip 2mm \raisebox{-0.4\height}{
\begin{texdraw}
\drawdim em \setunitscale 0.1 \linewd 0.5 \move(-10 0)\lvec(0
0)\lvec(0 10)\lvec(-10 10)\lvec(-10 0) \move(0 0)\lvec(2.5
2.5)\lvec(2.5 12.5)\lvec(-7.5 12.5)\lvec(-10 10) \move(0
10)\lvec(2.5 12.5) \lpatt(0.3 1) \move(2.5 2.5)\lvec(5 5)\lvec(2.5
5)
\end{texdraw}
} $=$ \raisebox{-0.4\height}{
\begin{texdraw}
\drawdim em \setunitscale 0.1 \linewd 0.5 \move(-10 0)\lvec(-10
10)\lvec(0 10)\lvec(-10 0)\lpatt(0.3 1)\move(-10 0)\lvec(0
0)\lvec(0 10)
\end{texdraw}
}

& 1 & $\frac{1}{2}$ & 1  \\

\end{tabular}
\end{center}
\renewcommand{\arraystretch}{1}
\vskip 3mm and we give a coloring on them by the index set $I$.

Given $\frak{g}$ and a dominant integral weight $\Lambda$ of level
1, we fix a frame $Y_{\Lambda}$ called the {\it ground state wall}
of weight $\Lambda$, and we stack the above blocks on
$Y_{\Lambda}$ following the pattern depending on $\frak{g}$ and
$\Lambda$, which will be given in later sections. A collection of
finite number of colored blocks added on $Y_{\Lambda}$ is called a
{\it Young wall on $Y_{\Lambda}$} if the heights of its columns
are weakly decreasing from right to left. We often write
$Y=(y_k)_{k=1}^{\infty}$ as an infinite sequence of its columns
where the columns are enumerated from right to left. We define
$|Y|=(|y_k|)_{k=1}^{\infty}$ to be the sequence, where $|y_k|$ is
the number of blocks in the $k$th column of $Y$ (except the one in
$Y_{\Lambda}$), and call it the {\it associated partition of $Y$}.

A column of a Young wall is called a {\it full column} if the
block at the top is a unit cube. For type $A_n^{(1)}$, every Young
wall is defined to be {\it proper}. For type $A_{2n-1}^{(2)}$,
 $A_{2n}^{(2)}$, $B_n^{(1)}$, $D_n^{(1)}$ and $D_{n+1}^{(2)}$, a
Young wall is said to be {\it proper} if none of the full columns
have the same heights.  We denote by $\Z(\Lambda)$ the set of all
proper Young walls on $Y_{\Lambda}$.

Let $\delta=d_0\alpha_0+\cdots +d_n\alpha_n$ be the null root of
$\frak{g}$, and set $a_i=d_i$ if $\frak{g}\neq D^{(2)}_{n+1}$,
$a_i=2d_i$ if $\frak{g}=D^{(2)}_{n+1}$. The part of a column with
$a_i$-many $i$-blocks (or $i$-colored block) for each $i\in I$ in
some cyclic order is called a {\it $\delta$-column}. A
$\delta$-column in a proper Young wall is called {\it removable}
if it can be removed to yield another proper Young wall. A proper
Young wall $Y$ is said to be {\it reduced} if none of its columns
contain a removable $\delta$-column. We denote by ${\mathcal
Y}(\Lambda)$ the set of all reduced proper Young walls on
$Y_{\Lambda}$.
\begin{ex}{\rm
If $\frak{g}=A_5^{(2)}$ and $\Lambda=\Lambda_0$, then the Young
wall $Y$ given below is a proper Young wall in $\Z(\Lambda_0)$. It
is reduced since it contains no removable $\delta$-column.

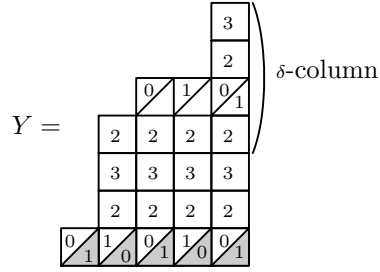
\begin{center}
$Y=$\raisebox{-0.5\height}{\begin{texdraw}\fontsize{7}{7}\drawdim
mm \setunitscale 0.5 \linewd 0.5

\move(51 70)\clvec(55 60)(55 40)(51 30)\htext(57
50){$\delta$-column}

\move(0 0)\lvec(10 0)\lvec(10 10)\lvec(0 10)\lvec(0 0)

\move(10 0)\lvec(20 0)\lvec(20 10)\lvec(10 10)\lvec(10 0)

\move(20 0)\lvec(30 0)\lvec(30 10)\lvec(20 10)\lvec(20 0)

\move(30 0)\lvec(40 0)\lvec(40 10)\lvec(30 10)\lvec(30 0)

\move(40 0)\lvec(50 0)\lvec(50 10)\lvec(40 10)\lvec(40 0)
\move(10 10)\lvec(20 10)\lvec(20 20)\lvec(10 20)\lvec(10
10)\htext(13 13){$2$}

\move(20 10)\lvec(30 10)\lvec(30 20)\lvec(20 20)\lvec(20
10)\htext(23 13){$2$}

\move(30 10)\lvec(40 10)\lvec(40 20)\lvec(30 20)\lvec(30
10)\htext(33 13){$2$}

\move(40 10)\lvec(50 10)\lvec(50 20)\lvec(40 20)\lvec(40
10)\htext(43 13){$2$}
\move(10 20)\lvec(20 20)\lvec(20 30)\lvec(10 30)\lvec(10
20)\htext(13 23){$3$}

\move(20 20)\lvec(30 20)\lvec(30 30)\lvec(20 30)\lvec(20
20)\htext(23 23){$3$}

\move(30 20)\lvec(40 20)\lvec(40 30)\lvec(30 30)\lvec(30
20)\htext(33 23){$3$}

\move(40 20)\lvec(50 20)\lvec(50 30)\lvec(40 30)\lvec(40
20)\htext(43 23){$3$}
\move(10 30)\lvec(20 30)\lvec(20 40)\lvec(10 40)\lvec(10
30)\htext(13 33){$2$}

\move(20 30)\lvec(30 30)\lvec(30 40)\lvec(20 40)\lvec(20
30)\htext(23 33){$2$}

\move(30 30)\lvec(40 30)\lvec(40 40)\lvec(30 40)\lvec(30
30)\htext(33 33){$2$}

\move(40 30)\lvec(50 30)\lvec(50 40)\lvec(40 40)\lvec(40
30)\htext(43 33){$2$}

\move(40 50)\lvec(50 50)\lvec(50 60)\lvec(40 60)\lvec(40
50)\htext(43 53){$2$}
\move(40 60)\lvec(50 60)\lvec(50 70)\lvec(40 70)\lvec(40
60)\htext(43 63){$3$}
\move(0 0)\lvec(10 10)\lvec(10 0)\lvec(0 0)\lfill f:0.8 \htext(1
5){$0$}\htext(6 1){$1$}

\move(10 0)\lvec(20 10)\lvec(20 0)\lvec(10 0)\lfill f:0.8
\htext(11 5){$1$}\htext(16 1){$0$}

\move(20 0)\lvec(30 10)\lvec(30 0)\lvec(20 0)\lfill f:0.8
\htext(22 5){$0$}\htext(26 2){$1$}

\move(30 0)\lvec(40 10)\lvec(40 0)\lvec(30 0)\lfill f:0.8
\htext(32 5){$1$}\htext(36 2){$0$}

\move(40 0)\lvec(50 10)\lvec(50 0)\lvec(40 0)\lfill f:0.8
\htext(42 5){$0$}\htext(46 2){$1$}
\move(20 40)\lvec(30 50)\lvec(20 50)\lvec(20 40)\htext(22 45){$0$}

\move(30 40)\lvec(40 50)\lvec(30 50)\lvec(30 40)\htext(32 45){$1$}

\move(40 40)\lvec(50 50)\lvec(50 40)\lvec(40 40)\lvec(40
50)\htext(42 45){$0$}\htext(46 42){$1$}
\end{texdraw}}
\end{center}\vskip 3mm

}
\end{ex}

\begin{df}{\rm Let $Y$ be a proper Young wall on $Y_{\Lambda}$.

{\rm (1)} An $i$-block in $Y$ is called a {\it removable
$i$-block} if $Y$ remains a proper Young wall after removing the
block.

{\rm (2)} A place in $Y$ is called an {\it admissible $i$-slot} if
one may add an $i$-block to obtain another proper Young wall.

{\rm (3)} A column in $Y$ is said to be {\it $i$-removable
{\rm(}resp. $i$-admissible{\rm )}} if there is a removable
$i$-block {\rm (}resp. an admissible $i$-slot{\rm )} in that
column.}
\end{df}

We now define the operators $\tilde{e}_i$, $\tilde{f}_i$ on
$\Z(\Lambda)$ as follows. Fix $i\in I$ and let
$Y=(y_k)_{k=1}^{\infty}$ be a proper Young wall on $Y_{\Lambda}$.

(1) To each column $y_k$ of $Y$, we assign
\begin{equation*}
\begin{cases}
-- & \text{if $y_k$ is twice $i$-removable,} \\
- & \text{if $y_k$ is once $i$-removable but not $i$-admissible,} \\
-+ & \text{if $y_k$ is once $i$-removable and once
$i$-admissible,} \\
+ & \text{if $y_k$ is once $i$-admissible but not $i$-removable,}
\\
++ & \text{if $y_k$ is twice $i$-admissible,} \\
\ \ \cdot & \text{otherwise}.
\end{cases}
\end{equation*}

(2) From this sequence of $+$'s and $-$'s, we cancel out every
$(+,-)$-pair to obtain a finite sequence of $-$'s followed by
$+$'s, reading from left to right. This finite sequence
$(-\cdots-,+\cdots+)$ is called the {\it $i$-signature} of $Y$.

(3) We define $\tilde{e}_i Y$ to be the proper Young wall obtained
from $Y$ by removing the $i$-block corresponding to the right-most
$-$ in the $i$-signature of $Y$. We define $\tilde{e}_i Y=0$ if
there is no $-$ in the $i$-signature of $Y$.

(4) We define $\tilde{f}_i Y$ to be the proper Young wall obtained
from $Y$ by adding an $i$-block to the column corresponding to the
left-most $+$ in the $i$-signature of $Y$. We define $\tilde{f}_i
Y=0$ if there is no $+$ in the $i$-signature of $Y$.

We also define
\begin{equation*}
\begin{split}
{\rm wt}(Y)&=\Lambda-\sum_{i\in I}k_i\alpha_i \in P, \\
\varepsilon_i(Y)&=\text{the number of $-$'s in the $i$-signature of $Y$}, \\
\varphi_i(Y)&=\text{the number of $+$'s in the $i$-signature of
$Y$},
\end{split}
\end{equation*}
where $k_i$ denotes the number of $i$-blocks in $Y$ that have been
added to $Y_{\Lambda}$. For a set $S$ consisting of some blocks in
$Y$, we define ${\rm cont}(S)=\sum_{i\in I}k_i\alpha_i\in Q_+$
where $k_i$ denotes the number of $i$-blocks in $S$, and call it
the {\it content of $S$}. For example, ${\rm cont}(Y)=\Lambda-{\rm
wt}(Y)$.

\begin{thm}{\rm (\cite{Ka2003})}
\begin{itemize}

\item[(1)] The set $\Z(\Lambda)$ together with
$\tilde{e}_i,\tilde{f}_i,{\rm wt},\varepsilon_i$ and $\varphi_i$
{\rm (}$i\in I${\rm )}, is an affine crystal.

\item[(2)] The set ${\mathcal Y}(\Lambda)$ is an affine subcrystal
of $\Z(\Lambda)$ and isomorphic to $B(\Lambda)$, where
$B(\Lambda)$ is the crystal of the basic representation
$V(\Lambda)$.
\end{itemize}
\qed
\end{thm}

Set $\epsilon=2$ if $\frak{g}$ is of type $D_{n+1}^{(2)}$, and
$\epsilon=1$ if otherwise. Then we have
\begin{cor}\label{decomp} There exists an isomorphism of
affine crystals
\begin{equation}
\Z(\Lambda)\simeq \bigoplus_{m\geq 0} B(\Lambda-\epsilon
m\delta)^{\oplus p(m)},
\end{equation}
where $B(\Lambda-m\delta)$ is the crystal of the highest weight
module $V(\Lambda-m\delta)$ and $p(m)$ is the number of partitions
of $m$.
\end{cor}
\pf For convenience, we assume that $\epsilon=1$, or $\frak{g}\neq
D_{n+1}^{(2)}$. Let $\lambda$ be a partition of $m$, that is, a
non-increasing sequence of positive integers $\lambda_1\geq
\lambda_2 \geq
 \cdots \geq \lambda_r$ whose sum is $m$. Let
$Y_{\Lambda,\lambda}$ be the proper Young wall which is obtained
by adding $\lambda_k$ many $\delta$-columns on the $k$th column of
$Y_{\Lambda}$. Then we can check that ${\rm
wt}(Y_{\Lambda,\lambda})=\Lambda-m\delta$ and $\tilde{e}_i
Y_{\Lambda,\lambda}=0$ for all $i\in I$ (see \cite{KK2}).

For each $Y\in{\mathcal Y}(\Lambda)$, define $S_{\lambda}(Y)$ to
be the proper Young wall which is obtained by adding the $k$th
column of $Y$ (except the block in $Y_{\Lambda}$) on the $k$th
column of $Y_{\Lambda,\lambda}$. Set
\begin{equation}
{\mathcal Y}(\Lambda,\lambda)=\{\,S_{\lambda}(Y)\,|\,Y\in{\mathcal
Y}(\Lambda)\,\}.
\end{equation}
Then it is not difficult to see that
\begin{equation}
\Z(\Lambda)=\bigsqcup_{\lambda}{\mathcal Y}(\Lambda,\lambda),
\end{equation}
where the union is taken over all partitions $\lambda$. Finally,
for $i\in I$, the map $S_{\lambda}:\mathcal{Y}(\Lambda)\rightarrow
\Z(\Lambda)$ commutes with $\tilde{e}_i$ and $\tilde{f}_i$, which
can be checked directly from the definitions of $\tilde{e}_i$ and
$\tilde{f}_i$. Therefore, we conclude that ${\mathcal
Y}(\Lambda,\lambda)$ is the connected component of
$Y_{\Lambda,\lambda}$, and hence it is isomorphic to
$B(\Lambda-m\delta)$. This completes the proof. \qed

\begin{rem}{\rm In \cite{KK2}, we constructed a
$U_q(\frak{g})$-module
\begin{equation}
{\mathcal F}(\Lambda)=\bigoplus_{Y\in \Z(\Lambda)}\mathbb{Q}(q)Y
\end{equation}
called the {\it Fock space representation}, and it was shown that
$\Z(\Lambda)$ is the crystal of ${\mathcal F}(\Lambda)$. From
this, the decomposition given in Corollary \ref{decomp} can be
obtained directly by verifying that
\begin{equation}
\{\,Y\in\Z(\Lambda)\,|\,\tilde{e}_iY=0\, \text{ for all $i\in
I$}\,\}=\{\,Y_{\Lambda,\lambda}\,|\,\lambda \, \text{:
partition}\,\}.
\end{equation}
}
\end{rem}
\vskip 3mm

Let $\lambda \in P$ be a maximal weight in $\Z(\Lambda)$, i.e.
$|\Z(\Lambda)_{\lambda}|\neq 0$ but
$|\Z(\Lambda)_{\lambda+\delta}|=0$. We define
\begin{equation}
\Sigma^{\Lambda}_{\lambda}(q)=\sum_{m\geq
0}|\Z(\Lambda)_{\lambda-m\delta}|q^m,
\end{equation}
and call it the {\it string function of $\mathcal{F}(\Lambda)$}
with respect to $\Lambda$ and $\lambda$. By Corollary
\ref{decomp}, we have
\begin{equation}
|\Z(\Lambda)_{\lambda-m\delta}|=\sum_{s\geq
0}p(s)|B(\Lambda-\epsilon s\delta)_{\lambda-m\delta}|
=\sum_{\epsilon s+t=m}p(s)|B(\Lambda)_{\lambda-t\delta}|.
\end{equation}
This implies that
\begin{equation}\label{Sigma}
\Sigma^{\Lambda}_{\lambda}(q)=\frac{1}{(q^{\epsilon})_{\infty}}\sigma^{\Lambda}_{\lambda}(q),
\end{equation}
where $(q)_{\infty}=\prod_{m\geq 1}(1-q^m)$ and
$\sigma^{\Lambda}_{\lambda}(q)=\sum_{m\geq
0}|B(\Lambda)_{\lambda-m\delta}|q^m$ is the string function of
$V(\Lambda)$ with respect to $\Lambda$ and $\lambda$
(cf.\cite{Kac90}). By using the Weyl group action on the crystal
$\Z(\Lambda)$ (cf.\cite{Kas93}), we have a bijection
\begin{equation}
w: \Z(\Lambda)_{\lambda-m\delta}\longrightarrow
\Z(\Lambda)_{w\lambda-m\delta},
\end{equation}
where $w$ is a Weyl group element. Thus when we consider the
weight multiplicities of the basic representation $V(\Lambda)$, it
is enough to compute
\begin{equation}
\begin{cases}
\Sigma^{\Lambda}_{\Lambda}(q) & \text{when $\frak{g}=A_n^{(1)}$,
$A_{2n-1}^{(2)}$, $A_{2n}^{(2)}$, $D_n^{(1)}$, $D_{n+1}^{(2)}$},
\\
\Sigma^{\Lambda_0}_{\Lambda_0}(q),
\Sigma^{\Lambda_0}_{\Lambda_1}(q) & \text{when
$\frak{g}=B_n^{(1)}$},
\end{cases}
\end{equation}
(see \cite{Kac90}).\vskip 3mm

\section{$A_{n-1}^{(1)}$-case $(n\geq 2)$}

This section is based on the arguments in \cite{JK}, and we will
rewrite them following our notations. First, let us recall some
basic terminologies (cf.\cite{Mac}). A {\it partition} is a
non-increasing sequence of non-negative integers,
$\lambda=(\lambda_k)_{k\geq 1}$ such that all but a finite number
of its terms are zero. Each $\lambda_k$ is called a {\it part of
$\lambda$} and the number of the non-zero parts is called the {\it
length} of $\lambda$, denoted by $\ell(\lambda)$. We also write
$\lambda=(1^{m_1},2^{m_2},3^{m_3},\cdots)$, where $m_i$ is the
number of the parts of $\lambda$ equal to $i$. We say that
$\lambda$ is a {\it partition of $m$} ($m\geq 0$) if $\sum_{k\geq
1}\lambda_k=m$ and write $|\lambda|=m$. When all non-zero parts of
$\lambda$ are distinct, $\lambda$ is said to be {\it strict}. For
each $m\geq 0$, let $\cP(m)$ be the set of partitions of $m$ and
set $\cP=\bigcup_{m\geq 0}\cP(m)$. We denote by $p(m)$ the number
of partitions of $m$ with $p(0)=1$ by convention.

A partition $\lambda=(\lambda_k)_{k\geq 1}$ is identified with a
{\it Young diagram} which is a collection of boxes stacked from
the bottom with $\lambda_k$ boxes in each $k$th column. We will
enumerate the columns of a Young diagram from {\it right to left}
so that the number of boxes are weakly decreasing from right to
left. Note that the boxes are stacked from the south-east corner
in our definition (cf.\cite{Mac}).

Let $\lambda$ be a Young diagram and let $b_1,\cdots,b_r$ be the
boxes in the main diagonal of $\lambda$, which are enumerated from
the south-east corner. Let $\lambda'_k$ ($1\leq k\leq r$) be the
number of boxes lying in the same row of $b_k$ and to the left of
$b_k$, and let $\lambda''_k$ ($1\leq k\leq r$) be the number of
boxes lying in the same column of $b_k$ and above $b_k$. Then we
have a pair of strict partitions $\lambda'=(\lambda'_1>\cdots
>\lambda'_r\geq 0)$ and $\lambda''=(\lambda''_1>\cdots
>\lambda''_r\geq 0)$
which are uniquely determined by $\lambda$. We may write
$\lambda=(\lambda'|\lambda'')$, which is called the {\it Frobenius
notation of $\lambda$}. For example, $(1^3,3,4^2)=((1,5)|(2,3))$.

Let $\frak{g}$ be of type $A_{n-1}^{(1)}$ $(n\geq 2)$. Fix a
dominant integral weight $\Lambda$ of level 1. The pattern for
$\Z(\Lambda)$ is given as follows: \vskip 5mm

\begin{center}
$\Lambda=\Lambda_i$\ \ \  \raisebox{-1\height}{\begin{texdraw}
\textref h:C v:C \fontsize{7}{7}\selectfont \drawdim em
\setunitscale 1.7 \move(0 0)\rlvec(-7.7 0) \move(0 1)\rlvec(-7.7
0) \move(0 2)\rlvec(-7.7 0) \move(0 3.5)\rlvec(-7.7 0) \move(0
4.5)\rlvec(-7.7 0) \move(0 5.5)\rlvec(-7.7 0) \move(0
6.5)\rlvec(-7.7 0) \move(0 0)\rlvec(0 6.7) \move(-1 0)\rlvec(0
6.7) \move(-2 0)\rlvec(0 6.7) \move(-3.5 0)\rlvec(0 6.7)
\move(-4.5 0)\rlvec(0 6.7) \move(-5.5 0)\rlvec(0 6.7) \move(-6.5
0)\rlvec(0 6.7) \move(-7.5 0)\rlvec(0 6.7) \move(-0.5 0.5)
\bsegment \htext(0 0){$i$} \htext(0 1){$i\!\!+\!\!1$} \vtext(0
2.25){$\cdots$} \htext(0 3.5){$n$} \htext(0 4.5){$0$} \htext(0
5.5){$1$} \htext(-1 0){$i\!\!-\!\!1$} \htext(-1 1){$i$} \vtext(-1
2.25){$\cdots$} \htext(-1 3.5){$n\!\!-\!\!1$} \htext(-1 4.5){$n$}
\htext(-1 5.5){$0$} \htext(-2.25 0){$\cdots$} \htext(-2.25
1){$\cdots$} \htext(-3.5 0){$2$} \htext(-3.5 1){$3$} \htext(-4.5
0){$1$} \htext(-4.5 1){$2$} \htext(-5.5 0){$0$} \htext(-5.5
1){$1$} \htext(-6.5 0){$n$} \htext(-6.5 1){$0$} \esegment
\end{texdraw}}
\end{center}%

\vskip 5mm For a given $Y=(y_k)_{k\geq 1}\in \Z(\Lambda)$, we may
identify $y_k$ with the number of blocks in the $k$th column of
$Y$. Then the set $\Z(\Lambda)$ can be identified with
$\mathscr{P}$, where $Y_{\Lambda}$ corresponds to the empty
partition.

We define the {\it abacus of type $A_{n-1}^{(1)}$} to be the
arrangement of positive integers in the following way:\vskip 5mm

\begin{center}
\begin{texdraw}
\drawdim em \setunitscale 0.13 \linewd 0.5 \fontsize{10}{10}

\htext(0 0){$1$}\htext(20 0){$2$}\htext(40 0){$\cdots$}\htext(60
0){$n-1$}\htext(90 0){$n$}

\htext(-7 -15){$n+1$}\htext(13 -15){$n+2$}\htext(40
-15){$\cdots$}\htext(57 -15){$2n-1$}\htext(88 -15){$2n$}

\htext(0 -30){$\vdots$}\htext(20 -30){$\vdots$}\htext(67
-30){$\vdots$}\htext(90 -30){$\vdots$}
\end{texdraw}
\end{center}

\vskip 5mm Let $R_k$ ($1\leq k\leq n$) be the set of integers $s
\equiv k \pmod n$ and call it the {\it $k$th runner}. Each
positive integer is called a {\it position}. Then we can put a
{\it bead}, denoted by $\bigcirc$, at each position and move a
bead along the runner which it belongs to. We suppose that there
is at most one bead at each position. So, we can move a bead at
$s$ one position up (resp. down) along the runner only when there
is no bead at $s-n$ (resp. $s+n$).

For a proper Young wall $Y=(y_k)_{k\geq 1}\in \Z(\Lambda)$, choose
an $r$ such that $y_k=0$ for all $k\geq r$. Consider the set of
$r$ distinct positive integers $\{\,z_k\,|\, z_k=y_k+r-k+1, 1\leq
k\leq r\,\}$. The {\it {\rm ($r$-)}bead configuration of $Y$} is
the set of $r$ beads assigned at the position $z_k$ ($1\leq k\leq
r$). Conversely, a given set of $r$ beads in the abacus represents
the $r$-bead configuration of a unique proper Young wall in
$\Z(\Lambda)$.

Consider a bead configuration of a proper Young wall
$Y\in\Z(\Lambda)$. Suppose that a bead $b$ is movable one position
up. Let $Z$ be the proper Young wall obtained by moving $b$ one
position up. Then we observe that $Y/Z$ forms a {\it border strip
of length $n$}; that is, $Y/Z$ is a skew Young diagram consisting
of $n$ boxes which is connected and contains no $2\times
2$-collection of boxes. Furthermore, if we consider a content of
$Y/Z$ following the pattern for $\Z(\Lambda)$, then it is
$\delta$. Therefore, we have $\wt(Z)=\wt(Y)+\delta$. Note that
this process of moving a bead one position up (equivalently,
removing a border strip of length $n$ from a Young diagram) is
reversible.

Let $\widetilde{Y}$ be a proper Young wall obtained by applying
the above processes until there is no bead movable up along the
runner. Note that $\widetilde{Y}$ is uniquely determined since
$\widetilde{Y}$ does not depend on the order in which we move up
the beads.  We call $\widetilde{Y}$ the {\it $n$-core of $Y$} and
denote it by ${\rm core}_{n}(Y)$. The total number of movements of
beads to obtain ${\rm core}_{n}(Y)$ from $Y$ is called the {\it
$n$-weight of $Y$} and denoted by ${\wt}_n(Y)$.

\begin{ex}{\rm Suppose $Y=(1,3,5,6)$ and $n=4$. Then the bead
configuration of $Y$ with $5$ beads is as follows:

\begin{center}
$Y=$\raisebox{-0.5\height}{
\begin{texdraw}\drawdim
em \setunitscale 0.13 \linewd 0.5

\move(10 0)\lvec(20 0)\lvec(20 10)\lvec(10 10)\lvec(10 0)

\move(20 0)\lvec(30 0)\lvec(30 10)\lvec(20 10)\lvec(20 0)

\move(30 0)\lvec(40 0)\lvec(40 10)\lvec(30 10)\lvec(30 0)

\move(40 0)\lvec(50 0)\lvec(50 10)\lvec(40 10)\lvec(40 0)

\move(20 10)\lvec(30 10)\lvec(30 20)\lvec(20 20)\lvec(20 10)

\move(30 10)\lvec(40 10)\lvec(40 20)\lvec(30 20)\lvec(30 10)

\move(40 10)\lvec(50 10)\lvec(50 20)\lvec(40 20)\lvec(40 10)

\move(20 20)\lvec(30 20)\lvec(30 30)\lvec(20 30)\lvec(20 20)

\move(30 20)\lvec(40 20)\lvec(40 30)\lvec(30 30)\lvec(30 20)

\move(40 20)\lvec(50 20)\lvec(50 30)\lvec(40 30)\lvec(40 20)

\move(30 30)\lvec(40 30)\lvec(40 40)\lvec(30 40)\lvec(30 30)

\move(40 30)\lvec(50 30)\lvec(50 40)\lvec(40 40)\lvec(40 30)

\move(30 40)\lvec(40 40)\lvec(40 50)\lvec(30 50)\lvec(30 40)

\move(40 40)\lvec(50 40)\lvec(50 50)\lvec(40 50)\lvec(40 40)

\move(40 50)\lvec(50 50)\lvec(50 60)\lvec(40 60)\lvec(40 50)

\linewd 0.3 \move(15 5)\lvec(25 5)\lvec(25 15)\lvec(35 15)

\move(25 25)\lvec(45 25)\lvec(45 35)

\move(35 35)\lvec(35 45)\lvec(45 45)\lvec(45 55)
\end{texdraw}}\ \ \ \ \ \ $\longleftrightarrow$ \ \ \ \ \ \raisebox{-0.5\height}{\begin{texdraw} \drawdim em \setunitscale
0.16 \linewd 0.5 \fontsize{8}{8}

\htext(10 0){$1$}\htext(20 0){$2$}\htext(30 0){$3$}\htext(40
0){$4$}

\htext(10 -10){$5$}\htext(20 -10){$6$}\htext(30 -10){$7$}\htext(40
-10){$8$}

\htext(10 -20){$9$}\htext(20 -20){$10$}\htext(30
-20){$11$}\htext(40 -20){$12$}

\htext(10 -30){$\vdots$}\htext(20 -30){$\vdots$}\htext(30
-30){$\vdots$}\htext(40 -30){$\vdots$}

\move(11 1)\linewd .3 \lcir r:3.5 \move(31 1)\linewd .3 \lcir
r:3.5 \move(21 -9)\linewd .3 \lcir r:3.5 \move(11 -19)\linewd .3
\lcir r:3.5 \move(32 -19)\linewd .3 \lcir r:3.5
\end{texdraw}}\vskip 5mm

 ${\rm core}_4(Y)=$ \raisebox{-0.5\height}{
\begin{texdraw}\drawdim
em \setunitscale 0.13 \linewd 0.5

\move(30 0)\lvec(40 0)\lvec(40 10)\lvec(30 10)\lvec(30 0)

\move(40 0)\lvec(50 0)\lvec(50 10)\lvec(40 10)\lvec(40 0)

\move(40 10)\lvec(50 10)\lvec(50 20)\lvec(40 20)\lvec(40 10)

\end{texdraw}}\ \ , \ \ \ \ \  ${\rm wt}_4(Y)=3$.\end{center}}
\end{ex}\vskip 5mm

Note that two proper Young walls in $\mathcal{Z}(\Lambda)$ have
the same $n$-core if and only if they have the same weight (or
content) (see \cite{JK}). For example, suppose that the weight of
$Y\in \mathcal{Z}(\Lambda)$ is $\Lambda-m\delta$ for some $m\geq
0$. On the other hand, consider
$Z=(1^{mn})\in\mathcal{Z}(\Lambda)$ as a partition. Since ${\rm
wt}(Z)=\Lambda-m\delta$ and ${\rm core}_n(Z)=Y_{\Lambda}$ (or
$(0)$), we have ${\rm core}_n(Y)=Y_{\Lambda}$. Conversely, any
$Y\in \mathcal{Z}(\Lambda)$ with ${\rm core}_n(Y)=Y_{\Lambda}$ has
weight $\Lambda-m\delta$ for some $m\geq 0$.

Therefore, for $m\geq 0$, we have
\begin{equation}
\Z(\Lambda)_{\Lambda-m\delta}=\{\,Y\in\Z(\Lambda)\,|\,{\rm core}_n
(Y)=Y_{\Lambda},\ {\rm wt}_n(Y)=m\,\}.
\end{equation}

Fix $m\geq 0$. Let $Y$ be a proper Young wall in
$\Z(\Lambda)_{\Lambda-m\delta}$. Set $M=mn$. Consider the $M$-bead
configuration of $Y$. Since ${\rm core}_n{Y}=Y_{\Lambda}$, it is
easy to see that there are $m$ beads in each runner $R_k$ ($1\leq
k\leq n$). For each $1\leq k\leq n$, let
$\{\,b^{(k)}_1,\cdots,b^{(k)}_m\,\}$ be the set of $m$ beads in
$R_k$ enumerated from the bottom. Suppose that $b^{(k)}_i$ is
located at $N^{(k)}_i$. Put $p^{(k)}_i=\frac{N^{(k)}_i-k}{n}$.
Then $p^{(k)}_i\geq m-i$ and $\{\,p^{(k)}_i-m+i \,|\,1\leq i\leq
m\,\}$ forms a unique partition $\lambda^{(k)}$ whose sum is the
number of all possible movements of beads in $R_k$ to obtain the
core of $Y$. Hence, we obtain an $n$-tuple of partitions
$\pi(Y)=(\lambda^{(1)},\cdots,\lambda^{(n)})$ with
$\sum_{i=1}^n|\lambda^{(i)}|=m$, which is called the {\it
$n$-quotient of $Y$}. Conversely, for a given $n$-tuple of
partitions $(\lambda^{(1)},\cdots,\lambda^{(n)})$ with
$\sum_{i=1}^n|\lambda^{(i)}|=m$, we can place $m$ beads in each
runner $R_k$ whose corresponding partition is $\lambda^{(k)}$
($1\leq k\leq n$). Then the resulting unique proper Young wall $Y$
is in $\Z(\Lambda)_{\Lambda-m\delta}$ with
$\pi(Y)=(\lambda^{(1)},\cdots,\lambda^{(n)})$. Therefore, we
obtain
\begin{thm}\label{ZAn1}{\rm (\cite{JK})}
For $m\geq 0$, the map
\begin{equation}
\pi : \Z(\Lambda)_{\Lambda-m\delta} \longrightarrow \cP^{(n)}(m)
\end{equation}
is a bijection where
$$\cP^{(n)}(m)=\{\,(\lambda^{(1)},\cdots,\lambda^{(n)})\,|\,\lambda^{(i)}\in\cP,\,
\sum_{i=1}^n|\lambda^{(i)}|=m\,\}.$$ \qed
\end{thm}

\begin{rem}\label{2core}{\rm  In \cite{JK}, the bijection in the above theorem is given in a
more general form, that is, a bijection between $\cP^{(n)}(m)$ and
the set of all partitions with a given $n$-core and an $n$-weight
$m$. This implies that a partition is uniquely determined by its
$n$-core and $n$-quotient.}
\end{rem}

By Theorem \ref{ZAn1} and \eqref{Sigma}, we obtain
$\sigma^{\Lambda}_{\Lambda}(q)=1/(q)_{\infty}^{n-1}$, or
\begin{equation}\label{2red}
\dim V(\Lambda)_{\Lambda-m\delta} =\sum_{\sum_{i=1}^{n-1}
m_i=m}p(m_1)\cdots p(m_{n-1})
\end{equation}
for $m\geq 0$ (cf. \cite{Kac90}).

A partition is called {\it $n$-reduced} (or {\it $n$-restricted})
if the difference of any two adjacent columns is less than $n$.
Note that $\mathcal{Y}(\Lambda)$ is the set of all $n$-reduced
partitions. In particular, we denote by $\mathscr{DP}_0(m)$ the
set of  all $2$-reduced partitions whose $2$-core is empty and
$2$-weight is $m$ (or the set of all strict partitions with empty
$2$-core and $2$-weight $m$), and set
$\mathscr{DP}_0=\bigcup_{m\geq 0}\mathscr{DP}_{0}(m)$. By
\eqref{2red}, we have
\begin{equation}\label{YAn1}
|\mathscr{DP}_0(m)|=p(m).
\end{equation}

\begin{rem}{\rm
The set of $n$-cores is in one-to-one correspondence with
$W/W_{\Lambda}$ where $W$ is the Weyl group of type
$A_{n-1}^{(1)}$ and $W_{\Lambda}$ is the stabilizer subgroup of
$\Lambda$. In fact, the bijection is given by ${\rm
wt}(Y)=w\Lambda$. Moreover, using this fact and the Weyl group
action on the crystal $\Z(\Lambda)$, Lascoux, Leclerc and Thibon
described a bijection between
$\mathcal{Y}(\Lambda)_{\Lambda-m\delta}$ and $\cP^{(n-1)}(m)$
\cite{LLT}.}
\end{rem}

\section{$A_{2n}^{(2)}$, $D_{n+1}^{(2)}$-case}

Suppose that $\frak{g}$ is of type $A_{2n}^{(2)}$ ($n\geq 1$) or
$D_{n+1}^{(2)}$ ($n\geq 2$) and $\Lambda$ is a dominant integral
weight of level $1$. The patterns for $\Z(\Lambda)$ are given as
follows:\vskip 5mm

$A_{2n}^{(2)}$,

\begin{center}
$\Lambda=\Lambda_0$ : \raisebox{-1\height}{
\begin{texdraw} \textref h:C v:C \fontsize{6}{6}\selectfont \drawdim mm
\setunitscale 5 \move(0 0)\lvec(-4 0)\lvec(-4 0.5)\lvec(0
0.5)\ifill f:0.7 \move(0 0)\rlvec(-4.3 0) \move(0 0.5)\rlvec(-4.3
0) \move(0 1)\rlvec(-4.3 0) \move(0 2)\rlvec(-4.3 0) \move(0
3.5)\rlvec(-4.3 0) \move(0 4.5)\rlvec(-4.3 0) \move(0
6)\rlvec(-4.3 0) \move(0 7)\rlvec(-4.3 0) \move(0 7.5)\rlvec(-4.3
0) \move(0 8)\rlvec(-4.3 0) \move(0 9)\rlvec(-4.3 0) \move(0
0)\rlvec(0 9.3) \move(-1 0)\rlvec(0 9.3) \move(-2 0)\rlvec(0 9.3)
\move(-3 0)\rlvec(0 9.3) \move(-4 0)\rlvec(0 9.3) \htext(-0.5
0.25){$0$} \htext(-1.5 0.25){$0$} \htext(-2.5 0.25){$0$}
\htext(-3.5 0.25){$0$} \htext(-0.5 0.75){$0$} \htext(-1.5
0.75){$0$} \htext(-2.5 0.75){$0$} \htext(-3.5 0.75){$0$}
\htext(-0.5 1.5){$1$} \htext(-1.5 1.5){$1$} \htext(-2.5 1.5){$1$}
\htext(-3.5 1.5){$1$} \htext(-0.5 4){$n$} \htext(-1.5 4){$n$}
\htext(-2.5 4){$n$} \htext(-3.5 4){$n$} \htext(-0.5 6.5){$1$}
\htext(-1.5 6.5){$1$} \htext(-2.5 6.5){$1$} \htext(-3.5 6.5){$1$}
\htext(-0.5 7.25){$0$} \htext(-1.5 7.25){$0$} \htext(-2.5
7.25){$0$} \htext(-3.5 7.25){$0$} \htext(-0.5 7.75){$0$}
\htext(-1.5 7.75){$0$} \htext(-2.5 7.75){$0$} \htext(-3.5
7.75){$0$} \htext(-0.5 8.5){$1$} \htext(-1.5 8.5){$1$} \htext(-2.5
8.5){$1$} \htext(-3.5 8.5){$1$} \vtext(-0.5 2.75){$\cdots$}
\vtext(-1.5 2.75){$\cdots$} \vtext(-2.5 2.75){$\cdots$}
\vtext(-3.5 2.75){$\cdots$} \vtext(-0.5 5.25){$\cdots$}
\vtext(-1.5 5.25){$\cdots$} \vtext(-2.5 5.25){$\cdots$}
\vtext(-3.5 5.25){$\cdots$}
\end{texdraw}}
\end{center}\vskip 5mm

$D_{n+1}^{(2)}$,
\begin{center}
$\Lambda=\Lambda_0$ : \raisebox{-1\height}{
\begin{texdraw} \textref h:C v:C \fontsize{6}{6}\selectfont \drawdim mm
\setunitscale 5 \move(0 0)\lvec(-4 0)\lvec(-4 0.5)\lvec(0
0.5)\ifill f:0.7 \move(0 0)\rlvec(-4.3 0) \move(0 0.5)\rlvec(-4.3
0) \move(0 1)\rlvec(-4.3 0) \move(0 2)\rlvec(-4.3 0) \move(0
3.5)\rlvec(-4.3 0) \move(0 4.5)\rlvec(-4.3 0) \move(0
4.5)\rlvec(-4.3 0) \move(0 5)\rlvec(-4.3 0) \move(0
5.5)\rlvec(-4.3 0) \move(0 6.5)\rlvec(-4.3 0) \move(0
8)\rlvec(-4.3 0) \move(0 9)\rlvec(-4.3 0) \move(0 9.5)\rlvec(-4.3
0) \move(0 10)\rlvec(-4.3 0) \move(0 11)\rlvec(-4.3 0) \move(0
0)\rlvec(0 11.3) \move(-1 0)\rlvec(0 11.3) \move(-2 0)\rlvec(0
11.3) \move(-3 0)\rlvec(0 11.3) \move(-4 0)\rlvec(0 11.3)
\htext(-0.5 0.25){$0$} \htext(-1.5 0.25){$0$} \htext(-2.5
0.25){$0$} \htext(-3.5 0.25){$0$} \htext(-0.5 0.75){$0$}
\htext(-1.5 0.75){$0$} \htext(-2.5 0.75){$0$} \htext(-3.5
0.75){$0$} \htext(-0.5 1.5){$1$} \htext(-1.5 1.5){$1$} \htext(-2.5
1.5){$1$} \htext(-3.5 1.5){$1$} \vtext(-0.5 2.75){$\cdots$}
\vtext(-1.5 2.75){$\cdots$} \vtext(-2.5 2.75){$\cdots$}
\vtext(-3.5 2.75){$\cdots$} \htext(-0.5 4){$n\!\!-\!\!1$}
\htext(-1.5 4){$n\!\!-\!\!1$} \htext(-2.5 4){$n\!\!-\!\!1$}
\htext(-3.5 4){$n\!\!-\!\!1$} \htext(-0.5 4.75){$n$} \htext(-1.5
4.75){$n$} \htext(-2.5 4.75){$n$} \htext(-3.5 4.75){$n$}
\htext(-0.5 5.25){$n$} \htext(-1.5 5.25){$n$} \htext(-2.5
5.25){$n$} \htext(-3.5 5.25){$n$} \htext(-0.5 6){$n\!\!-\!\!1$}
\htext(-1.5 6){$n\!\!-\!\!1$} \htext(-2.5 6){$n\!\!-\!\!1$}
\htext(-3.5 6){$n\!\!-\!\!1$} \vtext(-0.5 7.25){$\cdots$}
\vtext(-1.5 7.25){$\cdots$} \vtext(-2.5 7.25){$\cdots$}
\vtext(-3.5 7.25){$\cdots$} \htext(-0.5 8.5){$1$} \htext(-1.5
8.5){$1$} \htext(-2.5 8.5){$1$} \htext(-3.5 8.5){$1$} \htext(-0.5
9.25){$0$} \htext(-1.5 9.25){$0$} \htext(-2.5 9.25){$0$}
\htext(-3.5 9.25){$0$} \htext(-0.5 9.75){$0$} \htext(-1.5
9.75){$0$} \htext(-2.5 9.75){$0$} \htext(-3.5 9.75){$0$}
\htext(-0.5 10.5){$1$} \htext(-1.5 10.5){$1$} \htext(-2.5
10.5){$1$} \htext(-3.5 10.5){$1$}
\end{texdraw}}\hskip 2cm
$\Lambda=\Lambda_1$ : \raisebox{-1\height}{
\begin{texdraw}\textref h:C v:C \fontsize{6}{6}\selectfont \drawdim
mm \setunitscale 5 \move(0 0)\lvec(-4 0)\lvec(-4 0.5)\lvec(0
0.5)\ifill f:0.7 \move(0 0)\rlvec(-4.3 0) \move(0 0.5)\rlvec(-4.3
0) \move(0 1)\rlvec(-4.3 0) \move(0 2)\rlvec(-4.3 0) \move(0
3.5)\rlvec(-4.3 0) \move(0 4.5)\rlvec(-4.3 0) \move(0
4.5)\rlvec(-4.3 0) \move(0 5)\rlvec(-4.3 0) \move(0
5.5)\rlvec(-4.3 0) \move(0 6.5)\rlvec(-4.3 0) \move(0
8)\rlvec(-4.3 0) \move(0 9)\rlvec(-4.3 0) \move(0 9.5)\rlvec(-4.3
0) \move(0 10)\rlvec(-4.3 0) \move(0 11)\rlvec(-4.3 0) \move(0
0)\rlvec(0 11.3) \move(-1 0)\rlvec(0 11.3) \move(-2 0)\rlvec(0
11.3) \move(-3 0)\rlvec(0 11.3) \move(-4 0)\rlvec(0 11.3)
\htext(-0.5 0.25){$n$} \htext(-1.5 0.25){$n$} \htext(-2.5
0.25){$n$} \htext(-3.5 0.25){$n$} \htext(-0.5 0.75){$n$}
\htext(-1.5 0.75){$n$} \htext(-2.5 0.75){$n$} \htext(-3.5
0.75){$n$} \htext(-0.5 1.5){$n\!\!-\!\!1$} \htext(-1.5
1.5){$n\!\!-\!\!1$} \htext(-2.5 1.5){$n\!\!-\!\!1$} \htext(-3.5
1.5){$n\!\!-\!\!1$} \vtext(-0.5 2.75){$\cdots$} \vtext(-1.5
2.75){$\cdots$} \vtext(-2.5 2.75){$\cdots$} \vtext(-3.5
2.75){$\cdots$} \htext(-0.5 4){$1$} \htext(-1.5 4){$1$}
\htext(-2.5 4){$1$} \htext(-3.5 4){$1$} \htext(-0.5 4.75){$0$}
\htext(-1.5 4.75){$0$} \htext(-2.5 4.75){$0$} \htext(-3.5
4.75){$0$} \htext(-0.5 5.25){$0$} \htext(-1.5 5.25){$0$}
\htext(-2.5 5.25){$0$} \htext(-3.5 5.25){$0$} \htext(-0.5 6){$1$}
\htext(-1.5 6){$1$} \htext(-2.5 6){$1$} \htext(-3.5 6){$1$}
\vtext(-0.5 7.25){$\cdots$} \vtext(-1.5 7.25){$\cdots$}
\vtext(-2.5 7.25){$\cdots$} \vtext(-3.5 7.25){$\cdots$}
\htext(-0.5 8.5){$n\!\!-\!\!1$} \htext(-1.5 8.5){$n\!\!-\!\!1$}
\htext(-2.5 8.5){$n\!\!-\!\!1$} \htext(-3.5 8.5){$n\!\!-\!\!1$}
\htext(-0.5 9.25){$n$} \htext(-1.5 9.25){$n$} \htext(-2.5
9.25){$n$} \htext(-3.5 9.25){$n$} \htext(-0.5 9.75){$n$}
\htext(-1.5 9.75){$n$} \htext(-2.5 9.75){$n$} \htext(-3.5
9.75){$n$} \htext(-0.5 10.5){$n\!\!-\!\!1$} \htext(-1.5
10.5){$n\!\!-\!\!1$} \htext(-2.5 10.5){$n\!\!-\!\!1$} \htext(-3.5
10.5){$n\!\!-\!\!1$}
\end{texdraw}}
\end{center}\vskip 5mm

Set
\begin{equation}
\begin{split}
&\ell=
\begin{cases}
2n+1 & \text{if $\frak{g}=A_{2n}^{(2)}$}, \\
n+1 & \text{if $\frak{g}=D_{n+1}^{(2)}$},
\end{cases}
\hskip 3mm L=
\begin{cases}
\ell & \text{if $\frak{g}=A_{2n}^{(2)}$}, \\
2\ell & \text{if $\frak{g}=D_{n+1}^{(2)}$},
\end{cases}\\
&\epsilon=L/\ell
\end{split}
\end{equation}
Note that $L$ is the number of blocks in a $\delta$-column. We
define the {\it abacus of type $A_{2n}^{(2)}$ {\rm (}resp.
$D_{n+1}^{(2)}${\rm )}} to be the arrangement of positive integers
in the following way :\vskip 5mm

$A_{2n}^{(2)}$
\begin{center}
\begin{texdraw}
\drawdim em \setunitscale 0.13 \linewd 0.5 \fontsize{10}{10}

\htext(0 0){$1$}\htext(20 0){$2$}\htext(40 0){$\cdots$}\htext(60
0){$\ell-1$}\htext(90 0){$\ell$}

\htext(-7 -15){$\ell+1$}\htext(13 -15){$\ell+2$}\htext(40
-15){$\cdots$}\htext(57 -15){$2\ell-1$}\htext(88 -15){$2\ell$}

\htext(0 -30){$\vdots$}\htext(20 -30){$\vdots$}\htext(67
-30){$\vdots$}\htext(90 -30){$\vdots$}
\end{texdraw}
\end{center}

\vskip 5mm

$D_{n+1}^{(2)}$
\begin{center}
\begin{texdraw}
\drawdim em \setunitscale 0.13 \linewd 0.5 \fontsize{10}{10}

\htext(0 0){$1$}\htext(20 0){$\cdots$}\htext(40
0){$\ell-1$}\htext(60 0){$\ell+1$}\htext(80 0){$\cdots$}\htext(100
0){$2\ell-1$}

\htext(-10 -20){$2\ell+1$}\htext(20 -20){$\cdots$}\htext(37
-20){$3\ell-1$}\htext(58 -20){$3\ell+1$}\htext(80
-20){$\cdots$}\htext(100 -20){$4\ell-1$}

\htext(0 -35){$\vdots$}\htext(45 -35){$\vdots$}\htext(67
-35){$\vdots$}\htext(108 -35){$\vdots$}

\htext(130 10){$\ell$}\htext(128 0){$2\ell$}\htext(128
-10){$3\ell$}\htext(128 -20){$4\ell$}\htext(130 -35){$\vdots$}
\end{texdraw}
\end{center}\vskip 5mm

Let $R_k$ ($1\leq k<L$, $k\neq\ell$ ) be the set of all integers
$s \equiv k \pmod L$ and let $R_{\ell}$ be the set of all integers
$s \equiv 0 \pmod \ell$.  We call $R_k$ the {\it $k$th runner}.
There are $2n+1$ runners in each abacus. The rules of placing and
moving beads in $R_k$ ($k\neq \ell $) are the same as in the case
of type $A_n^{(1)}$, and we say that $R_k$ is of {\it type $\rm
I$}. On the other hand, we suppose that there can be more than one
bead at each position in $R_{\ell}$. We denote $k$ beads at $s$ by
\raisebox{-.25 \height}{\begin{texdraw} \drawdim em \setunitscale
0.13 \linewd 0.5 \move(2 2)\lcir r:4 \htext(0 0){$s$}\htext(7
3){$_k$}
\end{texdraw}}. Moreover, if $b$ is a bead at $m(\neq \ell)$ in
$R_{\ell}$, then we can always move up (resp. down) $b$ along the
runner by increasing the number of beads at $m-\ell$ (resp.
$m+\ell$) by one, and decreasing the number of beads at $m$ by
one. We say that $R_{\ell}$ is of {\it type $\rm II$}.

For $Y\in\Z(\Lambda)$, let $|Y|=(|y_k|)_{k\geq 1}$ be its
associated partition. Let $\{\,|y_1|,\cdots,|y_r|\,\}$ be the set
of all non-zero parts in $|Y|$. Then by definition of
$\Z(\Lambda)$, the numbers $|y_k|$'s ($1\leq k\leq r$) are
distinct except when $|y_k|\equiv 0\pmod \ell$. We define the {\it
bead configuration of $Y$} to be the set of $r$ beads
$b_1,\cdots,b_r$ in the above abacus where $b_k$ is placed at
$|y_k|$. Note that $Y$ is uniquely determined by its associated
partition, and hence by its bead configuration.

\begin{ex}\label{ex4.1}{\rm Suppose that $\frak{g}=A_{4}^{(2)}$. Then \vskip
3mm
\begin{center}
$Y=$\raisebox{-0.5\height}{
\begin{texdraw}
\drawdim em \setunitscale 0.13 \linewd 0.5

\move(-10 0)\lvec(0 0)\lvec(0 10)\lvec(-10 10)\lvec(-10
0)\htext(-7 1){\tiny $0$}

\move(0 0)\lvec(10 0)\lvec(10 10)\lvec(0 10)\lvec(0 0)\htext(3
1){\tiny $0$}

\move(10 0)\lvec(20 0)\lvec(20 10)\lvec(10 10)\lvec(10 0)\htext(13
1){\tiny $0$}

\move(20 0)\lvec(30 0)\lvec(30 10)\lvec(20 10)\lvec(20 0)\htext(23
1){\tiny $0$}

\move(30 0)\lvec(40 0)\lvec(40 10)\lvec(30 10)\lvec(30 0)\htext(33
1){\tiny $0$}

\move(40 0)\lvec(50 0)\lvec(50 10)\lvec(40 10)\lvec(40 0)\htext(43
1){\tiny $0$}

\move(50 0)\lvec(60 0)\lvec(60 10)\lvec(50 10)\lvec(50 0)\htext(53
1){\tiny $0$}
\move(-10 0)\lvec(0 0)\lvec(0 5)\lvec(-10 5)\lvec(-10 0)\lfill
f:0.8 \htext(-7 6){\tiny $0$}

\move(0 0)\lvec(10 0)\lvec(10 5)\lvec(0 5)\lvec(0 0)\lfill f:0.8
\htext(3 6){\tiny $0$}

\move(10 0)\lvec(20 0)\lvec(20 5)\lvec(10 5)\lvec(10 0)\lfill
f:0.8 \htext(13 6){\tiny $0$}

\move(20 0)\lvec(30 0)\lvec(30 5)\lvec(20 5)\lvec(20 0)\lfill
f:0.8 \htext(23 6){\tiny $0$}

\move(30 0)\lvec(40 0)\lvec(40 5)\lvec(30 5)\lvec(30 0)\lfill
f:0.8 \htext(33 6){\tiny $0$}

\move(40 0)\lvec(50 0)\lvec(50 5)\lvec(40 5)\lvec(40 0)\lfill
f:0.8 \htext(43 6){\tiny $0$}

\move(50 0)\lvec(60 0)\lvec(60 5)\lvec(50 5)\lvec(50 0)\lfill
f:0.8 \htext(53 6){\tiny $0$}

\move(0 10)\lvec(10 10)\lvec(10 20)\lvec(0 20)\lvec(0 10)\htext(3
13){$_1$}

\move(10 10)\lvec(20 10)\lvec(20 20)\lvec(10 20)\lvec(10
10)\htext(13 13){$_1$}

\move(20 10)\lvec(30 10)\lvec(30 20)\lvec(20 20)\lvec(20
10)\htext(23 13){$_1$}

\move(30 10)\lvec(40 10)\lvec(40 20)\lvec(30 20)\lvec(30
10)\htext(33 13){$_1$}

\move(40 10)\lvec(50 10)\lvec(50 20)\lvec(40 20)\lvec(40
10)\htext(43 13){$_1$}

\move(50 10)\lvec(60 10)\lvec(60 20)\lvec(50 20)\lvec(50
10)\htext(53 13){$_1$}
\move(0 20)\lvec(10 20)\lvec(10 30)\lvec(0 30)\lvec(0 20)\htext(3
23){$_2$}

\move(10 20)\lvec(20 20)\lvec(20 30)\lvec(10 30)\lvec(10
20)\htext(13 23){$_2$}

\move(20 20)\lvec(30 20)\lvec(30 30)\lvec(20 30)\lvec(20
20)\htext(23 23){$_2$}

\move(30 20)\lvec(40 20)\lvec(40 30)\lvec(30 30)\lvec(30
20)\htext(33 23){$_2$}

\move(40 20)\lvec(50 20)\lvec(50 30)\lvec(40 30)\lvec(40
20)\htext(43 23){$_2$}

\move(50 20)\lvec(60 20)\lvec(60 30)\lvec(50 30)\lvec(50
20)\htext(53 23){$_2$}

\move(0 30)\lvec(10 30)\lvec(10 40)\lvec(0 40)\lvec(0 30)\htext(3
33){$_1$}

\move(10 30)\lvec(20 30)\lvec(20 40)\lvec(10 40)\lvec(10
30)\htext(13 33){$_1$}

\move(20 30)\lvec(30 30)\lvec(30 40)\lvec(20 40)\lvec(20
30)\htext(23 33){$_1$}

\move(30 30)\lvec(40 30)\lvec(40 40)\lvec(30 40)\lvec(30
30)\htext(33 33){$_1$}

\move(40 30)\lvec(50 30)\lvec(50 40)\lvec(40 40)\lvec(40
30)\htext(43 33){$_1$}

\move(50 30)\lvec(60 30)\lvec(60 40)\lvec(50 40)\lvec(50
30)\htext(53 33){$_1$}
\move(10 40)\lvec(20 40)\lvec(20 45)\lvec(10 45)\lvec(10
40)\htext(13 41){\tiny $0$}

\move(20 40)\lvec(30 40)\lvec(30 45)\lvec(20 45)\lvec(20
40)\htext(23 41){\tiny $0$}

\move(30 40)\lvec(40 40)\lvec(40 45)\lvec(30 45)\lvec(30
40)\htext(33 41){\tiny $0$}

\move(40 40)\lvec(50 40)\lvec(50 45)\lvec(40 45)\lvec(40
40)\htext(43 41){\tiny $0$}

\move(50 40)\lvec(60 40)\lvec(60 45)\lvec(50 45)\lvec(50
40)\htext(53 41){\tiny $0$}

\move(50 45)\lvec(60 45)\lvec(60 50)\lvec(50 50)\lvec(50
45)\htext(53 46){\tiny $0$}

\move(50 50)\lvec(50 60)\lvec(60 60)\lvec(60 50)\htext(53
53){$_1$}

\move(50 60)\lvec(50 70)\lvec(60 70)\lvec(60 60)\htext(53
63){$_2$}

\move(50 70)\lvec(50 80)\lvec(60 80)\lvec(60 70)\htext(53
73){$_1$}

\move(50 80)\lvec(50 85)\lvec(60 85)\lvec(60 80)\htext(53
81){\tiny $0$}

\move(40 45)\lvec(50 45)\lvec(50 50)\lvec(40 50)\lvec(40
45)\htext(43 46){\tiny $0$}

\move(40 50)\lvec(40 60)\lvec(50 60)\lvec(50 50)\htext(43
53){$_1$}
\end{texdraw}}\ \ \ \ \ $\longleftrightarrow$ \ \ \ \ \
\raisebox{-0.5\height}{\begin{texdraw} \drawdim em \setunitscale
0.16 \linewd 0.5 \fontsize{8}{8}

\htext(10 0){$1$}\htext(20 0){$2$}\htext(30 0){$3$}\htext(40
0){$4$}\htext(50 0){$5$}

\htext(10 -10){$6$}\htext(20 -10){$7$}\htext(30 -10){$8$}\htext(40
-10){$9$}\htext(49 -10){$10$}

\htext(10 -20){$\vdots$}\htext(20 -20){$\vdots$}\htext(30
-20){$\vdots$}\htext(40 -20){$\vdots$}\htext(50 -20){$\vdots$}

\move(11 1)\linewd .3 \lcir r:3.5 \move(41 1)\linewd .3 \lcir
r:3.5 \move(51 1)\linewd .3 \lcir r:3.5 \move(21 -9)\linewd .3
\lcir r:3.5 \fontsize{4}{4}\move(54 4)\htext{$3$} \move(52
-9)\linewd .3 \lcir r:3.5
\end{texdraw}}\ \ \ .
\end{center}\vskip 3mm}
\end{ex}

We will describe an algorithm in the abacus which is an analogue
of removing a border strip in the abacus of type $A_{n-1}^{(1)}$,
and then characterize $\mathcal{Z}(\Lambda)_{\Lambda-m\delta}$
($m\geq 0$) in terms of bead configurations.

\begin{lem}\label{abacusA2n2}
Let $Y$ be a proper Young wall in $\Z(\Lambda)$ and let $Y'$ be
the proper Young wall obtained by applying one of the following
processes to the bead configuration of $Y$:
\begin{itemize}
\item[$(B_1)$] if $b$ is a bead at $s$ in a runner of type I and
there is no bead at $s-L$, then move $b$ one position up,

\item[$(B_2)$] if $b$ is a bead at $s$ in $R_{\ell}$ and
$s\neq\ell$, then move $b$ one position up,

\item[$(B_3)$] if $b$ and $b'$ are beads at $s$ and $L-s$ {\rm
($1\leq s\leq n$)}, respectively, then remove $b$ and $b'$
simultaneously,

\item[$(B_4)$] if $b$ is a bead at $\ell$, then remove $b$.
\end{itemize}
\noindent Then we have
\begin{equation}
\wt(Y')=
\begin{cases}
\wt(Y)+\epsilon\delta & \text{if $Y'$ is obtained by $(B_i)$
{\rm ($i=1,2,3$)}}, \\
\wt(Y)+\delta & \text{if $Y'$ is obtained by $(B_4)$}.
\end{cases}
\end{equation}

\end{lem}
\pf  If we apply $(B_i)$ ($i=1,2,3$) to $Y$, then $Y'$ is obtained
by removing some $L$ blocks from $Y$, say
$\{\,b_1,\cdots,b_{L}\,\}$. Let $i_{k}$ ($1\leq k\leq L$) be the
color of $b_k$. It follows directly from the pattern for
$\Z(\Lambda)$ that $\sum_{k=1}^{L}\alpha_{i_k}=\epsilon\delta$ and
therefore $\wt(Y')=\wt(Y)+\epsilon\delta$. The proof is similar
when $Y'$ is obtained by $(B_4)$. \qed\vskip 3mm

\begin{ex}{\rm Let $Y$ be the proper Young wall in Example
\ref{ex4.1}.

(1) Apply $(B_1)$ to \raisebox{-0.25 \height}{\begin{texdraw}
\drawdim em \setunitscale 0.16 \linewd 0.5 \fontsize{8}{8}
\htext(0 0){$7$} \move(1 1)\linewd .3 \lcir r:3
\end{texdraw}}. This means that we remove a $\delta$-column in the
2nd column of $Y$ and then shift the blocks which are placed in
the left of the $\delta$-column, to the right as far as possible.

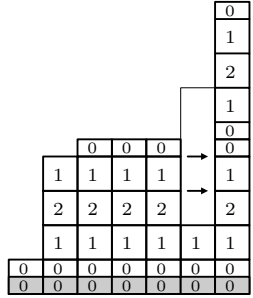
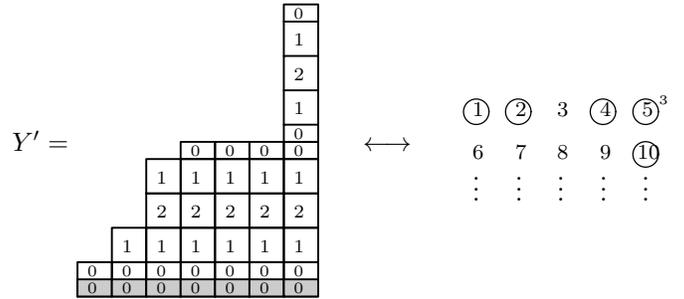
\begin{center}
\raisebox{-0.5\height}{
\begin{texdraw}
\drawdim em \setunitscale 0.13 \linewd 0.5 \arrowheadtype t:F
\arrowheadsize l:2 w:2

\move(42 30)\avec(48 30) \move(42 40)\avec(48 40)

\move(-10 0)\lvec(0 0)\lvec(0 10)\lvec(-10 10)\lvec(-10
0)\htext(-7 1){\tiny $0$}

\move(0 0)\lvec(10 0)\lvec(10 10)\lvec(0 10)\lvec(0 0)\htext(3
1){\tiny $0$}

\move(10 0)\lvec(20 0)\lvec(20 10)\lvec(10 10)\lvec(10 0)\htext(13
1){\tiny $0$}

\move(20 0)\lvec(30 0)\lvec(30 10)\lvec(20 10)\lvec(20 0)\htext(23
1){\tiny $0$}

\move(30 0)\lvec(40 0)\lvec(40 10)\lvec(30 10)\lvec(30 0)\htext(33
1){\tiny $0$}

\move(40 0)\lvec(50 0)\lvec(50 10)\lvec(40 10)\lvec(40 0)\htext(43
1){\tiny $0$}

\move(50 0)\lvec(60 0)\lvec(60 10)\lvec(50 10)\lvec(50 0)\htext(53
1){\tiny $0$}
\move(-10 0)\lvec(0 0)\lvec(0 5)\lvec(-10 5)\lvec(-10 0)\lfill
f:0.8 \htext(-7 6){\tiny $0$}

\move(0 0)\lvec(10 0)\lvec(10 5)\lvec(0 5)\lvec(0 0)\lfill f:0.8
\htext(3 6){\tiny $0$}

\move(10 0)\lvec(20 0)\lvec(20 5)\lvec(10 5)\lvec(10 0)\lfill
f:0.8 \htext(13 6){\tiny $0$}

\move(20 0)\lvec(30 0)\lvec(30 5)\lvec(20 5)\lvec(20 0)\lfill
f:0.8 \htext(23 6){\tiny $0$}

\move(30 0)\lvec(40 0)\lvec(40 5)\lvec(30 5)\lvec(30 0)\lfill
f:0.8 \htext(33 6){\tiny $0$}

\move(40 0)\lvec(50 0)\lvec(50 5)\lvec(40 5)\lvec(40 0)\lfill
f:0.8 \htext(43 6){\tiny $0$}

\move(50 0)\lvec(60 0)\lvec(60 5)\lvec(50 5)\lvec(50 0)\lfill
f:0.8 \htext(53 6){\tiny $0$}

\move(0 10)\lvec(10 10)\lvec(10 20)\lvec(0 20)\lvec(0 10)\htext(3
13){$_1$}

\move(10 10)\lvec(20 10)\lvec(20 20)\lvec(10 20)\lvec(10
10)\htext(13 13){$_1$}

\move(20 10)\lvec(30 10)\lvec(30 20)\lvec(20 20)\lvec(20
10)\htext(23 13){$_1$}

\move(30 10)\lvec(40 10)\lvec(40 20)\lvec(30 20)\lvec(30
10)\htext(33 13){$_1$}

\move(40 10)\lvec(50 10)\lvec(50 20)\lvec(40 20)\lvec(40
10)\htext(43 13){$_1$}

\move(50 10)\lvec(60 10)\lvec(60 20)\lvec(50 20)\lvec(50
10)\htext(53 13){$_1$}
\move(0 20)\lvec(10 20)\lvec(10 30)\lvec(0 30)\lvec(0 20)\htext(3
23){$_2$}

\move(10 20)\lvec(20 20)\lvec(20 30)\lvec(10 30)\lvec(10
20)\htext(13 23){$_2$}

\move(20 20)\lvec(30 20)\lvec(30 30)\lvec(20 30)\lvec(20
20)\htext(23 23){$_2$}

\move(30 20)\lvec(40 20)\lvec(40 30)\lvec(30 30)\lvec(30
20)\htext(33 23){$_2$}

\move(50 20)\lvec(60 20)\lvec(60 30)\lvec(50 30)\lvec(50
20)\htext(53 23){$_2$}

\move(0 30)\lvec(10 30)\lvec(10 40)\lvec(0 40)\lvec(0 30)\htext(3
33){$_1$}

\move(10 30)\lvec(20 30)\lvec(20 40)\lvec(10 40)\lvec(10
30)\htext(13 33){$_1$}

\move(20 30)\lvec(30 30)\lvec(30 40)\lvec(20 40)\lvec(20
30)\htext(23 33){$_1$}

\move(30 30)\lvec(40 30)\lvec(40 40)\lvec(30 40)\lvec(30
30)\htext(33 33){$_1$}


\move(50 30)\lvec(60 30)\lvec(60 40)\lvec(50 40)\lvec(50
30)\htext(53 33){$_1$}
\move(10 40)\lvec(20 40)\lvec(20 45)\lvec(10 45)\lvec(10
40)\htext(13 41){\tiny $0$}

\move(20 40)\lvec(30 40)\lvec(30 45)\lvec(20 45)\lvec(20
40)\htext(23 41){\tiny $0$}

\move(30 40)\lvec(40 40)\lvec(40 45)\lvec(30 45)\lvec(30
40)\htext(33 41){\tiny $0$}

\move(50 40)\lvec(60 40)\lvec(60 45)\lvec(50 45)\lvec(50
40)\htext(53 41){\tiny $0$}

\move(50 45)\lvec(60 45)\lvec(60 50)\lvec(50 50)\lvec(50
45)\htext(53 46){\tiny $0$}

\move(50 50)\lvec(50 60)\lvec(60 60)\lvec(60 50)\htext(53
53){$_1$}

\move(50 60)\lvec(50 70)\lvec(60 70)\lvec(60 60)\htext(53
63){$_2$}

\move(50 70)\lvec(50 80)\lvec(60 80)\lvec(60 70)\htext(53
73){$_1$}

\move(50 80)\lvec(50 85)\lvec(60 85)\lvec(60 80)\htext(53
81){\tiny $0$}

\linewd 0.2 \move(40 20)\lvec(40 60)\lvec(50 60)\lvec(50 20)
\end{texdraw}}\end{center}

Therefore, we have
\begin{center}
$Y'=$\raisebox{-0.5\height}{
\begin{texdraw}
\drawdim em \setunitscale 0.13 \linewd 0.5

\move(-10 0)\lvec(0 0)\lvec(0 10)\lvec(-10 10)\lvec(-10
0)\htext(-7 1){\tiny $0$}

\move(0 0)\lvec(10 0)\lvec(10 10)\lvec(0 10)\lvec(0 0)\htext(3
1){\tiny $0$}

\move(10 0)\lvec(20 0)\lvec(20 10)\lvec(10 10)\lvec(10 0)\htext(13
1){\tiny $0$}

\move(20 0)\lvec(30 0)\lvec(30 10)\lvec(20 10)\lvec(20 0)\htext(23
1){\tiny $0$}

\move(30 0)\lvec(40 0)\lvec(40 10)\lvec(30 10)\lvec(30 0)\htext(33
1){\tiny $0$}

\move(40 0)\lvec(50 0)\lvec(50 10)\lvec(40 10)\lvec(40 0)\htext(43
1){\tiny $0$}

\move(50 0)\lvec(60 0)\lvec(60 10)\lvec(50 10)\lvec(50 0)\htext(53
1){\tiny $0$}
\move(-10 0)\lvec(0 0)\lvec(0 5)\lvec(-10 5)\lvec(-10 0)\lfill
f:0.8 \htext(-7 6){\tiny $0$}

\move(0 0)\lvec(10 0)\lvec(10 5)\lvec(0 5)\lvec(0 0)\lfill f:0.8
\htext(3 6){\tiny $0$}

\move(10 0)\lvec(20 0)\lvec(20 5)\lvec(10 5)\lvec(10 0)\lfill
f:0.8 \htext(13 6){\tiny $0$}

\move(20 0)\lvec(30 0)\lvec(30 5)\lvec(20 5)\lvec(20 0)\lfill
f:0.8 \htext(23 6){\tiny $0$}

\move(30 0)\lvec(40 0)\lvec(40 5)\lvec(30 5)\lvec(30 0)\lfill
f:0.8 \htext(33 6){\tiny $0$}

\move(40 0)\lvec(50 0)\lvec(50 5)\lvec(40 5)\lvec(40 0)\lfill
f:0.8 \htext(43 6){\tiny $0$}

\move(50 0)\lvec(60 0)\lvec(60 5)\lvec(50 5)\lvec(50 0)\lfill
f:0.8 \htext(53 6){\tiny $0$}

\move(0 10)\lvec(10 10)\lvec(10 20)\lvec(0 20)\lvec(0 10)\htext(3
13){$_1$}

\move(10 10)\lvec(20 10)\lvec(20 20)\lvec(10 20)\lvec(10
10)\htext(13 13){$_1$}

\move(20 10)\lvec(30 10)\lvec(30 20)\lvec(20 20)\lvec(20
10)\htext(23 13){$_1$}

\move(30 10)\lvec(40 10)\lvec(40 20)\lvec(30 20)\lvec(30
10)\htext(33 13){$_1$}

\move(40 10)\lvec(50 10)\lvec(50 20)\lvec(40 20)\lvec(40
10)\htext(43 13){$_1$}

\move(50 10)\lvec(60 10)\lvec(60 20)\lvec(50 20)\lvec(50
10)\htext(53 13){$_1$}
%

\move(10 20)\lvec(20 20)\lvec(20 30)\lvec(10 30)\lvec(10
20)\htext(13 23){$_2$}

\move(20 20)\lvec(30 20)\lvec(30 30)\lvec(20 30)\lvec(20
20)\htext(23 23){$_2$}

\move(30 20)\lvec(40 20)\lvec(40 30)\lvec(30 30)\lvec(30
20)\htext(33 23){$_2$}

\move(40 20)\lvec(50 20)\lvec(50 30)\lvec(40 30)\lvec(40
20)\htext(43 23){$_2$}

\move(50 20)\lvec(60 20)\lvec(60 30)\lvec(50 30)\lvec(50
20)\htext(53 23){$_2$}
%


\move(10 30)\lvec(20 30)\lvec(20 40)\lvec(10 40)\lvec(10
30)\htext(13 33){$_1$}

\move(20 30)\lvec(30 30)\lvec(30 40)\lvec(20 40)\lvec(20
30)\htext(23 33){$_1$}

\move(30 30)\lvec(40 30)\lvec(40 40)\lvec(30 40)\lvec(30
30)\htext(33 33){$_1$}

\move(40 30)\lvec(50 30)\lvec(50 40)\lvec(40 40)\lvec(40
30)\htext(43 33){$_1$}

\move(50 30)\lvec(60 30)\lvec(60 40)\lvec(50 40)\lvec(50
30)\htext(53 33){$_1$}
%

\move(20 40)\lvec(30 40)\lvec(30 45)\lvec(20 45)\lvec(20
40)\htext(23 41){\tiny $0$}

\move(30 40)\lvec(40 40)\lvec(40 45)\lvec(30 45)\lvec(30
40)\htext(33 41){\tiny $0$}

\move(40 40)\lvec(50 40)\lvec(50 45)\lvec(40 45)\lvec(40
40)\htext(43 41){\tiny $0$}

\move(50 40)\lvec(60 40)\lvec(60 45)\lvec(50 45)\lvec(50
40)\htext(53 41){\tiny $0$}

\move(50 45)\lvec(60 45)\lvec(60 50)\lvec(50 50)\lvec(50
45)\htext(53 46){\tiny $0$}

\move(50 50)\lvec(50 60)\lvec(60 60)\lvec(60 50)\htext(53
53){$_1$}

\move(50 60)\lvec(50 70)\lvec(60 70)\lvec(60 60)\htext(53
63){$_2$}

\move(50 70)\lvec(50 80)\lvec(60 80)\lvec(60 70)\htext(53
73){$_1$}

\move(50 80)\lvec(50 85)\lvec(60 85)\lvec(60 80)\htext(53
81){\tiny $0$}


\end{texdraw}}\ \ \ \ \ $\longleftrightarrow$ \ \ \ \ \
\raisebox{-0.5\height}{\begin{texdraw} \drawdim em \setunitscale
0.16 \linewd 0.5 \fontsize{8}{8}

\htext(10 0){$1$}\htext(20 0){$2$}\htext(30 0){$3$}\htext(40
0){$4$}\htext(50 0){$5$}

\htext(10 -10){$6$}\htext(20 -10){$7$}\htext(30 -10){$8$}\htext(40
-10){$9$}\htext(49 -10){$10$}

\htext(10 -20){$\vdots$}\htext(20 -20){$\vdots$}\htext(30
-20){$\vdots$}\htext(40 -20){$\vdots$}\htext(50 -20){$\vdots$}

\move(11 1)\linewd .3 \lcir r:3 \move(41 1)\linewd .3 \lcir r:3
\move(51 1)\linewd .3 \lcir r:3 \move(21 1)\linewd .3 \lcir r:3
\fontsize{4}{4}\move(54 3)\htext{$3$} \move(51 -9)\linewd .3 \lcir
r:3
\end{texdraw}}\ \ \ .
\end{center}\vskip 3mm

(2) Similarly, if we apply $(B_3)$ to \raisebox{-0.25
\height}{\begin{texdraw} \drawdim em \setunitscale 0.16 \linewd
0.5 \fontsize{8}{8} \htext(0 0){$1$} \move(1 1)\linewd .3 \lcir
r:3
\end{texdraw}} and \raisebox{-0.25
\height}{\begin{texdraw} \drawdim em \setunitscale 0.16 \linewd
0.5 \fontsize{8}{8} \htext(0 0){$4$} \move(1 1)\linewd .3 \lcir
r:3
\end{texdraw}}, then we have
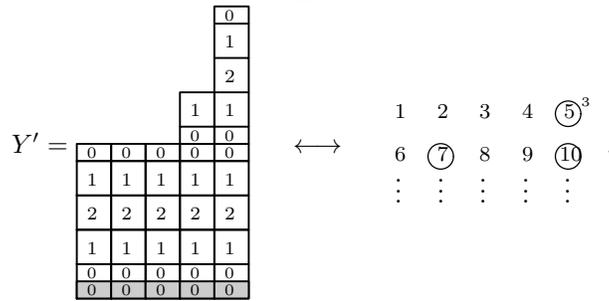
\begin{center}
$Y'=$\raisebox{-0.5\height}{
\begin{texdraw}
\drawdim em \setunitscale 0.13 \linewd 0.5



\move(10 0)\lvec(20 0)\lvec(20 10)\lvec(10 10)\lvec(10 0)\htext(13
1){\tiny $0$}

\move(20 0)\lvec(30 0)\lvec(30 10)\lvec(20 10)\lvec(20 0)\htext(23
1){\tiny $0$}

\move(30 0)\lvec(40 0)\lvec(40 10)\lvec(30 10)\lvec(30 0)\htext(33
1){\tiny $0$}

\move(40 0)\lvec(50 0)\lvec(50 10)\lvec(40 10)\lvec(40 0)\htext(43
1){\tiny $0$}

\move(50 0)\lvec(60 0)\lvec(60 10)\lvec(50 10)\lvec(50 0)\htext(53
1){\tiny $0$}
%


\move(10 0)\lvec(20 0)\lvec(20 5)\lvec(10 5)\lvec(10 0)\lfill
f:0.8 \htext(13 6){\tiny $0$}

\move(20 0)\lvec(30 0)\lvec(30 5)\lvec(20 5)\lvec(20 0)\lfill
f:0.8 \htext(23 6){\tiny $0$}

\move(30 0)\lvec(40 0)\lvec(40 5)\lvec(30 5)\lvec(30 0)\lfill
f:0.8 \htext(33 6){\tiny $0$}

\move(40 0)\lvec(50 0)\lvec(50 5)\lvec(40 5)\lvec(40 0)\lfill
f:0.8 \htext(43 6){\tiny $0$}

\move(50 0)\lvec(60 0)\lvec(60 5)\lvec(50 5)\lvec(50 0)\lfill
f:0.8 \htext(53 6){\tiny $0$}


\move(10 10)\lvec(20 10)\lvec(20 20)\lvec(10 20)\lvec(10
10)\htext(13 13){$_1$}

\move(20 10)\lvec(30 10)\lvec(30 20)\lvec(20 20)\lvec(20
10)\htext(23 13){$_1$}

\move(30 10)\lvec(40 10)\lvec(40 20)\lvec(30 20)\lvec(30
10)\htext(33 13){$_1$}

\move(40 10)\lvec(50 10)\lvec(50 20)\lvec(40 20)\lvec(40
10)\htext(43 13){$_1$}

\move(50 10)\lvec(60 10)\lvec(60 20)\lvec(50 20)\lvec(50
10)\htext(53 13){$_1$}
%

\move(10 20)\lvec(20 20)\lvec(20 30)\lvec(10 30)\lvec(10
20)\htext(13 23){$_2$}

\move(20 20)\lvec(30 20)\lvec(30 30)\lvec(20 30)\lvec(20
20)\htext(23 23){$_2$}

\move(30 20)\lvec(40 20)\lvec(40 30)\lvec(30 30)\lvec(30
20)\htext(33 23){$_2$}

\move(40 20)\lvec(50 20)\lvec(50 30)\lvec(40 30)\lvec(40
20)\htext(43 23){$_2$}

\move(50 20)\lvec(60 20)\lvec(60 30)\lvec(50 30)\lvec(50
20)\htext(53 23){$_2$}
%


\move(10 30)\lvec(20 30)\lvec(20 40)\lvec(10 40)\lvec(10
30)\htext(13 33){$_1$}

\move(20 30)\lvec(30 30)\lvec(30 40)\lvec(20 40)\lvec(20
30)\htext(23 33){$_1$}

\move(30 30)\lvec(40 30)\lvec(40 40)\lvec(30 40)\lvec(30
30)\htext(33 33){$_1$}

\move(40 30)\lvec(50 30)\lvec(50 40)\lvec(40 40)\lvec(40
30)\htext(43 33){$_1$}

\move(50 30)\lvec(60 30)\lvec(60 40)\lvec(50 40)\lvec(50
30)\htext(53 33){$_1$}
\move(10 40)\lvec(20 40)\lvec(20 45)\lvec(10 45)\lvec(10
40)\htext(13 41){\tiny $0$}

\move(20 40)\lvec(30 40)\lvec(30 45)\lvec(20 45)\lvec(20
40)\htext(23 41){\tiny $0$}

\move(30 40)\lvec(40 40)\lvec(40 45)\lvec(30 45)\lvec(30
40)\htext(33 41){\tiny $0$}

\move(40 40)\lvec(50 40)\lvec(50 45)\lvec(40 45)\lvec(40
40)\htext(43 41){\tiny $0$}

\move(50 40)\lvec(60 40)\lvec(60 45)\lvec(50 45)\lvec(50
40)\htext(53 41){\tiny $0$}

\move(50 45)\lvec(60 45)\lvec(60 50)\lvec(50 50)\lvec(50
45)\htext(53 46){\tiny $0$}

\move(50 50)\lvec(50 60)\lvec(60 60)\lvec(60 50)\htext(53
53){$_1$}

\move(50 60)\lvec(50 70)\lvec(60 70)\lvec(60 60)\htext(53
63){$_2$}

\move(50 70)\lvec(50 80)\lvec(60 80)\lvec(60 70)\htext(53
73){$_1$}

\move(50 80)\lvec(50 85)\lvec(60 85)\lvec(60 80)\htext(53
81){\tiny $0$}

\move(40 45)\lvec(50 45)\lvec(50 50)\lvec(40 50)\lvec(40
45)\htext(43 46){\tiny $0$}

\move(40 50)\lvec(40 60)\lvec(50 60)\lvec(50 50)\htext(43
53){$_1$}
\end{texdraw}}\ \ \ \ \ $\longleftrightarrow$ \ \ \ \ \
\raisebox{-0.5\height}{\begin{texdraw} \drawdim em \setunitscale
0.16 \linewd 0.5 \fontsize{8}{8}

\htext(10 0){$1$}\htext(20 0){$2$}\htext(30 0){$3$}\htext(40
0){$4$}\htext(50 0){$5$}

\htext(10 -10){$6$}\htext(20 -10){$7$}\htext(30 -10){$8$}\htext(40
-10){$9$}\htext(49 -10){$10$}

\htext(10 -20){$\vdots$}\htext(20 -20){$\vdots$}\htext(30
-20){$\vdots$}\htext(40 -20){$\vdots$}\htext(50 -20){$\vdots$}

\move(51 1)\linewd .3 \lcir r:3 \move(21 -9)\linewd .3 \lcir r:3
\fontsize{4}{4}\move(54 3)\htext{$3$} \move(51 -9)\linewd .3 \lcir
r:3
\end{texdraw}}\ \ \ .
\end{center}

}
\end{ex}

Let $\widetilde{Y}$ be the proper Young wall which is obtained
from $Y$ by applying $(B_i)$ ($i=1,2,3,4$) until there is no bead
movable up or removable. Note that $\widetilde{Y}$ does not depend
on the order of steps, hence is uniquely determined.

\begin{lem}\label{wtA2n2}
Let $Y$ be a proper Young wall in $\Z(\Lambda)$. Then
$Y\in\mathcal{Z}(\Lambda)_{\Lambda-m\delta}$ for some $m\geq 0$ if
and only if $\widetilde{Y}=Y_{\Lambda}$.
\end{lem}
\pf Suppose that $R_k$ is of type I and let $r_k$ be the number of
beads in $R_k$ occurring in the bead configuration of $Y$. Note
that $\widetilde{Y}=Y_{\Lambda}$ if and only if $r_k=r_{L-k}$ for
$1\leq k\leq n$.

Suppose that ${\rm wt}(Y)=Y-m\delta$ for some $m\geq 0$. We assume
that $\frak{g}$ is of type $A_{2n}^{(2)}$ (the proof for
$D_{n+1}^{(2)}$ is similar).

By considering the content of the blocks corresponding to each
bead, it is straightforward to check that
\begin{equation}
{\rm cont}(Y)=\sum_{i=0}^n c_i\alpha_i + M\delta,
\end{equation}
for some $M\geq 0$, where
\begin{equation}
c_i=
\begin{cases}
\sum_{j=1}^{2n}r_j & \text{if $i=0$}, \\
\sum_{j=i+1}^{2n-i}r_j + \sum_{j=2n-i+1}^{2n}2r_j & \text{if
$1\leq i \leq n-1$}, \\
\sum_{j=n+1}^{2n}r_j & \text{if $i=n$}.
\end{cases}
\end{equation}

 Since ${\rm cont}(Y)=m\delta$ and $\delta=\sum_{i=0}^{n-1}2\alpha_i+\alpha_n$, it
follows that $c_0=c_1=\cdots=c_{n-1}=2c_n$. From the equations
$c_{i-1}=c_{i}$ for $1 \leq i\leq n-1$ and $c_{n-1}=2c_n$, we
obtain $r_i=r_{2n-i+1}$ and $r_{n}=r_{n+1}$, respectively.

Conversely, it is clear by Lemma \ref{abacusA2n2} that
$\widetilde{Y}=Y_{\Lambda}$ implies that ${\rm wt}(Y)=Y-m\delta$
for some $m\geq 0$. \qed\vskip 3mm

Fix $m\geq 0$. Let $Y$ be a proper Young wall in
$\Z(\Lambda)_{\Lambda-m\delta}$. Consider its bead configuration.
Let $R_k$ be a runner of type I and $r_k$ the number of beads in
$R_k$. Let $b_{i}^{(k)}$ ($1\leq i\leq r_k$) be the beads in $R_k$
enumerated from the bottom to top. Suppose that $b_{i}^{(k)}$
($1\leq i\leq r_k$) is located at $N^{(k)}_i$ and put
$p_i^{(k)}=\frac{N^{(k)}_i-k}{L}$. Then $p_i^{(k)}\geq r_k-i$ and
the sequence $p_i^{(k)}$ ($1\leq i\leq r_k$) forms a strict
partition, say $\mu^{(k)}$. By Lemma \ref{wtA2n2}, $r_k=r_{L-k}$
for $1\leq k\leq n$ and the pair of $\mu^{(k)}$ and $\mu^{(L-k)}$
determines a unique partition
$\lambda^{(k)}=(\mu^{(k)}|\mu^{(L-k)})$ whose sum is the number of
all possible moving and removing steps in $R_k$ and $R_{L-k}$ to
obtain $\widetilde{Y}$.

In $R_{\ell}$, suppose that there are $m_k$ beads at $k\ell$
($k\geq 1$). Set $\lambda^{(0)}=(1^{m_1},2^{m_2},\cdots)$. We
define
\begin{equation}
\pi(Y)=(\lambda^{(0)},\cdots,\lambda^{(n)}).
\end{equation}
Note that
$|\lambda^{(0)}|+\epsilon\sum_{i=1}^{n}|\lambda^{(i)}|=m$.
Conversely, for a given ($n+1$)-tuple of partitions
$(\lambda^{(0)},\cdots,\lambda^{(n)})$ with
$|\lambda^{(0)}|+\epsilon\sum_{i=1}^{n}|\lambda^{(i)}|=m$, we can
associate a unique $Y\in\mathcal{Z}(\Lambda)_{\Lambda-m\delta}$ by
reversing the construction of $\pi$. Then
$\pi(Y)=(\lambda^{(0)},\cdots,\lambda^{(n)})$ and it follows that
$\pi$ is a bijection.

Summarizing the above argument, we obtain

\begin{thm}\label{ZA2n2}
For $m\geq 0$, the map
\begin{equation}
\pi : \Z(\Lambda)_{\Lambda-m\delta} \longrightarrow
\bigsqcup_{m_0+\epsilon m_1=m}\cP(m_0)\times \cP^{(n)}(m_{1})
\end{equation}
is a bijection. \qed
\end{thm}

\begin{ex}\label{exCh4}{\rm Suppose that $\frak{g}=A_4^{(2)}$. Let $Y$ be a proper Young wall in $\Z(\Lambda_0)$
whose bead configuration is as follows:\vskip 3mm

\begin{center}
\raisebox{-0.5\height}{\begin{texdraw} \drawdim em \setunitscale
0.16 \linewd 0.5 \fontsize{8}{8}

\htext(10 0){$1$}\htext(20 0){$2$}\htext(30 0){$3$}\htext(40
0){$4$}\htext(50 0){$5$}

\htext(10 -10){$6$}\htext(20 -10){$7$}\htext(30 -10){$8$}\htext(40
-10){$9$}\htext(49 -10){$10$}

\htext(10 -20){$11$}\htext(20 -20){$12$}\htext(30
-20){$13$}\htext(40 -20){$14$}\htext(49 -20){$15$}

\htext(10 -30){$16$}\htext(20 -30){$17$}\htext(30
-30){$18$}\htext(40 -30){$19$}\htext(49 -30){$20$}

\htext(10 -40){$\vdots$}\htext(20 -40){$\vdots$}\htext(30
-40){$\vdots$}\htext(40 -40){$\vdots$}\htext(50 -40){$\vdots$}

\move(12 -9)\linewd .3 \lcir r:3.5

\move(13 -29)\linewd .3 \lcir r:3.5

\move(22 1)\linewd .3 \lcir r:3.5

\move(23 -19)\linewd .3 \lcir r:3.5

\move(33 -19)\linewd .3 \lcir r:3.5

\move(33 -29)\linewd .3 \lcir r:3.5

\move(42 1)\linewd .3 \lcir r:3.5

\move(43 -29)\linewd .3 \lcir r:3.5

\move(52 -9)\linewd .3 \lcir r:3.5

\move(52 -19)\linewd .3 \lcir r:3.5

\fontsize{4}{4}\move(54 -6)\htext{$4$} \move(54 -16)\htext{$2$}
\end{texdraw}}
\end{center}
Then we have $\widetilde{Y}=Y_{\Lambda_0}$. Hence,
$Y\in\Z(\Lambda_0)_{\Lambda_0-32\delta}$ and
$\pi(Y)=(\lambda^{(0)},\lambda^{(1)},\lambda^{(2)})$, where
\begin{equation}
\begin{split}
\lambda^{(0)}&=(2^4,3^2), \\
\lambda^{(1)}&=((1,3)|(0,3))=(1,2^2,4), \\
\lambda^{(2)}&=((0,2)|(2,3))=(1,4^2).
\end{split}
\end{equation}

}
\end{ex}

Now, we recover the formula for the string functions in
\cite{Kac90}.
\begin{cor}{\rm }
We have
$\Sigma^{\Lambda}_{\Lambda}(q)=\dfrac{1}{(q)_{\infty}(q^{\epsilon})^n_{\infty}}$.
\qed
\end{cor}

\begin{rem}{\rm If $\frak{g}$ is of type $A_{2n}^{(2)}$, then
$\mathcal{Y}(\Lambda)$ can be identified with the set of
partitions satisfying the conditions:
\begin{itemize}
\item[(1)] only parts divisible by $\ell$ may be repeated,

\item[(2)] the smallest part is smaller than $\ell$,

\item[(3)] the difference between successive parts is at most
$\ell$ and strictly less than $\ell$ if either part is divisible
by $\ell$.
\end{itemize}
In \cite{Y}, Yamada also gave another combinatorial description of
weight vectors for the basic representation $V(\Lambda)$. By using
vertex operator construction, he showed that the weight vectors
for $V(\Lambda)$ can be parametrized by the set of strict
partitions whose parts are not divisible by $\ell$, say
$\cP'_{\ell}$. Then he described the bead configurations of
elements in $\cP'_{\ell}$ and computed the weight multiplicities
of $V(\Lambda)$ in a similar way.

On the other hand, in  \cite{Be}, Bessenrodt constructed a certain
bijection between two kinds of partition sets generalizing the
Andrews-Olsson partition identity. As a particular case of her
result, we can establish an explicit weight-preserving bijection
between $\mathcal{Y}(\Lambda)$ and $\cP'_{\ell}$. But, unlike
$\mathcal{Y}(\Lambda)$, it seems to be difficult to describe a
crystal graph structure on $\cP'_{\ell}$.
 }
\end{rem}

\section{$A_{2n-1}^{(2)}$, $D_{n+1}^{(1)}$-case}

Suppose that $\frak{g}$ is of type $A_{2n-1}^{(2)}$  or
$D_{n+1}^{(1)}$ ($n\geq 3$), and $\Lambda$ is a dominant integral
weight of level $1$. The patterns for $\Z(\Lambda)$ are given as
follows:\vskip 5mm

$A_{2n-1}^{(2)}$ ($n\geq 3$),\vskip 5mm

\begin{center}
$\Lambda=\Lambda_0$ : \raisebox{-1\height}{\begin{texdraw}
\textref h:C v:C \fontsize{6}{6}\selectfont \drawdim mm
\setunitscale 5
\newcommand{\dtri}{\bsegment \move(-1 0)\lvec(0 1)\lvec(0 0)\lvec(-1
0)\ifill f:0.7 \esegment } \move(0 0)\dtri \move(-1 0)\dtri
\move(-2 0)\dtri \move(-3 0)\dtri \move(0 0)\rlvec(-4.3 0) \move(0
1)\rlvec(-4.3 0) \move(0 2)\rlvec(-4.3 0) \move(0 3.5)\rlvec(-4.3
0) \move(0 4.5)\rlvec(-4.3 0) \move(0 6)\rlvec(-4.3 0) \move(0
7)\rlvec(-4.3 0) \move(0 8)\rlvec(-4.3 0) \move(0 9)\rlvec(-4.3 0)
\move(0 0)\rlvec(0 9.3) \move(-1 0)\rlvec(0 9.3) \move(-2
0)\rlvec(0 9.3) \move(-3 0)\rlvec(0 9.3) \move(-4 0)\rlvec(0 9.3)
\move(-1 0)\rlvec(1 1) \move(-2 0)\rlvec(1 1) \move(-3 0)\rlvec(1
1) \move(-4 0)\rlvec(1 1) \move(-1 7)\rlvec(1 1) \move(-2
7)\rlvec(1 1) \move(-3 7)\rlvec(1 1) \move(-4 7)\rlvec(1 1)
\vtext(-0.5 2.75){$\cdots$} \vtext(-0.5 5.25){$\cdots$}
\vtext(-1.5 2.75){$\cdots$} \vtext(-1.5 5.25){$\cdots$}
\vtext(-2.5 2.75){$\cdots$} \vtext(-2.5 5.25){$\cdots$}
\vtext(-3.5 2.75){$\cdots$} \vtext(-3.5 5.25){$\cdots$}
\htext(-0.25 7.27){$1$} \htext(-0.75 7.75){$0$} \htext(-1.25
7.27){$0$} \htext(-1.75 7.75){$1$} \htext(-2.25 7.27){$1$}
\htext(-2.75 7.75){$0$} \htext(-3.25 7.27){$0$} \htext(-3.75
7.75){$1$} \htext(-0.25 0.27){$1$} \htext(-0.75 0.75){$0$}
\htext(-1.25 0.27){$0$} \htext(-1.75 0.75){$1$} \htext(-2.25
0.27){$1$} \htext(-2.75 0.75){$0$} \htext(-3.25 0.27){$0$}
\htext(-3.75 0.75){$1$} \htext(-0.5 1.5){$2$} \htext(-1.5
1.5){$2$} \htext(-2.5 1.5){$2$} \htext(-3.5 1.5){$2$} \htext(-0.5
6.5){$2$} \htext(-1.5 6.5){$2$} \htext(-2.5 6.5){$2$} \htext(-3.5
6.5){$2$} \htext(-0.5 8.5){$2$} \htext(-1.5 8.5){$2$} \htext(-2.5
8.5){$2$} \htext(-3.5 8.5){$2$} \htext(-0.5 4){$n$} \htext(-1.5
4){$n$} \htext(-2.5 4){$n$} \htext(-3.5 4){$n$}
\end{texdraw}}\hskip 2cm
$\Lambda=\Lambda_1$ : \raisebox{-1\height}{\begin{texdraw}
\textref h:C v:C \fontsize{6}{6}\selectfont \drawdim mm
\setunitscale 5
\newcommand{\dtri}{ \bsegment \move(-1 0)\lvec(0 1)\lvec(0 0)\lvec(-1
0)\ifill f:0.7 \esegment } \move(0 0)\dtri \move(-1 0)\dtri
\move(-2 0)\dtri \move(-3 0)\dtri \move(0 0)\rlvec(-4.3 0) \move(0
1)\rlvec(-4.3 0) \move(0 2)\rlvec(-4.3 0) \move(0 3.5)\rlvec(-4.3
0) \move(0 4.5)\rlvec(-4.3 0) \move(0 6)\rlvec(-4.3 0) \move(0
7)\rlvec(-4.3 0) \move(0 8)\rlvec(-4.3 0) \move(0 9)\rlvec(-4.3 0)
\move(0 0)\rlvec(0 9.3) \move(-1 0)\rlvec(0 9.3) \move(-2
0)\rlvec(0 9.3) \move(-3 0)\rlvec(0 9.3) \move(-4 0)\rlvec(0 9.3)
\move(-1 0)\rlvec(1 1) \move(-2 0)\rlvec(1 1) \move(-3 0)\rlvec(1
1) \move(-4 0)\rlvec(1 1) \move(-1 7)\rlvec(1 1) \move(-2
7)\rlvec(1 1) \move(-3 7)\rlvec(1 1) \move(-4 7)\rlvec(1 1)
\vtext(-0.5 2.75){$\cdots$} \vtext(-0.5 5.25){$\cdots$}
\vtext(-1.5 2.75){$\cdots$} \vtext(-1.5 5.25){$\cdots$}
\vtext(-2.5 2.75){$\cdots$} \vtext(-2.5 5.25){$\cdots$}
\vtext(-3.5 2.75){$\cdots$} \vtext(-3.5 5.25){$\cdots$}
\htext(-0.25 7.27){$0$} \htext(-0.75 7.75){$1$} \htext(-1.25
7.27){$1$} \htext(-1.75 7.75){$0$} \htext(-2.25 7.27){$0$}
\htext(-2.75 7.75){$1$} \htext(-3.25 7.27){$1$} \htext(-3.75
7.75){$0$} \htext(-0.25 0.27){$0$} \htext(-0.75 0.75){$1$}
\htext(-1.25 0.27){$1$} \htext(-1.75 0.75){$0$} \htext(-2.25
0.27){$0$} \htext(-2.75 0.75){$1$} \htext(-3.25 0.27){$1$}
\htext(-3.75 0.75){$0$} \htext(-0.5 1.5){$2$} \htext(-1.5
1.5){$2$} \htext(-2.5 1.5){$2$} \htext(-3.5 1.5){$2$} \htext(-0.5
6.5){$2$} \htext(-1.5 6.5){$2$} \htext(-2.5 6.5){$2$} \htext(-3.5
6.5){$2$} \htext(-0.5 8.5){$2$} \htext(-1.5 8.5){$2$} \htext(-2.5
8.5){$2$} \htext(-3.5 8.5){$2$} \htext(-0.5 4){$n$} \htext(-1.5
4){$n$} \htext(-2.5 4){$n$} \htext(-3.5 4){$n$}
\end{texdraw}}
\end{center}\vskip 5mm

$D_{n+1}^{(1)}$ ($n\geq 3$),
\begin{center}
$\Lambda=\Lambda_0$ : \raisebox{-1\height}{\begin{texdraw}\textref
h:C v:C \fontsize{6}{6}\selectfont \drawdim mm \setunitscale 5
\newcommand{\dtri}{ \bsegment \move(-1 0)\lvec(0 1)\lvec(0 0)\lvec(-1
0)\ifill f:0.7 \esegment } \move(0 0)\dtri \move(-1 0)\dtri
\move(-2 0)\dtri \move(-3 0)\dtri \move(0 0)\rlvec(-4.3 0) \move(0
1)\rlvec(-4.3 0) \move(0 2)\rlvec(-4.3 0) \move(0 3.5)\rlvec(-4.3
0) \move(0 4.5)\rlvec(-4.3 0) \move(0 5.5)\rlvec(-4.3 0) \move(0
6.5)\rlvec(-4.3 0) \move(0 8)\rlvec(-4.3 0) \move(0 9)\rlvec(-4.3
0) \move(0 10)\rlvec(-4.3 0) \move(0 11)\rlvec(-4.3 0) \move(0
0)\rlvec(0 11.3) \move(-1 0)\rlvec(0 11.3) \move(-2 0)\rlvec(0
11.3) \move(-3 0)\rlvec(0 11.3) \move(-4 0)\rlvec(0 11.3) \move(-1
0)\rlvec(1 1) \move(-2 0)\rlvec(1 1) \move(-3 0)\rlvec(1 1)
\move(-4 0)\rlvec(1 1) \move(-1 9)\rlvec(1 1) \move(-2 9)\rlvec(1
1) \move(-3 9)\rlvec(1 1) \move(-4 9)\rlvec(1 1) \htext(-0.3
0.25){$1$} \htext(-0.75 0.75){$0$} \htext(-0.5 1.5){$2$}
\vtext(-0.5 2.75){$\cdots$} \htext(-0.5 4){$n\!\!-\!\!1$}
\htext(-0.5 6){$n\!\!-\!\!1$} \htext(-0.5 8.5){$2$} \htext(-0.3
9.25){$1$} \htext(-0.75 9.75){$0$} \htext(-0.5 10.5){$2$}
\htext(-2.3 0.25){$1$} \htext(-2.75 0.75){$0$} \htext(-2.5
1.5){$2$} \vtext(-2.5 2.75){$\cdots$} \htext(-2.5
4){$n\!\!-\!\!1$} \htext(-2.5 6){$n\!\!-\!\!1$} \htext(-2.5
8.5){$2$} \htext(-2.3 9.25){$1$} \htext(-2.75 9.75){$0$}
\htext(-2.5 10.5){$2$} \htext(-1.3 0.25){$0$} \htext(-1.75
0.75){$1$} \htext(-1.5 1.5){$2$} \vtext(-1.5 2.75){$\cdots$}
\htext(-1.5 4){$n\!\!-\!\!1$} \htext(-1.5 6){$n\!\!-\!\!1$}
\htext(-1.5 8.5){$2$} \htext(-1.3 9.25){$0$} \htext(-1.75
9.75){$1$} \htext(-1.5 10.5){$2$} \htext(-3.3 0.25){$0$}
\htext(-3.75 0.75){$1$} \htext(-3.5 1.5){$2$} \vtext(-3.5
2.75){$\cdots$} \htext(-3.5 4){$n\!\!-\!\!1$} \htext(-3.5
6){$n\!\!-\!\!1$} \htext(-3.5 8.5){$2$} \htext(-3.3 9.25){$0$}
\htext(-3.75 9.75){$1$} \htext(-3.5 10.5){$2$} \htext(-0.4
4.75){$n\!\!+\!\!1$} \htext(-2.4 4.75){$n\!\!+\!\!1$} \htext(-1.5
5.25){$n\!\!+\!\!1$} \htext(-3.5 5.25){$n\!\!+\!\!1$} \move(-0.5
5)\rlvec(0.5 0.5) \move(-2.5 4)\rlvec(0.05 0.05)\rmove(0.45
0.45)\rlvec(0.5 0.5) \move(-2.5 5)\rlvec(0.5 0.5)\rmove(0.45
0.45)\rlvec(0.05 0.05) \move(-4.5 4)\rlvec(0.05 0.05)\rmove(0.45
0.45)\rlvec(0.5 0.5) \htext(-1.4 4.75){$n$} \htext(-3.4 4.75){$n$}
\htext(-0.6 5.25){$n$} \htext(-2.6 5.25){$n$} \vtext(-0.5
7.25){$\cdots$} \vtext(-1.5 7.25){$\cdots$} \vtext(-2.5
7.25){$\cdots$} \vtext(-3.5 7.25){$\cdots$}
\end{texdraw}}\hskip 2cm
$\Lambda=\Lambda_1$ : \raisebox{-1\height}{
\begin{texdraw} \textref h:C v:C \fontsize{6}{6}\selectfont \drawdim mm
\setunitscale 5 \newcommand{\dtri}{ \bsegment \move(-1 0)\lvec(0
1)\lvec(0 0)\lvec(-1 0)\ifill f:0.7 \esegment } \move(0 0)\dtri
\move(-1 0)\dtri \move(-2 0)\dtri \move(-3 0)\dtri \move(0
0)\rlvec(-4.3 0) \move(0 1)\rlvec(-4.3 0) \move(0 2)\rlvec(-4.3 0)
\move(0 3.5)\rlvec(-4.3 0) \move(0 4.5)\rlvec(-4.3 0) \move(0
5.5)\rlvec(-4.3 0) \move(0 6.5)\rlvec(-4.3 0) \move(0
8)\rlvec(-4.3 0) \move(0 9)\rlvec(-4.3 0) \move(0 10)\rlvec(-4.3
0) \move(0 11)\rlvec(-4.3 0) \move(0 0)\rlvec(0 11.3) \move(-1
0)\rlvec(0 11.3) \move(-2 0)\rlvec(0 11.3) \move(-3 0)\rlvec(0
11.3) \move(-4 0)\rlvec(0 11.3) \move(-1 0)\rlvec(1 1) \move(-2
0)\rlvec(1 1) \move(-3 0)\rlvec(1 1) \move(-4 0)\rlvec(1 1)
\move(-1 9)\rlvec(1 1) \move(-2 9)\rlvec(1 1) \move(-3 9)\rlvec(1
1) \move(-4 9)\rlvec(1 1) \htext(-0.3 0.25){$0$} \htext(-0.75
0.75){$1$} \htext(-0.5 1.5){$2$} \vtext(-0.5 2.75){$\cdots$}
\htext(-0.5 4){$n\!\!-\!\!1$} \htext(-0.5 6){$n\!\!-\!\!1$}
\htext(-0.5 8.5){$2$} \htext(-0.3 9.25){$0$} \htext(-0.75
9.75){$1$} \htext(-0.5 10.5){$2$} \htext(-2.3 0.25){$0$}
\htext(-2.75 0.75){$1$} \htext(-2.5 1.5){$2$} \vtext(-2.5
2.75){$\cdots$} \htext(-2.5 4){$n\!\!-\!\!1$} \htext(-2.5
6){$n\!\!-\!\!1$} \htext(-2.5 8.5){$2$} \htext(-2.3 9.25){$0$}
\htext(-2.75 9.75){$1$} \htext(-2.5 10.5){$2$} \htext(-1.3
0.25){$1$} \htext(-1.75 0.75){$0$} \htext(-1.5 1.5){$2$}
\vtext(-1.5 2.75){$\cdots$} \htext(-1.5 4){$n\!\!-\!\!1$}
\htext(-1.5 6){$n\!\!-\!\!1$} \htext(-1.5 8.5){$2$} \htext(-1.3
9.25){$1$} \htext(-1.75 9.75){$0$} \htext(-1.5 10.5){$2$}
\htext(-3.3 0.25){$1$} \htext(-3.75 0.75){$0$} \htext(-3.5
1.5){$2$} \vtext(-3.5 2.75){$\cdots$} \htext(-3.5
4){$n\!\!-\!\!1$} \htext(-3.5 6){$n\!\!-\!\!1$} \htext(-3.5
8.5){$2$} \htext(-3.3 9.25){$1$} \htext(-3.75 9.75){$0$}
\htext(-3.5 10.5){$2$} \htext(-0.4 4.75){$n\!\!+\!\!1$}
\htext(-2.4 4.75){$n\!\!+\!\!1$} \htext(-1.5 5.25){$n\!\!+\!\!1$}
\htext(-3.5 5.25){$n\!\!+\!\!1$} \move(-0.5 5)\rlvec(0.5 0.5)
\move(-2.5 4)\rlvec(0.05 0.05)\rmove(0.45 0.45)\rlvec(0.5 0.5)
\move(-2.5 5)\rlvec(0.5 0.5)\rmove(0.45 0.45)\rlvec(0.05 0.05)
\move(-4.5 4)\rlvec(0.05 0.05)\rmove(0.45 0.45)\rlvec(0.5 0.5)
\htext(-1.4 4.75){$n$} \htext(-3.4 4.75){$n$} \htext(-0.6
5.25){$n$} \htext(-2.6 5.25){$n$} \vtext(-0.5 7.25){$\cdots$}
\vtext(-1.5 7.25){$\cdots$} \vtext(-2.5 7.25){$\cdots$}
\vtext(-3.5 7.25){$\cdots$}
\end{texdraw}}
\end{center}
\begin{center}
$\Lambda=\Lambda_{n}$ : \raisebox{-1\height}{\begin{texdraw}
\textref h:C v:C \fontsize{6}{6}\selectfont \drawdim mm
\setunitscale 5 \newcommand{\dtri}{ \bsegment \move(-1 0)\lvec(0
1)\lvec(0 0)\lvec(-1 0)\ifill f:0.7 \esegment } \move(0 0)\dtri
\move(-1 0)\dtri \move(-2 0)\dtri \move(-3 0)\dtri \move(0
0)\rlvec(-4.3 0) \move(0 1)\rlvec(-4.3 0) \move(0 2)\rlvec(-4.3 0)
\move(0 3.5)\rlvec(-4.3 0) \move(0 4.5)\rlvec(-4.3 0) \move(0
5.5)\rlvec(-4.3 0) \move(0 6.5)\rlvec(-4.3 0) \move(0
8)\rlvec(-4.3 0) \move(0 9)\rlvec(-4.3 0) \move(0 10)\rlvec(-4.3
0) \move(0 11)\rlvec(-4.3 0) \move(0 0)\rlvec(0 11.3) \move(-1
0)\rlvec(0 11.3) \move(-2 0)\rlvec(0 11.3) \move(-3 0)\rlvec(0
11.3) \move(-4 0)\rlvec(0 11.3) \move(-1 4.5)\rlvec(1 1) \move(-2
4.5)\rlvec(1 1) \move(-3 4.5)\rlvec(1 1) \move(-4 4.5)\rlvec(1 1)
\htext(-0.4 0.25){$n\!\!+\!\!1$} \htext(-0.5 1.5){$n\!\!-\!\!1$}
\vtext(-0.5 2.75){$\cdots$} \htext(-0.5 4){$2$} \htext(-0.3
4.75){$1$} \htext(-0.75 5.25){$0$} \htext(-0.5 6){$2$} \vtext(-0.5
7.25){$\cdots$} \htext(-0.5 8.5){$n\!\!-\!\!1$} \htext(-0.4
9.25){$n\!\!+\!\!1$} \htext(-0.5 10.5){$n\!\!-\!\!1$} \htext(-2.4
0.25){$n\!\!+\!\!1$} \htext(-2.5 1.5){$n\!\!-\!\!1$} \vtext(-2.5
2.75){$\cdots$} \htext(-2.5 4){$2$} \htext(-2.3 4.75){$1$}
\htext(-2.75 5.25){$0$} \htext(-2.5 6){$2$} \vtext(-2.5
7.25){$\cdots$} \htext(-2.5 8.5){$n\!\!-\!\!1$} \htext(-2.4
9.25){$n\!\!+\!\!1$} \htext(-2.5 10.5){$n\!\!-\!\!1$} \htext(-1.55
0.75){$n\!\!+\!\!1$} \htext(-1.5 1.5){$n\!\!-\!\!1$} \vtext(-1.5
2.75){$\cdots$} \htext(-1.5 4){$2$} \htext(-1.3 4.75){$0$}
\htext(-1.75 5.25){$1$} \htext(-1.5 6){$2$} \vtext(-1.5
7.25){$\cdots$} \htext(-1.5 8.5){$n\!\!-\!\!1$} \htext(-1.55
9.75){$n\!\!+\!\!1$} \htext(-1.5 10.5){$n\!\!-\!\!1$} \htext(-3.55
0.75){$n\!\!+\!\!1$} \htext(-3.5 1.5){$n\!\!-\!\!1$} \vtext(-3.5
2.75){$\cdots$} \htext(-3.5 4){$2$} \htext(-3.3 4.75){$0$}
\htext(-3.75 5.25){$1$} \htext(-3.5 6){$2$} \vtext(-3.5
7.25){$\cdots$} \htext(-3.5 8.5){$n\!\!-\!\!1$} \htext(-3.55
9.75){$n\!\!+\!\!1$} \htext(-3.5 10.5){$n\!\!-\!\!1$} \move(-0.5
0.5)\rlvec(0.5 0.5)\rmove(-0.95 -0.95)\rlvec(0.05 0.05) \move(-2
0)\rlvec(0.5 0.5)\rmove(0.45 0.45)\rlvec(0.05 0.05) \move(-3
0)\rlvec(0.05 0.05)\rmove(0.45 0.45)\rlvec(0.5 0.5) \move(-4
0)\rlvec(0.5 0.5)\rmove(0.45 0.45)\rlvec(0.05 0.05) \move(-1
9)\rlvec(0.05 0.05)\rmove(0.45 0.45)\rlvec(0.5 0.5) \move(-2
9)\rlvec(0.5 0.5)\rmove(0.45 0.45)\rlvec(0.05 0.05) \move(-3
9)\rlvec(0.05 0.05)\rmove(0.45 0.45)\rlvec(0.5 0.5) \move(-4
9)\rlvec(0.5 0.5) \rmove(0.45 0.45)\rlvec(0.05 0.05)\htext(-3.4
9.25){$n$} \htext(-3.4 0.25){$n$} \htext(-2.6 9.75){$n$}
\htext(-2.6 0.75){$n$} \htext(-1.4 0.25){$n$} \htext(-1.4
9.25){$n$} \htext(-0.6 9.75){$n$} \htext(-0.6 0.75){$n$}
\end{texdraw}}\hskip 2cm
$\Lambda=\Lambda_{n+1}$ : \raisebox{-1\height}{\begin{texdraw}
\textref h:C v:C \fontsize{6}{6}\selectfont \drawdim mm
\setunitscale 5 \newcommand{\dtri}{ \bsegment \move(-1 0)\lvec(0
1)\lvec(0 0)\lvec(-1 0)\ifill f:0.7 \esegment } \move(0 0)\dtri
\move(-1 0)\dtri \move(-2 0)\dtri \move(-3 0)\dtri \move(0
0)\rlvec(-4.3 0) \move(0 1)\rlvec(-4.3 0) \move(0 2)\rlvec(-4.3 0)
\move(0 3.5)\rlvec(-4.3 0) \move(0 4.5)\rlvec(-4.3 0) \move(0
5.5)\rlvec(-4.3 0) \move(0 6.5)\rlvec(-4.3 0) \move(0
8)\rlvec(-4.3 0) \move(0 9)\rlvec(-4.3 0) \move(0 10)\rlvec(-4.3
0) \move(0 11)\rlvec(-4.3 0) \move(0 0)\rlvec(0 11.3) \move(-1
0)\rlvec(0 11.3) \move(-2 0)\rlvec(0 11.3) \move(-3 0)\rlvec(0
11.3) \move(-4 0)\rlvec(0 11.3) \move(-1 4.5)\rlvec(1 1) \move(-2
4.5)\rlvec(1 1) \move(-3 4.5)\rlvec(1 1) \move(-4 4.5)\rlvec(1 1)
\htext(-0.3 0.25){$n$} \htext(-0.5 1.5){$n\!\!-\!\!1$} \vtext(-0.5
2.75){$\cdots$} \htext(-0.5 4){$2$} \htext(-0.3 4.75){$1$}
\htext(-0.75 5.25){$0$} \htext(-0.5 6){$2$} \vtext(-0.5
7.25){$\cdots$} \htext(-0.5 8.5){$n\!\!-\!\!1$} \htext(-0.3
9.25){$n$} \htext(-0.5 10.5){$n\!\!-\!\!1$} \htext(-2.3 0.25){$n$}
\htext(-2.5 1.5){$n\!\!-\!\!1$} \vtext(-2.5 2.75){$\cdots$}
\htext(-2.5 4){$2$} \htext(-2.3 4.75){$1$} \htext(-2.75 5.25){$0$}
\htext(-2.5 6){$2$} \vtext(-2.5 7.25){$\cdots$} \htext(-2.5
8.5){$n\!\!-\!\!1$} \htext(-2.3 9.25){$n$} \htext(-2.5
10.5){$n\!\!-\!\!1$} \htext(-1.75 0.75){$n$} \htext(-1.5
1.5){$n\!\!-\!\!1$} \vtext(-1.5 2.75){$\cdots$} \htext(-1.5
4){$2$} \htext(-1.3 4.75){$0$} \htext(-1.75 5.25){$1$} \htext(-1.5
6){$2$} \vtext(-1.5 7.25){$\cdots$} \htext(-1.5
8.5){$n\!\!-\!\!1$} \htext(-1.75 9.75){$n$} \htext(-1.5
10.5){$n\!\!-\!\!1$} \htext(-3.75 0.75){$n$} \htext(-3.5
1.5){$n\!\!-\!\!1$} \vtext(-3.5 2.75){$\cdots$} \htext(-3.5
4){$2$} \htext(-3.3 4.75){$0$} \htext(-3.75 5.25){$1$} \htext(-3.5
6){$2$} \vtext(-3.5 7.25){$\cdots$} \htext(-3.5
8.5){$n\!\!-\!\!1$} \htext(-3.75 9.75){$n$} \htext(-3.5
10.5){$n\!\!-\!\!1$} \move(-1 0)\rlvec(0.5 0.5)\rmove(0.45
0.45)\rlvec(0.05 0.05) \move(-2 0)\rlvec(0.05 0.05)\rmove(0.45
0.45)\rlvec(0.5 0.5) \move(-3 0)\rlvec(0.5 0.5)\rmove(0.45
0.45)\rlvec(0.05 0.05) \move(-4 0)\rlvec(0.05 0.05)\rmove(0.45
0.45)\rlvec(0.5 0.5) \move(-1 9)\rlvec(0.5 0.5)\rmove(0.45
0.45)\rlvec(0.05 0.05) \move(-2 9)\rlvec(0.05 0.05)\rmove(0.45
0.45)\rlvec(0.5 0.5) \move(-3 9)\rlvec(0.5 0.5)\rmove(0.45
0.45)\rlvec(0.05 0.05) \move(-4 9)\rlvec(0.05 0.05)\rmove(0.45
0.45)\rlvec(0.5 0.5) \htext(-3.5 9.25){$n\!\!+\!\!1$} \htext(-3.5
0.25){$n\!\!+\!\!1$} \htext(-2.5 9.75){$n\!\!+\!\!1$} \htext(-2.5
0.75){$n\!\!+\!\!1$} \htext(-1.5 0.25){$n\!\!+\!\!1$} \htext(-1.5
9.25){$n\!\!+\!\!1$} \htext(-0.5 9.75){$n\!\!+\!\!1$} \htext(-0.5
0.75){$n\!\!+\!\!1$}
\end{texdraw}}\hskip 2cm
\end{center}\vskip 5mm

Set
\begin{equation}
\begin{split}
&\ell=
\begin{cases}
2n-1 & \text{if $\frak{g}=A_{2n-1}^{(2)}$}, \\
n & \text{if $\frak{g}=D_{n+1}^{(1)}$},
\end{cases}
\hskip 3mm L=
\begin{cases}
\ell & \text{if $\frak{g}=A_{2n-1}^{(2)}$}, \\
2\ell & \text{if $\frak{g}=D_{n+1}^{(1)}$},
\end{cases}\\
&\epsilon=L/\ell.
\end{split}
\end{equation}
Note that $L$ is the number of blocks in a $\delta$-column or the
Coxeter number.
 Then we define the {\it abacus of type
$A_{2n-1}^{(2)}$ {\rm (}resp. $D_{n+1}^{(1)}${\rm )}} to be the
arrangement of positive integers in the following way:\vskip 5mm

$A_{2n-1}^{(2)}$
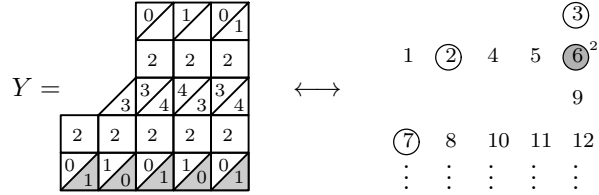
\begin{center}
\begin{texdraw}
\drawdim em \setunitscale 0.13 \linewd 0.5 \fontsize{10}{10}

\htext(0 0){$1$}\htext(20 0){$2$}\htext(40 0){$\cdots$}\htext(60
0){$\ell-1$}\htext(90 0){$\ell$}

\htext(-7 -15){$\ell+1$}\htext(13 -15){$\ell+2$}\htext(40
-15){$\cdots$}\htext(57 -15){$2\ell-1$}\htext(88 -15){$2\ell$}

\htext(0 -30){$\vdots$}\htext(20 -30){$\vdots$}\htext(67
-30){$\vdots$}\htext(90 -30){$\vdots$}
\end{texdraw}
\end{center}

\vskip 5mm

$D_{n+1}^{(1)}$
\begin{center}
\begin{texdraw}
\drawdim em \setunitscale 0.13 \linewd 0.5 \fontsize{10}{10}

\htext(0 0){$1$}\htext(20 0){$\cdots$}\htext(40
0){$\ell-1$}\htext(60 0){$\ell+1$}\htext(80 0){$\cdots$}\htext(100
0){$2\ell-1$}

\htext(-10 -20){$2\ell+1$}\htext(20 -20){$\cdots$}\htext(37
-20){$3\ell-1$}\htext(58 -20){$3\ell+1$}\htext(80
-20){$\cdots$}\htext(100 -20){$4\ell-1$}

\htext(0 -35){$\vdots$}\htext(45 -35){$\vdots$}\htext(67
-35){$\vdots$}\htext(108 -35){$\vdots$}

\htext(130 10){$\ell$}\htext(128 0){$2\ell$}\htext(128
-10){$3\ell$}\htext(128 -20){$4\ell$}\htext(130 -35){$\vdots$}
\end{texdraw}
\end{center}\vskip 5mm

Let $R_k$ ($1\leq k<L$, $k\neq\ell$ ) be the set of all integers
$s \equiv k \pmod L$ and let $R_{\ell}$ be the set of all integers
$s \equiv 0 \pmod \ell$. For $1\leq k <L$ ($k\neq \ell$), we
assume that the $k$th runner $R_k$ is of type I (cf. Section 5).
In $R_{\ell}$, we may put two kinds of beads, white $\bigcirc$ and
gray \raisebox{-.25 \height}{\begin{texdraw} \drawdim em
\setunitscale 0.13 \linewd 0.5 \fcir f:0.5 r:4 \lcir r:4
\end{texdraw}}, and place more than one bead with the same color at each position.
We call $R_{\ell}$ a {\it runner of type ${\rm III}$}. But we will
not move any bead in $R_{\ell}$.

For $Y\in\Z(\Lambda)$, let $|Y|=(|y_k|)_{k\geq 1}$ be its
associated partition. Let $\{\,|y_1|,\cdots,|y_r|\,\}$ be the set
of all non-zero parts in $|Y|$. Then by definition of
$\Z(\Lambda)$, the numbers $|y_k|$'s ($1\leq k\leq r$) are
distinct except when $|y_k|\equiv 0\pmod \ell$. We define the {\it
bead configuration of $Y$} to be the set of $r$ beads
$b_1,\cdots,b_r$ placed in the above abacus, where $b_k$ is at
$|y_k|$ and the color $c(b_k)$ of $b_k$ is determined by
\begin{equation}\label{color of bead}
c(b_k)=
\begin{cases}
\text {white} & \text{if the block at the top of $y_k$ is
\raisebox{-0.2\height}{
\begin{texdraw}
\drawdim em \setunitscale 0.1 \linewd 0.5 \move(-10 0)\lvec(0
0)\lvec(0 10)\lvec(-10 0)
\end{texdraw}}}\ ,
\\
\text{gray} & \text{if the block at the top of $y_k$ is
\raisebox{-0.2\height}{
\begin{texdraw}
\drawdim em \setunitscale 0.1 \linewd 0.5 \move(-10 0)\lvec(-10
10)\lvec(0 10)\lvec(-10 0)
\end{texdraw}
}}, \\
\text {white} & \text{otherwise}.
\end{cases}
\end{equation}
Then it is easy to see that $Y$ is uniquely determined by its bead
configuration.

\begin{ex}\label{ex5.1}{\rm (1) Suppose that $\frak{g}=A_5^{(2)}$ and $\Lambda=\Lambda_0$.

\begin{center}
$Y=$\raisebox{-0.5\height}{\begin{texdraw}\fontsize{7}{7}\drawdim
mm \setunitscale 0.5 \linewd 0.5

\move(0 0)\lvec(10 0)\lvec(10 10)\lvec(0 10)\lvec(0 0)

\move(10 0)\lvec(20 0)\lvec(20 10)\lvec(10 10)\lvec(10 0)

\move(20 0)\lvec(30 0)\lvec(30 10)\lvec(20 10)\lvec(20 0)

\move(30 0)\lvec(40 0)\lvec(40 10)\lvec(30 10)\lvec(30 0)

\move(40 0)\lvec(50 0)\lvec(50 10)\lvec(40 10)\lvec(40 0)
\move(10 10)\lvec(20 10)\lvec(20 20)\lvec(10 20)\lvec(10
10)\htext(13 13){$2$}

\move(20 10)\lvec(30 10)\lvec(30 20)\lvec(20 20)\lvec(20
10)\htext(23 13){$2$}

\move(30 10)\lvec(40 10)\lvec(40 20)\lvec(30 20)\lvec(30
10)\htext(33 13){$2$}

\move(40 10)\lvec(50 10)\lvec(50 20)\lvec(40 20)\lvec(40
10)\htext(43 13){$2$}
\move(10 20)\lvec(20 20)\lvec(20 30)\lvec(10 30)\lvec(10
20)\htext(13 23){$3$}

\move(20 20)\lvec(30 20)\lvec(30 30)\lvec(20 30)\lvec(20
20)\htext(23 23){$3$}

\move(30 20)\lvec(40 20)\lvec(40 30)\lvec(30 30)\lvec(30
20)\htext(33 23){$3$}

\move(40 20)\lvec(50 20)\lvec(50 30)\lvec(40 30)\lvec(40
20)\htext(43 23){$3$}
\move(10 30)\lvec(20 30)\lvec(20 40)\lvec(10 40)\lvec(10
30)\htext(13 33){$2$}

\move(20 30)\lvec(30 30)\lvec(30 40)\lvec(20 40)\lvec(20
30)\htext(23 33){$2$}

\move(30 30)\lvec(40 30)\lvec(40 40)\lvec(30 40)\lvec(30
30)\htext(33 33){$2$}

\move(40 30)\lvec(50 30)\lvec(50 40)\lvec(40 40)\lvec(40
30)\htext(43 33){$2$}
%


%

\move(40 50)\lvec(50 50)\lvec(50 60)\lvec(40 60)\lvec(40
50)\htext(43 53){$2$}
\move(40 60)\lvec(50 60)\lvec(50 70)\lvec(40 70)\lvec(40
60)\htext(43 63){$3$}
\move(0 0)\lvec(10 10)\lvec(10 0)\lvec(0 0)\lfill f:0.8 \htext(1
5){$0$}\htext(6 1){$1$}

\move(10 0)\lvec(20 10)\lvec(20 0)\lvec(10 0)\lfill f:0.8
\htext(11 5){$1$}\htext(16 1){$0$}

\move(20 0)\lvec(30 10)\lvec(30 0)\lvec(20 0)\lfill f:0.8
\htext(22 5){$0$}\htext(26 2){$1$}

\move(30 0)\lvec(40 10)\lvec(40 0)\lvec(30 0)\lfill f:0.8
\htext(32 5){$1$}\htext(36 2){$0$}

\move(40 0)\lvec(50 10)\lvec(50 0)\lvec(40 0)\lfill f:0.8
\htext(42 5){$0$}\htext(46 2){$1$}
\move(20 40)\lvec(30 50)\lvec(20 50)\lvec(20 40)\htext(22 45){$0$}

\move(30 40)\lvec(40 50)\lvec(30 50)\lvec(30 40)\htext(32 45){$1$}

\move(40 40)\lvec(50 50)\lvec(50 40)\lvec(40 40)\lvec(40
50)\htext(42 45){$0$}\htext(46 42){$1$}
\end{texdraw}}\ \ \ \ \ $\longleftrightarrow$ \ \ \ \ \
\raisebox{-0.5\height}{\begin{texdraw} \drawdim em \setunitscale
0.16 \linewd 0.5 \fontsize{8}{8}

\htext(10 0){$1$}\htext(20 0){$2$}\htext(30 0){$3$}\htext(40
0){$4$}\htext(50 0){$5$}

\htext(10 -10){$6$}\htext(20 -10){$7$}\htext(30 -10){$8$}\htext(40
-10){$9$}\htext(49 -10){$10$}

\htext(10 -20){$\vdots$}\htext(20 -20){$\vdots$}\htext(30
-20){$\vdots$}\htext(40 -20){$\vdots$}\htext(50 -20){$\vdots$}

\move(11 1)\linewd .3 \lcir r:3 \move(41 1)\linewd .3 \lcir r:3
\move(51 1)\linewd .3 \fcir f:0.7 r:3 \lcir r:3 \move(31
-9)\linewd .3 \lcir r:3 \fontsize{4}{4}\move(54 3)\htext{$2$}
\end{texdraw}}
\end{center}\vskip 3mm

(2) Suppose that $\frak{g}=D_4^{(1)}$ and $\Lambda=\Lambda_0$.
\begin{center}
$Y=$\raisebox{-0.5\height}{\begin{texdraw}\fontsize{7}{7}\drawdim
mm \setunitscale 0.5 \linewd 0.5

\move(0 0)\lvec(10 0)\lvec(10 10)\lvec(0 10)\lvec(0 0)

\move(10 0)\lvec(20 0)\lvec(20 10)\lvec(10 10)\lvec(10 0)

\move(20 0)\lvec(30 0)\lvec(30 10)\lvec(20 10)\lvec(20 0)

\move(30 0)\lvec(40 0)\lvec(40 10)\lvec(30 10)\lvec(30 0)

\move(40 0)\lvec(50 0)\lvec(50 10)\lvec(40 10)\lvec(40 0)
\move(0 10)\lvec(10 10)\lvec(10 20)\lvec(0 20)\lvec(0 10)\htext(3
13){$2$}

\move(10 10)\lvec(20 10)\lvec(20 20)\lvec(10 20)\lvec(10
10)\htext(13 13){$2$}

\move(20 10)\lvec(30 10)\lvec(30 20)\lvec(20 20)\lvec(20
10)\htext(23 13){$2$}

\move(30 10)\lvec(40 10)\lvec(40 20)\lvec(30 20)\lvec(30
10)\htext(33 13){$2$}

\move(40 10)\lvec(50 10)\lvec(50 20)\lvec(40 20)\lvec(40
10)\htext(43 13){$2$}
\move(10 20)\lvec(20 30)\lvec(20 20)\lvec(10 20)\htext(16 21){$3$}

\move(20 20)\lvec(30 20)\lvec(30 30)\lvec(20 30)\lvec(20
20)\lvec(30 30)\htext(21 25){$3$}\htext(26 21){$4$}

\move(30 20)\lvec(40 20)\lvec(40 30)\lvec(30 30)\lvec(30
20)\lvec(40 30)\htext(36 21){$3$}\htext(31 25){$4$}

\move(40 20)\lvec(50 20)\lvec(50 30)\lvec(40 30)\lvec(40
20)\lvec(50 30)\htext(41 25){$3$}\htext(46 21){$4$}

\move(20 30)\lvec(30 30)\lvec(30 40)\lvec(20 40)\lvec(20
30)\htext(23 33){$2$}

\move(30 30)\lvec(40 30)\lvec(40 40)\lvec(30 40)\lvec(30
30)\htext(33 33){$2$}

\move(40 30)\lvec(50 30)\lvec(50 40)\lvec(40 40)\lvec(40
30)\htext(43 33){$2$}
\move(0 0)\lvec(10 10)\lvec(10 0)\lvec(0 0)\lfill f:0.8 \htext(1
5){$0$}\htext(6 1){$1$}

\move(10 0)\lvec(20 10)\lvec(20 0)\lvec(10 0)\lfill f:0.8
\htext(11 5){$1$}\htext(16 1){$0$}

\move(20 0)\lvec(30 10)\lvec(30 0)\lvec(20 0)\lfill f:0.8
\htext(22 5){$0$}\htext(26 2){$1$}

\move(30 0)\lvec(40 10)\lvec(40 0)\lvec(30 0)\lfill f:0.8
\htext(32 5){$1$}\htext(36 2){$0$}

\move(40 0)\lvec(50 10)\lvec(50 0)\lvec(40 0)\lfill f:0.8
\htext(42 5){$0$}\htext(46 2){$1$}
\move(20 40)\lvec(30 50)\lvec(20 50)\lvec(20 40)\htext(22 45){$0$}

\move(30 40)\lvec(40 50)\lvec(30 50)\lvec(30 40)\htext(32 45){$1$}

\move(40 40)\lvec(50 50)\lvec(50 40)\lvec(40 40)\lvec(40
50)\lvec(50 50)\htext(42 45){$0$}\htext(46 42){$1$}
\end{texdraw}}\ \ \ \ \ $\longleftrightarrow$ \ \ \ \ \
\raisebox{-0.5\height}{\begin{texdraw} \drawdim em \setunitscale
0.16 \linewd 0.5 \fontsize{8}{8}

\htext(10 0){$1$}\htext(20 0){$2$}\htext(30 0){$4$}\htext(40
0){$5$}\htext(50 10){$3$}\htext(50 0){$6$}

\htext(10 -20){$7$}\htext(20 -20){$8$}\htext(30
-20){$10$}\htext(40 -20){$11$}\htext(50 -10){$9$}\htext(50
-20){$12$}

\htext(10 -30){$\vdots$}\htext(20 -30){$\vdots$}\htext(30
-30){$\vdots$}\htext(40 -30){$\vdots$}\htext(50 -30){$\vdots$}

\move(21 1)\linewd .3 \lcir r:3 \move(51 11)\linewd .3 \lcir r:3
\move(51 1)\linewd .3 \fcir f:0.7 r:3 \lcir r:3 \move(11
-19)\linewd .3 \lcir r:3 \fontsize{4}{4}\move(54 3)\htext{$2$}
\end{texdraw}}
\end{center}\vskip 3mm}
\end{ex}

Let us describe the algorithm of moving and removing beads in the
abacus.
\begin{lem}\label{abacusA(2n-1)(2)}
Let $Y$ be a proper Young wall in $\Z(\Lambda)$ and let $Y'$ be
the proper Young wall obtained by applying one of the following
processes to the bead configuration of $Y$:
\begin{itemize}
\item[$(B_1)$] if $b$ is a bead at $s$ in a runner of type I and
there is no bead at $s-L$, then move $b$ one position up and
change the color of the beads at $k$ {\rm ($s-L< k < s$)} in
$R_{\ell}$,

\item[$(B_2)$] if $b$ and $b'$ are beads at $s$ and $L-s$ {\rm
($1\leq s\leq n-1$)}, respectively, then remove $b$ and $b'$
simultaneously. Also, if $s<\ell<L-s$, then change the color of
the beads at $\ell$.
\end{itemize}
\noindent Then we have $\wt(Y')=\wt(Y)+\delta$.
\end{lem}
\pf  Note that $Y'$ is obtained by removing some $L$ blocks from
$Y$, say $\{\,b_1,\cdots,b_{L}\,\}$. Let $i_{k}$ ($1\leq k\leq L$)
be the color of $b_k$. It follows directly from the pattern for
$\Z(\Lambda)$ that $\sum_{k=1}^{L}\alpha_{i_k}=\delta$ and
therefore $\wt(Y')=\wt(Y)+\delta$.  \qed

\begin{ex}{\rm \mbox{}

(1) Let $Y$ be the proper Young wall in Example \ref{ex5.1} (1).
Apply $(B_1)$ to \raisebox{-0.25 \height}{\begin{texdraw} \drawdim
em \setunitscale 0.16 \linewd 0.5 \fontsize{8}{8} \htext(0 0){$8$}
\move(1 1)\linewd .3 \lcir r:3
\end{texdraw}}. This means that we remove a $\delta$-column in the first
column of $Y$ and shift the blocks which are placed in the left of
the $\delta$-column, to the right as far as possible.

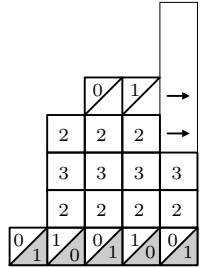
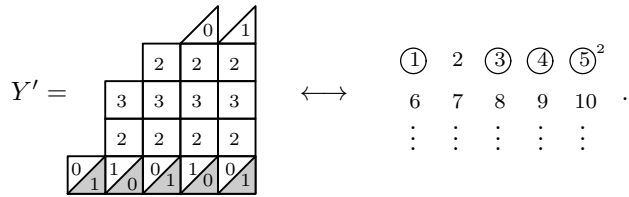
\begin{center}
\raisebox{-0.5\height}{\begin{texdraw}\fontsize{7}{7}\drawdim mm
\setunitscale 0.5 \linewd 0.5 \arrowheadtype t:F \arrowheadsize
l:2 w:2

\move(0 0)\lvec(10 0)\lvec(10 10)\lvec(0 10)\lvec(0 0)

\move(10 0)\lvec(20 0)\lvec(20 10)\lvec(10 10)\lvec(10 0)

\move(20 0)\lvec(30 0)\lvec(30 10)\lvec(20 10)\lvec(20 0)

\move(30 0)\lvec(40 0)\lvec(40 10)\lvec(30 10)\lvec(30 0)

\move(40 0)\lvec(50 0)\lvec(50 10)\lvec(40 10)\lvec(40 0)
\move(10 10)\lvec(20 10)\lvec(20 20)\lvec(10 20)\lvec(10
10)\htext(13 13){$2$}

\move(20 10)\lvec(30 10)\lvec(30 20)\lvec(20 20)\lvec(20
10)\htext(23 13){$2$}

\move(30 10)\lvec(40 10)\lvec(40 20)\lvec(30 20)\lvec(30
10)\htext(33 13){$2$}

\move(40 10)\lvec(50 10)\lvec(50 20)\lvec(40 20)\lvec(40
10)\htext(43 13){$2$}
\move(10 20)\lvec(20 20)\lvec(20 30)\lvec(10 30)\lvec(10
20)\htext(13 23){$3$}

\move(20 20)\lvec(30 20)\lvec(30 30)\lvec(20 30)\lvec(20
20)\htext(23 23){$3$}

\move(30 20)\lvec(40 20)\lvec(40 30)\lvec(30 30)\lvec(30
20)\htext(33 23){$3$}

\move(40 20)\lvec(50 20)\lvec(50 30)\lvec(40 30)\lvec(40
20)\htext(43 23){$3$}
\move(10 30)\lvec(20 30)\lvec(20 40)\lvec(10 40)\lvec(10
30)\htext(13 33){$2$}

\move(20 30)\lvec(30 30)\lvec(30 40)\lvec(20 40)\lvec(20
30)\htext(23 33){$2$}

\move(30 30)\lvec(40 30)\lvec(40 40)\lvec(30 40)\lvec(30
30)\htext(33 33){$2$}

\move(0 0)\lvec(10 10)\lvec(10 0)\lvec(0 0)\lfill f:0.8 \htext(1
5){$0$}\htext(6 1){$1$}

\move(10 0)\lvec(20 10)\lvec(20 0)\lvec(10 0)\lfill f:0.8
\htext(11 5){$1$}\htext(16 1){$0$}

\move(20 0)\lvec(30 10)\lvec(30 0)\lvec(20 0)\lfill f:0.8
\htext(22 5){$0$}\htext(26 2){$1$}

\move(30 0)\lvec(40 10)\lvec(40 0)\lvec(30 0)\lfill f:0.8
\htext(32 5){$1$}\htext(36 2){$0$}

\move(40 0)\lvec(50 10)\lvec(50 0)\lvec(40 0)\lfill f:0.8
\htext(42 5){$0$}\htext(46 2){$1$}
\move(20 40)\lvec(30 50)\lvec(20 50)\lvec(20 40)\htext(22 45){$0$}

\move(30 40)\lvec(40 50)\lvec(30 50)\lvec(30 40)\htext(32 45){$1$}

\move(42 35)\avec(48 35) \move(42 45)\avec(48 45)

\linewd 0.2 \move(40 20)\lvec(40 70)\lvec(50 70)\lvec(50 20)

\end{texdraw}}\end{center}\vskip 3mm

Therefore, we have
\begin{center}
$Y'=$\raisebox{-0.5\height}{\begin{texdraw}\fontsize{7}{7}\drawdim
mm \setunitscale 0.5 \linewd 0.5

\move(0 0)\lvec(10 0)\lvec(10 10)\lvec(0 10)\lvec(0 0)

\move(10 0)\lvec(20 0)\lvec(20 10)\lvec(10 10)\lvec(10 0)

\move(20 0)\lvec(30 0)\lvec(30 10)\lvec(20 10)\lvec(20 0)

\move(30 0)\lvec(40 0)\lvec(40 10)\lvec(30 10)\lvec(30 0)

\move(40 0)\lvec(50 0)\lvec(50 10)\lvec(40 10)\lvec(40 0)
\move(10 10)\lvec(20 10)\lvec(20 20)\lvec(10 20)\lvec(10
10)\htext(13 13){$2$}

\move(20 10)\lvec(30 10)\lvec(30 20)\lvec(20 20)\lvec(20
10)\htext(23 13){$2$}

\move(30 10)\lvec(40 10)\lvec(40 20)\lvec(30 20)\lvec(30
10)\htext(33 13){$2$}

\move(40 10)\lvec(50 10)\lvec(50 20)\lvec(40 20)\lvec(40
10)\htext(43 13){$2$}
\move(10 20)\lvec(20 20)\lvec(20 30)\lvec(10 30)\lvec(10
20)\htext(13 23){$3$}

\move(20 20)\lvec(30 20)\lvec(30 30)\lvec(20 30)\lvec(20
20)\htext(23 23){$3$}

\move(30 20)\lvec(40 20)\lvec(40 30)\lvec(30 30)\lvec(30
20)\htext(33 23){$3$}

\move(40 20)\lvec(50 20)\lvec(50 30)\lvec(40 30)\lvec(40
20)\htext(43 23){$3$}

\move(20 30)\lvec(30 30)\lvec(30 40)\lvec(20 40)\lvec(20
30)\htext(23 33){$2$}

\move(30 30)\lvec(40 30)\lvec(40 40)\lvec(30 40)\lvec(30
30)\htext(33 33){$2$}

\move(40 30)\lvec(50 30)\lvec(50 40)\lvec(40 40)\lvec(40
30)\htext(43 33){$2$}

\move(0 0)\lvec(10 10)\lvec(10 0)\lvec(0 0)\lfill f:0.8 \htext(1
5){$0$}\htext(6 1){$1$}

\move(10 0)\lvec(20 10)\lvec(20 0)\lvec(10 0)\lfill f:0.8
\htext(11 5){$1$}\htext(16 1){$0$}

\move(20 0)\lvec(30 10)\lvec(30 0)\lvec(20 0)\lfill f:0.8
\htext(22 5){$0$}\htext(26 2){$1$}

\move(30 0)\lvec(40 10)\lvec(40 0)\lvec(30 0)\lfill f:0.8
\htext(32 5){$1$}\htext(36 2){$0$}

\move(40 0)\lvec(50 10)\lvec(50 0)\lvec(40 0)\lfill f:0.8
\htext(42 5){$0$}\htext(46 2){$1$}

\move(30 40)\lvec(40 50)\lvec(40 40)\lvec(30 40)\htext(36 42){$0$}

\move(40 40)\lvec(50 50)\lvec(50 40)\lvec(40 40)\htext(46 42){$1$}
\end{texdraw}}\ \ \ \ \ $\longleftrightarrow$ \ \ \ \ \
\raisebox{-0.5\height}{\begin{texdraw} \drawdim em \setunitscale
0.16 \linewd 0.5 \fontsize{8}{8}

\htext(10 0){$1$}\htext(20 0){$2$}\htext(30 0){$3$}\htext(40
0){$4$}\htext(50 0){$5$}

\htext(10 -10){$6$}\htext(20 -10){$7$}\htext(30 -10){$8$}\htext(40
-10){$9$}\htext(49 -10){$10$}

\htext(10 -20){$\vdots$}\htext(20 -20){$\vdots$}\htext(30
-20){$\vdots$}\htext(40 -20){$\vdots$}\htext(50 -20){$\vdots$}

\move(11 1)\linewd .3 \lcir r:3 \move(41 1)\linewd .3 \lcir r:3
\move(51 1)\linewd .3 \lcir r:3 \move(31 1)\linewd .3 \lcir r:3
\fontsize{4}{4}\move(54 3)\htext{$2$}
\end{texdraw}}\ \ \ .
\end{center}\vskip 3mm

(2) Let $Y$ be the proper Young wall in Example \ref{ex5.1} (2).
Similar to (1), if we apply $(B_1)$ to \raisebox{-0.25
\height}{\begin{texdraw} \drawdim em \setunitscale 0.16 \linewd
0.5 \fontsize{8}{8} \htext(0 0){$7$} \move(1 1)\linewd .3 \lcir
r:3
\end{texdraw}}, then we have
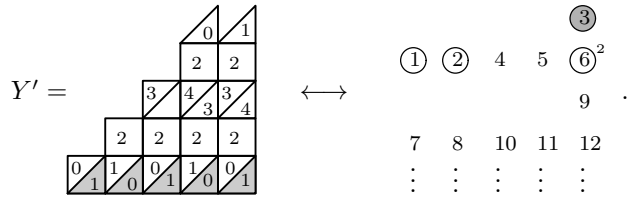
\begin{center}
$Y'=$\raisebox{-0.5\height}{\begin{texdraw}\fontsize{7}{7}\drawdim
mm \setunitscale 0.5 \linewd 0.5

\move(0 0)\lvec(10 0)\lvec(10 10)\lvec(0 10)\lvec(0 0)

\move(10 0)\lvec(20 0)\lvec(20 10)\lvec(10 10)\lvec(10 0)

\move(20 0)\lvec(30 0)\lvec(30 10)\lvec(20 10)\lvec(20 0)

\move(30 0)\lvec(40 0)\lvec(40 10)\lvec(30 10)\lvec(30 0)

\move(40 0)\lvec(50 0)\lvec(50 10)\lvec(40 10)\lvec(40 0)

\move(10 10)\lvec(20 10)\lvec(20 20)\lvec(10 20)\lvec(10
10)\htext(13 13){$2$}

\move(20 10)\lvec(30 10)\lvec(30 20)\lvec(20 20)\lvec(20
10)\htext(23 13){$2$}

\move(30 10)\lvec(40 10)\lvec(40 20)\lvec(30 20)\lvec(30
10)\htext(33 13){$2$}

\move(40 10)\lvec(50 10)\lvec(50 20)\lvec(40 20)\lvec(40
10)\htext(43 13){$2$}
\move(20 20)\lvec(30 20)\lvec(30 30)\lvec(20 30)\lvec(20
20)\lvec(30 30)\htext(21 25){$3$}

\move(30 20)\lvec(40 20)\lvec(40 30)\lvec(30 30)\lvec(30
20)\lvec(40 30)\htext(36 21){$3$}\htext(31 25){$4$}

\move(40 20)\lvec(50 20)\lvec(50 30)\lvec(40 30)\lvec(40
20)\lvec(50 30)\htext(41 25){$3$}\htext(46 21){$4$}
\move(30 30)\lvec(40 30)\lvec(40 40)\lvec(30 40)\lvec(30
30)\htext(33 33){$2$}

\move(40 30)\lvec(50 30)\lvec(50 40)\lvec(40 40)\lvec(40
30)\htext(43 33){$2$}
\move(0 0)\lvec(10 10)\lvec(10 0)\lvec(0 0)\lfill f:0.8 \htext(1
5){$0$}\htext(6 1){$1$}

\move(10 0)\lvec(20 10)\lvec(20 0)\lvec(10 0)\lfill f:0.8
\htext(11 5){$1$}\htext(16 1){$0$}

\move(20 0)\lvec(30 10)\lvec(30 0)\lvec(20 0)\lfill f:0.8
\htext(22 5){$0$}\htext(26 2){$1$}

\move(30 0)\lvec(40 10)\lvec(40 0)\lvec(30 0)\lfill f:0.8
\htext(32 5){$1$}\htext(36 2){$0$}

\move(40 0)\lvec(50 10)\lvec(50 0)\lvec(40 0)\lfill f:0.8
\htext(42 5){$0$}\htext(46 2){$1$}
\move(30 40)\lvec(40 50)\lvec(40 40)\lvec(30 40)\htext(36 41){$0$}

\move(40 40)\lvec(50 50)\lvec(50 40)\lvec(40 40)\htext(46 42){$1$}
\end{texdraw}}\ \ \ \ \ $\longleftrightarrow$ \ \ \ \ \
\raisebox{-0.5\height}{\begin{texdraw} \drawdim em \setunitscale
0.16 \linewd 0.5 \fontsize{8}{8}

\htext(10 0){$1$}\htext(20 0){$2$}\htext(30 0){$4$}\htext(40
0){$5$}\htext(50 10){$3$}\htext(50 0){$6$}

\htext(10 -20){$7$}\htext(20 -20){$8$}\htext(30
-20){$10$}\htext(40 -20){$11$}\htext(50 -10){$9$}\htext(50
-20){$12$}

\htext(10 -30){$\vdots$}\htext(20 -30){$\vdots$}\htext(30
-30){$\vdots$}\htext(40 -30){$\vdots$}\htext(50 -30){$\vdots$}

\move(21 1)\linewd .3 \lcir r:3 \move(51 11)\linewd .3 \fcir f:0.7
r:3 \lcir r:3 \move(51 1)\linewd .3  \lcir r:3 \move(11 1)\linewd
.3 \lcir r:3 \fontsize{4}{4}\move(54 3)\htext{$2$}
\end{texdraw}}\ \ \ .
\end{center}\vskip 3mm}
\end{ex}

Let $\widetilde{Y}$ be the proper Young wall obtained from $Y$ by
applying $(B_1)$ and $(B_2)$ until there is no bead movable up or
removable. Note that $\widetilde{Y}$ does not depend on the order
of steps, hence is uniquely determined.

Set
\begin{equation}
\begin{aligned}
\mathcal{Z}(\Lambda)'&=\{\,Y=(y_k)_{k\geq
1}\in\mathcal{Z}(\Lambda)\,|\,|y_k|\equiv 0\pmod \ell\,\}, \\
\mathcal{Z}(\Lambda)'_{\lambda}&=\mathcal{Z}(\Lambda)'\cap
\mathcal{Z}(\Lambda)_{\lambda} \ \ \ \text{for $\lambda\leq
\Lambda$}.
\end{aligned}
\end{equation}
Note that for $Y\in \Z(\Lambda)$, $Y\in \Z(\Lambda)'$ if and only
if there exists no bead in a runner of type I in its bead
configuration.

\begin{lem}\label{wtA(2n-1)(2)}
Let $Y$ be a proper Young wall in $\Z(\Lambda)$.  Then
$Y\in\mathcal{Z}(\Lambda)_{\Lambda-m\delta}$ for some $m\geq 0$ if
and only if $\widetilde{Y}\in
\mathcal{Z}(\Lambda)'_{\Lambda-m'\delta}$ for some $m'\geq 0$.
\end{lem}
\pf Let $R_k$ be a runner of type I and let $r_k$ be the number of
beads in $R_k$ occurring in $Y$. Note that $\widetilde{Y}\in
\mathcal{Z}(\Lambda)'$ if and only if $r_k=r_{L-k}$ for $1\leq
k\leq n-1$.

We assume that $\frak{g}$ is of type $A_{2n-1}^{(2)}$ (the proof
for $D_{n+1}^{(1)}$ is similar). Suppose that ${\rm
wt}(Y)=\Lambda-m\delta$ for some $m\geq 0$.  By considering the
content of each bead, we see that ${\rm cont}(Y)=\sum_{i=0}^n
c_i\alpha_i + M\delta$ for some $M\geq 0$ where
\begin{equation}
\begin{aligned}
&c_0+c_1= \sum_{j=1}^{2n-2}r_j, \\
&c_i=
\begin{cases}
\sum_{j=i}^{2n-i-1}r_j + \sum_{j=2n-i}^{2n-2}2r_j & \text{if
$2\leq i \leq n-1$}, \\
\sum_{j=n}^{2n-2}r_j & \text{if $i=n$}.
\end{cases}
\end{aligned}
\end{equation}
Since ${\rm cont}(Y)=m\delta$ and
$\delta=\alpha_0+\alpha_1+\sum_{i=2}^{n-1}2\alpha_i+\alpha_n$, it
follows that $c_0+c_1=c_2=\cdots=c_{n-1}=2c_n$, and hence
$r_i=r_{2n-i-1}$ for $1\leq i\leq n-1$. By Lemma
\ref{abacusA(2n-1)(2)}, we have that $\widetilde{Y}\in
\mathcal{Z}(\Lambda)'_{\Lambda-m'\delta}$ for some $m'\geq 0$.

The converse is clear from Lemma \ref{abacusA(2n-1)(2)}. \qed
\vskip 5mm

Fix $m\geq 0$. Let $Y$ be a proper Young wall in
$\mathcal{Z}(\Lambda)_{\Lambda-m\delta}$. Consider its bead
configuration. Let $R_k$ be a runner of type I and let $r_k$ be
the number of beads in $R_k$. By Lemma \ref{wtA(2n-1)(2)},
$r_k=r_{L-k}$ for $1\leq k\leq n-1$ and  as in the case of
$A_{2n}^{(2)}$ or $D_{n+1}^{(2)}$, we can associate a unique
partition $\lambda^{(k)}$ from the beads in $R_k$ and $R_{L-k}$
($1\leq k\leq n-1$) using Frobenius notation. We define
\begin{equation}
\pi_0(Y)=(\lambda^{(1)},\cdots,\lambda^{(n-1)}).
\end{equation}
Note that ${\rm cont}(\widetilde{Y})=m'\delta$ where
$m'=m-\sum_{i=1}^{n-1}|\lambda^{(i)}|$ and $Y$ is uniquely
determined by $\pi_0(Y)$ and $\widetilde{Y}$.

Conversely, suppose that we are given an $(n-1)$-tuple of
partitions $(\lambda^{(1)},\cdots,\lambda^{(n-1)})$ and
$Z\in\mathcal{Z}(\Lambda)'_{\Lambda-m'\delta}$ with
$m=m'+\sum_{i=0}^{n-1}|\lambda^{(i)}|$. For $1\leq k\leq n-1$, by
applying the inverse steps in Lemma \ref{abacusA(2n-1)(2)} to $Z$,
we can place beads at $R_k$ and $R_{L-k}$ whose corresponding
partition is $\lambda^{(k)}$. Then it is easy to see that the
resulting proper Young wall $Y$ satisfies (i)
$Y\in\mathcal{Z}(\Lambda)_{\Lambda-m\delta}$, (ii)
$\pi_0(Y)=(\lambda^{(1)},\cdots,\lambda^{(n-1)})$ and (iii)
$\widetilde{Y}=Z$.

Summarizing the above arguments, we obtain
\begin{prop}\label{psi} For $m\geq 0$, the map
\begin{equation}
\psi : \Z(\Lambda)_{\Lambda-m\delta} \longrightarrow
\bigsqcup_{m_1+m_{2}=m}\cP^{(n-1)}(m_1)\times
\mathcal{Z}(\Lambda)'_{\Lambda-m_2\delta}
\end{equation}
defined by $\psi(Y)=(\pi_0(Y),\widetilde{Y})$ is a bijection.\qed
\end{prop}\vskip 3mm

Before characterizing $\mathcal{Z}(\Lambda)'_{\Lambda-m\delta}$,
we recall some notions of two-colored partitions. Let $\mathbb{N}$
be the set of positive integers. We say that an element in
$\mathbb{N}$ is colored {\it white}. Let
$\underline{\mathbb{N}}=\{
\underline{1},\underline{2},\underline{3},\cdots \}$ be the set of
positive integers colored {\it gray}. Set
$N=\mathbb{N}\cup\underline{\mathbb{N}}$. For $a$ and $b\in N$, we
define $a+b$ to be the colored integer whose value is given by the
ordinary sum of their values, and whose color is white if they
have the same color, gray if not. For example,
$\underline{1}+\underline{2}=3$ and
$2+\underline{3}=\underline{5}$. Also, we assume that $a+0=0+a=a$
for all $a\in N$. A {\it two-colored partition} is a sequence
$\lambda=(\lambda_k)_{k\geq 1}$ of elements in $N\cup\{0\}$ such
that (i) all but a finite number of $\lambda_k$'s are zero (ii)
$\lambda_{k}\geq\lambda_{k+1}$ as ordinary integers ($\lambda$
does not depend on the order of the same-valued integers).  We
also write $\lambda=(1^{m_1},\underline{1}^{m_{\underline{1}}},
2^{m_2},\underline{2}^{m_{\underline{2}}},\cdots)$ where $m_{k}$
is the multiplicity of $k\in N$ in $\lambda$. We say that
$\lambda$ is a two-colored partition of $m$ ($m\geq 0$) if
$\sum_{k\geq 1}k(m_k+m_{\underline{k}})=m$, and write
$|\lambda|=m$.

Let $\mathscr{P}'$ be the set of all two-colored partitions such
that for each $k\geq 1$, both $k$ and $\underline{k}$ do not
appear as a part simultaneously. Note that $\cP'$ is closed under
addition, i.e. $\lambda+\mu=(\lambda_k+\mu_k)_{k\geq 1}\in \cP'$
for $\lambda$, $\mu\in\cP'$. For $Y\in\mathcal{Z}(\Lambda)'$,
consider its bead configuration where all the beads lie in
$R_{\ell}$. Define $\lambda_Y$ to be the two-colored partition in
$\cP'$, where the multiplicity of $k$ (resp. $\underline{k}$) is
determined by the number of white (resp. gray) beads at position
$k\ell$. Note that either one of $m_{k}$ and $m_{\underline{k}}$
is zero from the patterns of proper Young walls. Then the map $Y
\mapsto \lambda_Y$ for $Y\in\mathcal{Z}(\Lambda)'$ is a bijection
between $\mathcal{Z}(\Lambda)'$ and $\mathscr{P}'$.

Also, we may identify $\lambda\in\mathscr{P}'$ with a {\it
two-colored Young diagram} as follows: for each part $k\in
\mathbb{N}$ of $\lambda$, we associate a column with $k$ boxes,
and for each part $\underline{k}\in \underline{\mathbb{N}}$ of
$\lambda$, we associate a column with $k$ boxes where the top box
is colored with gray.

\begin{ex}\label{2colored partition}{\rm (1) If
$\lambda=(1,\underline{2},\underline{3}^2,4,\underline{5})$, then
the corresponding two-colored Young diagram is
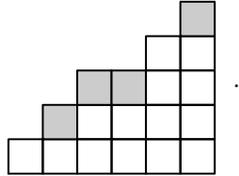
\begin{center}
\raisebox{-0.5\height}{
\begin{texdraw}\drawdim
em \setunitscale 0.13 \linewd 0.5

\move(0 0)\lvec(10 0)\lvec(10 10)\lvec(0 10)\lvec(0 0)

\move(10 0)\lvec(20 0)\lvec(20 10)\lvec(10 10)\lvec(10 0)

\move(20 0)\lvec(30 0)\lvec(30 10)\lvec(20 10)\lvec(20 0)

\move(30 0)\lvec(40 0)\lvec(40 10)\lvec(30 10)\lvec(30 0)

\move(40 0)\lvec(50 0)\lvec(50 10)\lvec(40 10)\lvec(40 0)

\move(50 0)\lvec(60 0)\lvec(60 10)\lvec(50 10)\lvec(50 0)

\move(10 10)\lvec(20 10)\lvec(20 20)\lvec(10 20)\lvec(10 10)\lfill
f:0.8

\move(20 10)\lvec(30 10)\lvec(30 20)\lvec(20 20)\lvec(20 10)

\move(30 10)\lvec(40 10)\lvec(40 20)\lvec(30 20)\lvec(30 10)

\move(40 10)\lvec(50 10)\lvec(50 20)\lvec(40 20)\lvec(40 10)

\move(50 10)\lvec(60 10)\lvec(60 20)\lvec(50 20)\lvec(50 10)

\move(20 20)\lvec(30 20)\lvec(30 30)\lvec(20 30)\lvec(20 20)\lfill
f:0.8

\move(30 20)\lvec(40 20)\lvec(40 30)\lvec(30 30)\lvec(30 20)\lfill
f:0.8

\move(40 20)\lvec(50 20)\lvec(50 30)\lvec(40 30)\lvec(40 20)

\move(50 20)\lvec(60 20)\lvec(60 30)\lvec(50 30)\lvec(50 20)

\move(40 30)\lvec(50 30)\lvec(50 40)\lvec(40 40)\lvec(40 30)

\move(50 30)\lvec(60 30)\lvec(60 40)\lvec(50 40)\lvec(50 30)

\move(50 40)\lvec(60 40)\lvec(60 50)\lvec(50 50)\lvec(50 40)
\lfill f:0.8
\end{texdraw}}\ \ .\end{center}\vskip 3mm

(2) Suppose that $\frak{g}=D_4^{(1)}$ and
$\Lambda=\Lambda_0$.\vskip 3mm

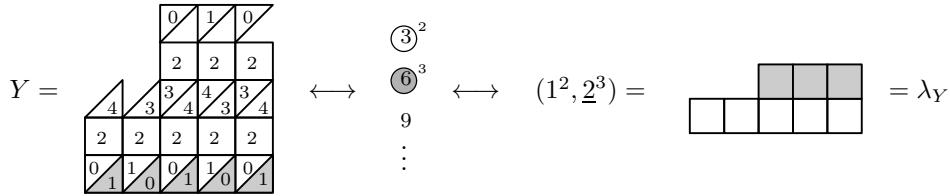
\begin{center}
$Y=$\ \
\ \raisebox{-0.5\height}{\begin{texdraw}\fontsize{7}{7}\drawdim mm
\setunitscale 0.5 \linewd 0.5

\move(0 0)\lvec(10 0)\lvec(10 10)\lvec(0 10)\lvec(0 0)

\move(10 0)\lvec(20 0)\lvec(20 10)\lvec(10 10)\lvec(10 0)

\move(20 0)\lvec(30 0)\lvec(30 10)\lvec(20 10)\lvec(20 0)

\move(30 0)\lvec(40 0)\lvec(40 10)\lvec(30 10)\lvec(30 0)

\move(40 0)\lvec(50 0)\lvec(50 10)\lvec(40 10)\lvec(40 0)
\move(0 10)\lvec(10 10)\lvec(10 20)\lvec(0 20)\lvec(0 10)\htext(3
13){$2$}

\move(10 10)\lvec(20 10)\lvec(20 20)\lvec(10 20)\lvec(10
10)\htext(13 13){$2$}

\move(20 10)\lvec(30 10)\lvec(30 20)\lvec(20 20)\lvec(20
10)\htext(23 13){$2$}

\move(30 10)\lvec(40 10)\lvec(40 20)\lvec(30 20)\lvec(30
10)\htext(33 13){$2$}

\move(40 10)\lvec(50 10)\lvec(50 20)\lvec(40 20)\lvec(40
10)\htext(43 13){$2$}
\move(0 20)\lvec(10 30)\lvec(10 20)\lvec(0 20)\htext(6 21){$4$}

\move(10 20)\lvec(20 30)\lvec(20 20)\lvec(10 20)\htext(16 21){$3$}

\move(20 20)\lvec(30 20)\lvec(30 30)\lvec(20 30)\lvec(20
20)\lvec(30 30)\htext(21 25){$3$}\htext(26 21){$4$}

\move(30 20)\lvec(40 20)\lvec(40 30)\lvec(30 30)\lvec(30
20)\lvec(40 30)\htext(36 21){$3$}\htext(31 25){$4$}

\move(40 20)\lvec(50 20)\lvec(50 30)\lvec(40 30)\lvec(40
20)\lvec(50 30)\htext(41 25){$3$}\htext(46 21){$4$}

\move(20 30)\lvec(30 30)\lvec(30 40)\lvec(20 40)\lvec(20
30)\htext(23 33){$2$}

\move(30 30)\lvec(40 30)\lvec(40 40)\lvec(30 40)\lvec(30
30)\htext(33 33){$2$}

\move(40 30)\lvec(50 30)\lvec(50 40)\lvec(40 40)\lvec(40
30)\htext(43 33){$2$}
\move(0 0)\lvec(10 10)\lvec(10 0)\lvec(0 0)\lfill f:0.8 \htext(1
5){$0$}\htext(6 1){$1$}

\move(10 0)\lvec(20 10)\lvec(20 0)\lvec(10 0)\lfill f:0.8
\htext(11 5){$1$}\htext(16 1){$0$}

\move(20 0)\lvec(30 10)\lvec(30 0)\lvec(20 0)\lfill f:0.8
\htext(22 5){$0$}\htext(26 2){$1$}

\move(30 0)\lvec(40 10)\lvec(40 0)\lvec(30 0)\lfill f:0.8
\htext(32 5){$1$}\htext(36 2){$0$}

\move(40 0)\lvec(50 10)\lvec(50 0)\lvec(40 0)\lfill f:0.8
\htext(42 5){$0$}\htext(46 2){$1$}
\move(20 40)\lvec(30 50)\lvec(20 50)\lvec(20 40)\htext(22 45){$0$}

\move(30 40)\lvec(40 50)\lvec(30 50)\lvec(30 40)\htext(32 45){$1$}

\move(40 40)\lvec(50 50)\lvec(40 50)\lvec(40 40)\htext(42 45){$0$}
\end{texdraw}} \ \ \ $\longleftrightarrow$ \ \ \
\raisebox{-0.5\height}{\begin{texdraw} \drawdim em \setunitscale
0.16 \linewd 0.5 \fontsize{8}{8}

\htext(50 10){$3$}\htext(50 0){$6$}

\htext(50 -10){$9$}

\htext(50 -20){$\vdots$}

\move(51 11)\linewd .3 \lcir r:3 \move(51 1)\linewd .3 \fcir f:0.7
r:3 \lcir r:3 \fontsize{4}{4}\move(54
3)\htext{$3$}\fontsize{4}{4}\move(54 13)\htext{$2$}
\end{texdraw}} \ \ \ $\longleftrightarrow$ \ \ \
$(1^2,\underline{2}^3)=$   \ \ \ \raisebox{-0.5\height}{
\begin{texdraw}\drawdim
em \setunitscale 0.13 \linewd 0.5

\move(-10 0)\lvec(0 0)\lvec(0 10)\lvec(-10 10)\lvec(-10 0)

\move(0 0)\lvec(10 0)\lvec(10 10)\lvec(0 10)\lvec(0 0)

\move(10 0)\lvec(20 0)\lvec(20 10)\lvec(10 10)\lvec(10 0)

\move(20 0)\lvec(30 0)\lvec(30 10)\lvec(20 10)\lvec(20 0)

\move(30 0)\lvec(40 0)\lvec(40 10)\lvec(30 10)\lvec(30 0)

\move(10 10)\lvec(20 10)\lvec(20 20)\lvec(10 20)\lvec(10 10)\lfill
f:0.8

\move(20 10)\lvec(30 10)\lvec(30 20)\lvec(20 20)\lvec(20 10)\lfill
f:0.8

\move(30 10)\lvec(40 10)\lvec(40 20)\lvec(30 20)\lvec(30 10)\lfill
f:0.8
\end{texdraw}}\ \ \ $=\lambda_Y$
\end{center}\vskip 5mm

In general, for a proper Young wall $Y\in \Z(\Lambda)'$,
$\lambda_Y$ can be obtained by identifying
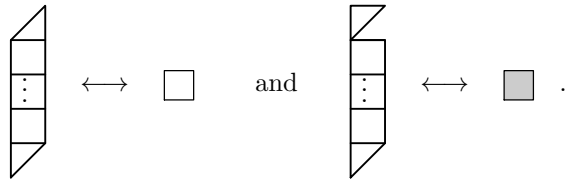
\begin{center}
\raisebox{-0.5\height}{\begin{texdraw}\drawdim em \setunitscale
0.13 \linewd 0.5 \fontsize{5}{5}

\move(0 0)\lvec(0 40)\lvec(10 50)\lvec(10 10)\lvec(0 0)\move(0
10)\lvec(10 10)\move(0 40)\lvec(10 40)\move(0 20)\lvec(10
20)\move(0 30)\lvec(10 30)\htext(3 22){$\vdots$}

\end{texdraw}} \ \ \ $\longleftrightarrow$  \ \ \ \raisebox{-0.3\height}{\begin{texdraw}\drawdim em \setunitscale
0.11 \linewd 0.5 \fontsize{7}{7} \move(0 0)\lvec(0 10)\lvec(10
10)\lvec(10 0)\lvec(0 0)
\end{texdraw} } \ \ \ \ \  and \ \ \ \ \
\raisebox{-0.5\height}{\begin{texdraw}\drawdim em \setunitscale
0.13 \linewd 0.5 \fontsize{7}{7}

\move(0 0)\lvec(0 50)\lvec(10 50)\lvec(0 40)\lvec(10 40)\lvec(10
10)\lvec(0 0)\move(0 10)\lvec(10 10)\move(0 40)\lvec(10 40)\move(0
20)\lvec(10 20)\move(0 30)\lvec(10 30)\htext(3 22){$\vdots$}

\end{texdraw}} \ \ \ $\longleftrightarrow$  \ \ \ \raisebox{-0.3\height}{\begin{texdraw}\drawdim em \setunitscale
0.11 \linewd 0.5 \fontsize{7}{7} \move(0 0)\lvec(0 10)\lvec(10
10)\lvec(10 0)\lvec(0 0)\lfill f:0.8
\end{texdraw}}\ \ \ .
\end{center}\vskip 3mm

 }\end{ex}

\begin{df}\label{lambda01}{\rm
For $\lambda=(\lambda_k)\in\mathscr{P}'$, we define
$\lambda^{0}=(\lambda^{0}_k)$ and $\lambda^{1}=(\lambda^{1}_k)$ to
be the unique partitions satisfying the following conditions:
\begin{itemize}
\item[(i)] $\lambda^{0}\in \mathscr{P}$ and
$\lambda^{1}\in\mathscr{P}'$,

\item[(ii)] $\lambda=\lambda^0+\lambda^1$,

\item[(iii)] $\lambda^1=(\underline{1}^{m_{\underline{1}}},
2^{m_2},\underline{3}^{m_{\underline{3}}},\cdots)$ with
$m_{\underline{1}}\neq 0$, and it is $2$-reduced,

\end{itemize}}
\end{df}

Note that the set of all partitions in $\cP'$ satisfying the
condition (3) can be identified with the set of all $2$-reduced
partitions since the way of coloring a number is unique.

\begin{ex}{\rm
If $\lambda=(1,\underline{2},\underline{3}^2,4,\underline{5})$
given in Example \ref{2colored partition} (1), then
$\lambda^0=(1^2,2^4)$ and
$\lambda^1=(\underline{1}^3,2,\underline{3})$. That is,\vskip 3mm

\begin{center}
$\lambda^0=$\raisebox{-0.5\height}{
\begin{texdraw}\drawdim
em \setunitscale 0.13 \linewd 0.5

\move(0 0)\lvec(10 0)\lvec(10 10)\lvec(0 10)\lvec(0 0)

\move(10 0)\lvec(20 0)\lvec(20 10)\lvec(10 10)\lvec(10 0)

\move(20 0)\lvec(30 0)\lvec(30 10)\lvec(20 10)\lvec(20 0)

\move(30 0)\lvec(40 0)\lvec(40 10)\lvec(30 10)\lvec(30 0)

\move(40 0)\lvec(50 0)\lvec(50 10)\lvec(40 10)\lvec(40 0)

\move(50 0)\lvec(60 0)\lvec(60 10)\lvec(50 10)\lvec(50 0)

\move(20 10)\lvec(30 10)\lvec(30 20)\lvec(20 20)\lvec(20 10)

\move(30 10)\lvec(40 10)\lvec(40 20)\lvec(30 20)\lvec(30 10)

\move(40 10)\lvec(50 10)\lvec(50 20)\lvec(40 20)\lvec(40 10)

\move(50 10)\lvec(60 10)\lvec(60 20)\lvec(50 20)\lvec(50 10)

\end{texdraw}} \ \ \ and \ \ \
$\lambda^1=$\raisebox{-0.5\height}{
\begin{texdraw}\drawdim
em \setunitscale 0.13 \linewd 0.5

\move(0 0)\lvec(10 0)\lvec(10 10)\lvec(0 10)\lvec(0 0)\lfill f:0.8

\move(10 0)\lvec(20 0)\lvec(20 10)\lvec(10 10)\lvec(10 0)\lfill
f:0.8

\move(20 0)\lvec(30 0)\lvec(30 10)\lvec(20 10)\lvec(20 0)\lfill
f:0.8

\move(30 0)\lvec(40 0)\lvec(40 10)\lvec(30 10)\lvec(30 0)

\move(40 0)\lvec(50 0)\lvec(50 10)\lvec(40 10)\lvec(40 0)

\move(30 10)\lvec(40 10)\lvec(40 20)\lvec(30 20)\lvec(30 10)

\move(40 10)\lvec(50 10)\lvec(50 20)\lvec(40 20)\lvec(40 10)

\move(40 20)\lvec(50 20)\lvec(50 30)\lvec(40 30)\lvec(40 20)\lfill
f:0.8

\end{texdraw}}\ \ \ .
\end{center}\vskip 5mm}
\end{ex}

Let $\lambda\in\mathscr{P}'$ be a two-colored partition (or Young
diagram). For a box $b$ in $\lambda$, suppose that $b$ lies in the
$p$th column and the $q$th row. We define the {\it residue of $b$}
to be $p+q \pmod 2$ if $b$ is white, and $p+1 \pmod 2$ if $b$ is
gray. For $\epsilon=0,1$, we define $r_{\epsilon}$ (resp.
$\underline{r}_{\epsilon}$) to be the number of white (resp. gray)
boxes in $\lambda$ with residue $\epsilon$.

\begin{lem}\label{residue} Let $Y$ be a proper Young wall in
$\mathcal{Z}(\Lambda)'$ and let $\lambda=\lambda_Y$ be the
corresponding 2-colored partition in $\cP'$.
\begin{itemize}
\item[(1)] Suppose that $\frak{g}$ is of type $A_{2n-1}^{(2)}$.
Then $Y\in\mathcal{Z}(\Lambda)'_{\Lambda-m\delta}$ for some $m\geq
0$ if and only if $\underline{r}_0=\underline{r}_1$.

\item[(2)] Suppose that $\frak{g}$ is of type $D_{n+1}^{(1)}$.
Then $Y\in\mathcal{Z}(\Lambda)'_{\Lambda-m\delta}$ for some $m\geq
0$ if and only if $r_0=r_1$ and $\underline{r}_0=\underline{r}_1$.
\end{itemize}
\end{lem}
\pf (1) Under the correspondence between $Y$ and $\lambda_Y$, we
observe from the pattern for $\Z(\Lambda)$ that the content of a
white box is $\delta$ and the content of two gray boxes with
different residues is $2\delta$. This proves (1).

(2) Similarly, we see that (i) the content of two boxes of the
same color with different residues is $2\delta$ and (ii) no pair
of two boxes of different color makes a content which is a
multiple of $\delta$. This proves (2). \qed\vskip 3mm

For $Y\in\mathcal{Z}(\Lambda)'_{\Lambda-m\delta}$, we define
\begin{equation}
\pi_1(Y)=((\lambda_Y)^{0},(\lambda_Y)^{1}),
\end{equation}
where we view $(\lambda_Y)^{1}$ as an ordinary $2$-reduced
partition.

\begin{prop}\label{pi1}\mbox{}
\begin{itemize}
\item[(1)] If $\frak{g}$ is of type $A_{2n-1}^{(2)}$, then for
$m\geq 0$, the map
\begin{equation*}
\pi_1: \mathcal{Z}(\Lambda)'_{\Lambda-m\delta} \longrightarrow
\bigsqcup_{m_1+2m_2=m}\mathscr{P}(m_1)\times \mathscr{DP}_0(m_2)
\end{equation*}
is a bijection.

\item[(2)] If $\frak{g}$ is of type $D_{n+1}^{(1)}$, then for
$m\geq 0$, the map
\begin{equation*}
\pi_1: \mathcal{Z}(\Lambda)'_{\Lambda-m\delta} \longrightarrow
\bigsqcup_{m_1+m_2=m}\mathscr{P}_0(m_1)\times \mathscr{DP}_0(m_2),
\end{equation*}
is a bijection where $\mathscr{P}_0(k)$ is the set of all
partitions with empty $2$-core and $2$-weight $k$.
\end{itemize}
\end{prop}
\pf For $Y\in\mathcal{Z}(\Lambda)'_{\Lambda-m\delta}$, let
$\lambda_Y$ be the two-colored partition in $\mathscr{P}'$
corresponding to $Y$. Then there exist unique proper Young walls
$Y_0$ and $Y_1\in\mathcal{Z}(\Lambda)'$ such that
$\lambda_{Y_0}=(\lambda_Y)^{0}$ and
$\lambda_{Y_1}=(\lambda_Y)^{1}$. And we can check that
\begin{equation}\label{cont}
{\rm cont}(Y)={\rm cont}(Y_0)+{\rm cont}(Y_1).
\end{equation}
If $\frak{g}$ is of type $A_{2n-1}^{(2)}$, then each white box in
$\lambda_{Y_0}=(\lambda_Y)^0$ corresponds to a $\delta$-column in
$Y_0$, which implies that ${\rm cont}(Y_0)\in\mathbb{Z}_{\geq
0}\delta$. Suppose that $\frak{g}$ is of type $D_{n+1}^{(1)}$ and
consider the residues of the white boxes in $\lambda_{Y_0}$ and
$\lambda_{Y_1}$. Since there exist even number of white boxes in
each column in $\lambda_{Y_1}$, the number of white boxes in
$\lambda_{Y_1}$ with residue 0 is equal to the number  of white
boxes in $\lambda_{Y_1}$ with residue 1. By Lemma \ref{residue}
(2) and \eqref{cont}, the number of white boxes in $\lambda_{Y_0}$
with residue 0 is equal to the number  of white boxes in
$\lambda_{Y_0}$ with residue 1, which implies that the $2$-core of
$\lambda_{Y_0}$ is empty and ${\rm cont}(Y_0)\in\mathbb{Z}_{\geq
0}\delta$.

Note that for each box $b$ in $\lambda_{Y_1}$, white or gray, the
residue of $b$ is equal to $p+q \pmod 2$ if it is placed in the
$p$th row and the $q$th column. By Lemma \ref{residue}, we have
\begin{equation}
{\rm cont}(Y_1)\in\mathbb{Z}_{\geq 0}\delta  \ \ \ \text{if and
only if}\ \ \ {\rm core}_2(\lambda_{Y_1})=\emptyset.
\end{equation}
Therefore,
$\pi_1(Y)=((\lambda_Y)^{0},(\lambda_Y)^{1})\in\mathscr{P}(m_1)\times
\mathscr{DP}_0(m_2)$ (resp. $\mathscr{P}_0(m_1)\times
\mathscr{DP}_0(m_2)$) for $m_1+2m_2=m$ (resp. $m_1+m_2=m$) if
$\frak{g}$ is of type $A_{2n-1}^{(2)}$ (resp. $D_{n+1}^{(1)}$).
Also, $Y$ is uniquely determined by $\pi_1(Y)$ by definition.

On the other hand, for a given $(\lambda^{(0)},\lambda^{(1)})\in
\mathscr{P}(m_1)\times \mathscr{DP}_0(m_2)$ (or
$\mathscr{P}_0(m_1)\times \mathscr{DP}_0(m_2)$), put
$\lambda=\lambda^{(0)}+\lambda^{(1)}$, where we view
$\lambda^{(1)}$ as a two-colored partition with the color of even
(resp. odd) part white (resp. gray). Let $Y$ be the unique proper
Young wall in $\Z(\Lambda)'$ corresponding to $\lambda$. By the
same argument, we have ${\rm cont}(Y)=m\delta$, where $m_1+2m_2=m$
(resp. $m_1+m_2=m$) if $\frak{g}$ is of type $A_{2n-1}^{(2)}$
(resp. $D_{n+1}^{(1)}$). This correspondence is the inverse map of
$\pi_1$, and hence it is a bijection. \qed

Now, for each $m\geq 0$ and $Y\in
\mathcal{Z}(\Lambda)_{\Lambda-m\delta}$, we define
\begin{equation}\label{pi}
\pi(Y)=(\pi_0(Y),\pi_1(\widetilde{Y})).
\end{equation}
Then by Proposition \ref{psi} and \ref{pi1}, we obtain
\begin{thm}\mbox{}
\begin{itemize}
\item [(1)] If $\frak{g}$ is of type $A_{2n-1}^{(2)}$, then for
$m\geq 0$, the map
\begin{equation*}
\pi : \Z(\Lambda)_{\Lambda-m\delta} \longrightarrow
\bigsqcup_{m_1+2m_2=m}\cP^{(n)}(m_1)\times \mathscr{DP}_0(m_{2})
\end{equation*}
is a bijection.

\item[(2)] If $\frak{g}$ is of type $D_{n+1}^{(1)}$, then  for
$m\geq 0$, the map
\begin{equation*}
\pi : \Z(\Lambda)_{\Lambda-m\delta} \longrightarrow
\bigsqcup_{\sum_{i=1}^3m_i=m}\cP^{(n-1)}(m_1)\times\cP_0(m_2)\times
\mathscr{DP}_0(m_{3})
\end{equation*}
is a bijection.
\end{itemize}
\qed
\end{thm}

\begin{ex}\label{exCh5}{\rm Suppose that $\frak{g}=A_5^{(2)}$.
Consider the following proper Young wall in $\Z(\Lambda)$. \vskip
3mm

\begin{center}
\ \ \  $Y=$\ \ \ \raisebox{-0.5\height}{\begin{texdraw} \drawdim
em \setunitscale 0.16 \linewd 0.5 \fontsize{8}{8}

\htext(10 0){$1$}\htext(20 0){$2$}\htext(30 0){$3$}\htext(40
0){$4$}\htext(50 0){$5$}

\htext(10 -10){$6$}\htext(20 -10){$7$}\htext(30 -10){$8$}\htext(40
-10){$9$}\htext(49 -10){$10$}

\htext(10 -20){$11$}\htext(20 -20){$12$}\htext(30
-20){$13$}\htext(40 -20){$14$}\htext(49 -20){$15$}

\htext(10 -30){$16$}\htext(20 -30){$17$}\htext(30
-30){$18$}\htext(40 -30){$19$}\htext(49 -30){$20$}

\htext(10 -40){$\vdots$}\htext(20 -40){$\vdots$}\htext(30
-40){$\vdots$}\htext(40 -40){$\vdots$}\htext(50 -40){$\vdots$}

\move(12 -9)\linewd .3 \lcir r:3.5

\move(13 -29)\linewd .3 \lcir r:3.5


\move(23 -19)\linewd .3 \lcir r:3.5

\move(32 -9)\linewd .3 \lcir r:3.5


\move(42 1)\linewd .3 \lcir r:3.5

\move(43 -19)\linewd .3 \lcir r:3.5

\move(52 1)\linewd .3 \lcir r:3.5

\move(52 -9)\linewd .3  \fcir f:0.7 r:3.5 \lcir r:3.5

\move(52 -19)\linewd .3 \fcir f:0.7 r:3.5 \lcir r:3.5

\move(52 -29)\linewd .3 \fcir f:0.7 r:3.5 \lcir r:3.5

\fontsize{4}{4}\move(54 4)\htext{$5$} \move(54 -16)\htext{$3$}
\end{texdraw}}
\end{center}\vskip 5mm

Then we have  \ \ \ \ $\widetilde{Y}=$ \ \ \
\raisebox{-0.5\height}{\begin{texdraw} \drawdim em \setunitscale
0.16 \linewd 0.5 \fontsize{8}{8} \htext(50 0){$5$}

\htext(49 -10){$10$}

\htext(49 -20){$15$}

\htext(49 -30){$20$}

\htext(50 -40){$\vdots$}

\move(52 1)\linewd .3 \fcir f:0.7 r:3.5 \lcir r:3.5

\move(52 -9)\linewd .3  \lcir r:3.5

\move(52 -19)\linewd .3 \lcir r:3.5

\move(52 -29)\linewd .3 \fcir f:0.7 r:3.5 \lcir r:3.5

\fontsize{4}{4}\move(54 4)\htext{$5$} \move(54 -16)\htext{$3$}
\end{texdraw}}\ \ \
\ \ \ and \ \ \  $\pi_0(Y)=(\lambda^{(1)},\lambda^{(2)})$, \ \ \
where\vskip 5mm
\begin{equation}
\begin{split}
\lambda^{(1)}&=((1,3)|(0,2))=(1,2^2,3), \\
\lambda^{(2)}&=((2)|(1))=(1^2,2).
\end{split}
\end{equation}
Also, the 2-colored partition in $\cP'$ corresponding to
$\widetilde{Y}$ is
$\lambda_{\widetilde{Y}}=(\underline{1}^5,2,3^3,\underline{4})$
with $(\lambda_{\widetilde{Y}})^0=(1^4)$ and
$(\lambda_{\widetilde{Y}})^1=(\underline{1}^5,2^4,\underline{3})$.
Since the $2$-core of $(\lambda_{\widetilde{Y}})^1$ is empty and
$|\lambda_{\widetilde{Y}}|=20$, we have
$\widetilde{Y}\in\Z(\Lambda)'_{\Lambda-20\delta}$ by Proposition
\ref{pi1} (1), and $\pi_1(\widetilde{Y})=((1^5,2^4,3),(1^4))$.
Hence, we have $Y\in\Z(\Lambda)_{\Lambda-32\delta}$ and
\begin{equation}
\pi(Y)=((1,2^2,3),(1^2,2),(1^4),(1^5,2^4,3)).
\end{equation}
}
\end{ex}

Hence, we recover the formulas for the string functions in
\cite{Kac90} by a new combinatorial way.
\begin{cor}\mbox{}

\begin{itemize}
\item [(1)] If $\frak{g}$ is of type $A_{2n-1}^{(2)}$, then we
have
$\Sigma^{\Lambda}_{\Lambda}(q)=\dfrac{1}{(q)_{\infty}^n(q^2)_{\infty}}$.

\item[(2)] If $\frak{g}$ is of type $D_{n+1}^{(1)}$, then we have
$\Sigma^{\Lambda}_{\Lambda}(q)=\dfrac{1}{(q)_{\infty}^{n+2}}$.
\end{itemize}
\end{cor}
\pf By Theorem \ref{ZAn1} and \eqref{YAn1}, we have
\begin{equation}
\sum_{m\geq 0}|\mathscr{DP}_0(m)|q^m=\frac{1}{(q)_{\infty}}, \ \ \
\ \sum_{m\geq 0}|\mathscr{P}_0(m)|q^m=\frac{1}{(q)^2_{\infty}}.
\end{equation}
Thus, we obtain the desired formula for
$\Sigma^{\Lambda}_{\Lambda}(q)$. \qed

\section{$B_{n}^{(1)}$-case}

Suppose that $\frak{g}$ is of type $B_{n}^{(1)}$, and $\Lambda$ is
a dominant integral weight of level $1$. The patterns for
$\Z(\Lambda)$ are given as follows:\vskip 5mm

\begin{center}
$\Lambda=\Lambda_0$ : \raisebox{-1\height}{
\begin{texdraw}\textref h:C v:C \fontsize{6}{6}\selectfont \drawdim mm
\setunitscale 5 \newcommand{\dtri}{ \bsegment \move(-1 0)\lvec(0
1)\lvec(0 0)\lvec(-1 0)\ifill f:0.7 \esegment } \move(0 0)\dtri
\move(-1 0)\dtri \move(-2 0)\dtri \move(-3 0)\dtri \move(0
0)\rlvec(-4.3 0) \move(0 1)\rlvec(-4.3 0) \move(0 2)\rlvec(-4.3 0)
\move(0 3.5)\rlvec(-4.3 0) \move(0 4.5)\rlvec(-4.3 0) \move(0
5.5)\rlvec(-4.3 0) \move(0 6.5)\rlvec(-4.3 0) \move(0
8)\rlvec(-4.3 0) \move(0 9)\rlvec(-4.3 0) \move(0 10)\rlvec(-4.3
0) \move(0 11)\rlvec(-4.3 0) \move(0 0)\rlvec(0 11.3) \move(-1
0)\rlvec(0 11.3) \move(-2 0)\rlvec(0 11.3) \move(-3 0)\rlvec(0
11.3) \move(-4 0)\rlvec(0 11.3) \move(-1 0)\rlvec(1 1) \move(-2
0)\rlvec(1 1) \move(-3 0)\rlvec(1 1) \move(-4 0)\rlvec(1 1)
\move(-1 9)\rlvec(1 1) \move(-2 9)\rlvec(1 1) \move(-3 9)\rlvec(1
1) \move(-4 9)\rlvec(1 1) \move(0 5)\rlvec(-4.3 0) \htext(-0.3
0.25){$1$} \htext(-0.75 0.75){$0$} \htext(-0.5 1.5){$2$}
\vtext(-0.5 2.75){$\cdots$} \htext(-0.5 4){$n\!\!-\!\!1$}
\htext(-0.5 6){$n\!\!-\!\!1$} \htext(-0.5 8.5){$2$} \htext(-0.3
9.25){$1$} \htext(-0.75 9.75){$0$} \htext(-0.5 10.5){$2$}
\htext(-2.3 0.25){$1$} \htext(-2.75 0.75){$0$} \htext(-2.5
1.5){$2$} \vtext(-2.5 2.75){$\cdots$} \htext(-2.5
4){$n\!\!-\!\!1$} \htext(-2.5 6){$n\!\!-\!\!1$} \htext(-2.5
8.5){$2$} \htext(-2.3 9.25){$1$} \htext(-2.75 9.75){$0$}
\htext(-2.5 10.5){$2$} \htext(-1.3 0.25){$0$} \htext(-1.75
0.75){$1$} \htext(-1.5 1.5){$2$} \vtext(-1.5 2.75){$\cdots$}
\htext(-1.5 4){$n\!\!-\!\!1$} \htext(-1.5 6){$n\!\!-\!\!1$}
\htext(-1.5 8.5){$2$} \htext(-1.3 9.25){$0$} \htext(-1.75
9.75){$1$} \htext(-1.5 10.5){$2$} \htext(-3.3 0.25){$0$}
\htext(-3.75 0.75){$1$} \htext(-3.5 1.5){$2$} \vtext(-3.5
2.75){$\cdots$} \htext(-3.5 4){$n\!\!-\!\!1$} \htext(-3.5
6){$n\!\!-\!\!1$} \htext(-3.5 8.5){$2$} \htext(-3.3 9.25){$0$}
\htext(-3.75 9.75){$1$} \htext(-3.5 10.5){$2$} \htext(-0.5
4.75){$n$} \htext(-2.5 4.75){$n$} \htext(-1.5 5.25){$n$}
\htext(-3.5 5.25){$n$} \htext(-1.5 4.75){$n$} \htext(-3.5
4.75){$n$} \htext(-0.5 5.25){$n$} \htext(-2.5 5.25){$n$}
\vtext(-0.5 7.25){$\cdots$} \vtext(-1.5 7.25){$\cdots$}
\vtext(-2.5 7.25){$\cdots$} \vtext(-3.5 7.25){$\cdots$}
\end{texdraw}}\hskip 1cm $\Lambda=\Lambda_1$:\raisebox{-1\height}{
\begin{texdraw} \textref h:C v:C \fontsize{6}{6}\selectfont \drawdim mm
\setunitscale 5 \newcommand{\dtri}{ \bsegment \move(-1 0)\lvec(0
1)\lvec(0 0)\lvec(-1 0)\ifill f:0.7 \esegment } \move(0 0)\dtri
\move(-1 0)\dtri \move(-2 0)\dtri \move(-3 0)\dtri \move(0
0)\rlvec(-4.3 0) \move(0 1)\rlvec(-4.3 0) \move(0 2)\rlvec(-4.3 0)
\move(0 3.5)\rlvec(-4.3 0) \move(0 4.5)\rlvec(-4.3 0) \move(0
5.5)\rlvec(-4.3 0) \move(0 6.5)\rlvec(-4.3 0) \move(0
8)\rlvec(-4.3 0) \move(0 9)\rlvec(-4.3 0) \move(0 10)\rlvec(-4.3
0) \move(0 11)\rlvec(-4.3 0) \move(0 0)\rlvec(0 11.3) \move(-1
0)\rlvec(0 11.3) \move(-2 0)\rlvec(0 11.3) \move(-3 0)\rlvec(0
11.3) \move(-4 0)\rlvec(0 11.3) \move(-1 0)\rlvec(1 1) \move(-2
0)\rlvec(1 1) \move(-3 0)\rlvec(1 1) \move(-4 0)\rlvec(1 1)
\move(-1 9)\rlvec(1 1) \move(-2 9)\rlvec(1 1) \move(-3 9)\rlvec(1
1) \move(-4 9)\rlvec(1 1) \move(0 5)\rlvec(-4.3 0) \htext(-0.3
0.25){$0$} \htext(-0.75 0.75){$1$} \htext(-0.5 1.5){$2$}
\vtext(-0.5 2.75){$\cdots$} \htext(-0.5 4){$n\!\!-\!\!1$}
\htext(-0.5 6){$n\!\!-\!\!1$} \htext(-0.5 8.5){$2$} \htext(-0.3
9.25){$0$} \htext(-0.75 9.75){$1$} \htext(-0.5 10.5){$2$}
\htext(-2.3 0.25){$0$} \htext(-2.75 0.75){$1$} \htext(-2.5
1.5){$2$} \vtext(-2.5 2.75){$\cdots$} \htext(-2.5
4){$n\!\!-\!\!1$} \htext(-2.5 6){$n\!\!-\!\!1$} \htext(-2.5
8.5){$2$} \htext(-2.3 9.25){$0$} \htext(-2.75 9.75){$1$}
\htext(-2.5 10.5){$2$} \htext(-1.3 0.25){$1$} \htext(-1.75
0.75){$0$} \htext(-1.5 1.5){$2$} \vtext(-1.5 2.75){$\cdots$}
\htext(-1.5 4){$n\!\!-\!\!1$} \htext(-1.5 6){$n\!\!-\!\!1$}
\htext(-1.5 8.5){$2$} \htext(-1.3 9.25){$1$} \htext(-1.75
9.75){$0$} \htext(-1.5 10.5){$2$} \htext(-3.3 0.25){$1$}
\htext(-3.75 0.75){$0$} \htext(-3.5 1.5){$2$} \vtext(-3.5
2.75){$\cdots$} \htext(-3.5 4){$n\!\!-\!\!1$} \htext(-3.5
6){$n\!\!-\!\!1$} \htext(-3.5 8.5){$2$} \htext(-3.3 9.25){$1$}
\htext(-3.75 9.75){$0$} \htext(-3.5 10.5){$2$} \htext(-0.5
4.75){$n$} \htext(-2.5 4.75){$n$} \htext(-1.5 5.25){$n$}
\htext(-3.5 5.25){$n$} \htext(-1.5 4.75){$n$} \htext(-3.5
4.75){$n$} \htext(-0.5 5.25){$n$} \htext(-2.5 5.25){$n$}
\vtext(-0.5 7.25){$\cdots$} \vtext(-1.5 7.25){$\cdots$}
\vtext(-2.5 7.25){$\cdots$} \vtext(-3.5 7.25){$\cdots$}
\end{texdraw}}\vskip 5mm

$\Lambda=\Lambda_n$ : \raisebox{-1\height}{
\begin{texdraw}\textref h:C v:C \fontsize{6}{6}\selectfont \drawdim mm
\setunitscale 5 \move(0 0)\lvec(-4 0)\lvec(-4 0.5)\lvec(0
0.5)\ifill f:0.7 \move(0 0)\rlvec(-4.3 0) \move(0 1)\rlvec(-4.3 0)
\move(0 2)\rlvec(-4.3 0) \move(0 3.5)\rlvec(-4.3 0) \move(0
4.5)\rlvec(-4.3 0) \move(0 5.5)\rlvec(-4.3 0) \move(0
6.5)\rlvec(-4.3 0) \move(0 8)\rlvec(-4.3 0) \move(0 9)\rlvec(-4.3
0) \move(0 10)\rlvec(-4.3 0) \move(0 11)\rlvec(-4.3 0) \move(0
0)\rlvec(0 11.3) \move(-1 0)\rlvec(0 11.3) \move(-2 0)\rlvec(0
11.3) \move(-3 0)\rlvec(0 11.3) \move(-4 0)\rlvec(0 11.3) \move(-1
4.5)\rlvec(1 1) \move(-2 4.5)\rlvec(1 1) \move(-3 4.5)\rlvec(1 1)
\move(-4 4.5)\rlvec(1 1) \htext(-0.5 0.25){$n$} \htext(-0.5
1.5){$n\!\!-\!\!1$} \vtext(-0.5 2.75){$\cdots$} \htext(-0.5
4){$2$} \htext(-0.3 4.75){$1$} \htext(-0.75 5.25){$0$} \htext(-0.5
6){$2$} \vtext(-0.5 7.25){$\cdots$} \htext(-0.5
8.5){$n\!\!-\!\!1$} \htext(-0.5 9.25){$n$} \htext(-0.5
10.5){$n\!\!-\!\!1$} \htext(-2.5 0.25){$n$} \htext(-2.5
1.5){$n\!\!-\!\!1$} \vtext(-2.5 2.75){$\cdots$} \htext(-2.5
4){$2$} \htext(-2.3 4.75){$1$} \htext(-2.75 5.25){$0$} \htext(-2.5
6){$2$} \vtext(-2.5 7.25){$\cdots$} \htext(-2.5
8.5){$n\!\!-\!\!1$} \htext(-2.5 9.25){$n$} \htext(-2.5
10.5){$n\!\!-\!\!1$} \htext(-1.5 0.75){$n$} \htext(-1.5
1.5){$n\!\!-\!\!1$} \vtext(-1.5 2.75){$\cdots$} \htext(-1.5
4){$2$} \htext(-1.3 4.75){$0$} \htext(-1.75 5.25){$1$} \htext(-1.5
6){$2$} \vtext(-1.5 7.25){$\cdots$} \htext(-1.5
8.5){$n\!\!-\!\!1$} \htext(-1.5 9.75){$n$} \htext(-1.5
10.5){$n\!\!-\!\!1$} \htext(-3.5 0.75){$n$} \htext(-3.5
1.5){$n\!\!-\!\!1$} \vtext(-3.5 2.75){$\cdots$} \htext(-3.5
4){$2$} \htext(-3.3 4.75){$0$} \htext(-3.75 5.25){$1$} \htext(-3.5
6){$2$} \vtext(-3.5 7.25){$\cdots$} \htext(-3.5
8.5){$n\!\!-\!\!1$} \htext(-3.5 9.75){$n$} \htext(-3.5
10.5){$n\!\!-\!\!1$} \move(0 0.5)\rlvec(-4.3 0) \move(0
9.5)\rlvec(-4.3 0) \htext(-3.5 9.25){$n$} \htext(-3.5 0.25){$n$}
\htext(-2.5 9.75){$n$} \htext(-2.5 0.75){$n$} \htext(-1.5
0.25){$n$} \htext(-1.5 9.25){$n$} \htext(-0.5 9.75){$n$}
\htext(-0.5 0.75){$n$}
\end{texdraw}}
\end{center}\vskip 5mm
Set $\ell=2n$, which is the number of blocks in a $\delta$-column.
We define the {\it abacus of type $B_n^{(1)}$} as follows:\vskip
5mm
\begin{center}
\begin{texdraw}
\drawdim em \setunitscale 0.13 \linewd 0.5 \fontsize{10}{10}

\htext(0 0){$1$}\htext(20 0){$2$}\htext(35 0){$\cdots$}\htext(50
0){$n$}\htext(65 0){$\cdots$}\htext(80 0){$\ell-1$}\htext(110
0){$\ell$}

\htext(-7 -15){$\ell+1$}\htext(13 -15){$\ell+2$}\htext(35
-15){$\cdots$}\htext(50 -15){$3n$}\htext(65
-15){$\cdots$}\htext(80 -15){$2\ell-1$}\htext(110 -15){$2\ell$}

\htext(0 -30){$\vdots$}\htext(20 -30){$\vdots$}\htext(90
-30){$\vdots$}\htext(110 -30){$\vdots$}\htext(50 -30){$\vdots$}
\end{texdraw}
\end{center}\vskip 5mm

For $1\leq k <\ell$, let $R_k$ be the set of all integers $s\equiv
k \pmod \ell$. For $1\leq k <\ell$ ($k\neq n$), we assume that the
$k$th runner $R_k$ is of type I. We also assume that $R_n$ is of
type III (resp. II) and $R_{\ell}$ is of type II (resp. III) if
$\Lambda=\Lambda_n$ (resp. $\Lambda=\Lambda_0,\Lambda_1$) (cf.
Section 5 and 6).

For $Y\in\Z(\Lambda)$, let $\{\,|y_1|,\cdots,|y_r|\,\}$ be the set
of all non-zero parts in $|Y|$. Then by definition of
$\Z(\Lambda)$, the numbers $|y_k|$'s ($1\leq k\leq r$) are
distinct except when $|y_k|\equiv 0\pmod n$. We define the {\it
bead configuration of $Y$} to be the set of $r$ beads
$b_1,\cdots,b_r$ placed in the above abacus where $b_k$ is placed
at $|y_k|$, and the color $c(b_k)$ of $b_k$ is determined by
\eqref{color of bead}. Then $Y$ is uniquely determined by its bead
configuration.

The algorithms of moving and removing beads in the abacus are as
follows (the proof is similar to those of Lemma \ref{abacusA2n2}
and Lemma \ref{abacusA(2n-1)(2)}):
\begin{lem}\label{abacusB(n)(1)}
Let $Y$ be a proper Young wall in $\Z(\Lambda)$ and let $Y'$ be
the proper Young wall which is obtained by applying one of the
following processes to the bead configuration of $Y$:
\begin{itemize}
\item[$(B_1)$] if $b$ is a bead at $s$ in a runner of type I and
there is no bead at $s-\ell$, then move $b$ one position up and
change the color of the beads at $k$ {\rm ($s-\ell< k < s$)} in
the runner of type III,

\item[$(B_2)$] if $b$ is a bead at $s$ in the runner of type II,
then move $b$ one position up along the runner and change the
color of the beads at $k$ {\rm ($s-\ell< k < s$)} in the runner of
type III,

\item[$(B_3)$] if $b$ and $b'$ are beads at $s$ and $\ell-s$ {\rm
($1\leq s\leq n-1$)} respectively, then remove $b$ and $b'$
simultaneously. Also if $R_n$ is of type III, then change the
color of the beads at $n$.

\item[$(B_4)$] if there exists at least one bead at $\ell$ and
$R_{\ell}$ is of type II, then remove one bead at $\ell$ and
change the color of the beads at $n$.

\item[$(B_5)$] if there exists at least two bead at $n$ and
$R_{n}$ is of type II, then remove two beads at $n$.
\end{itemize}
\noindent Then we have $\wt(Y')=\wt(Y)+\delta$.\qed
\end{lem}
\begin{ex}{\rm Suppose that $\frak{g}=B_3^{(1)}$.

(1) Let $Y$ be a proper Young wall in $\Z(\Lambda_3)$ given below.

\begin{center}
$Y=$\raisebox{-0.4\height}{
\begin{texdraw}
\drawdim em \setunitscale 0.13 \linewd 0.5

\move(10 0)\lvec(20 0)\lvec(20 10)\lvec(10 10)\lvec(10 0)\htext(13
1){\tiny $3$}

\move(20 0)\lvec(30 0)\lvec(30 10)\lvec(20 10)\lvec(20 0)\htext(23
1){\tiny $3$}

\move(30 0)\lvec(40 0)\lvec(40 10)\lvec(30 10)\lvec(30 0)\htext(33
1){\tiny $3$}

\move(40 0)\lvec(50 0)\lvec(50 10)\lvec(40 10)\lvec(40 0)\htext(43
1){\tiny $3$}

\move(50 0)\lvec(60 0)\lvec(60 10)\lvec(50 10)\lvec(50 0)\htext(53
1){\tiny $3$}
\move(10 0)\lvec(20 0)\lvec(20 5)\lvec(10 5)\lvec(10 0)\lfill
f:0.8 \htext(13 6){\tiny $3$}

\move(20 0)\lvec(30 0)\lvec(30 5)\lvec(20 5)\lvec(20 0)\lfill
f:0.8 \htext(23 6){\tiny $3$}

\move(30 0)\lvec(40 0)\lvec(40 5)\lvec(30 5)\lvec(30 0)\lfill
f:0.8 \htext(33 6){\tiny $3$}

\move(40 0)\lvec(50 0)\lvec(50 5)\lvec(40 5)\lvec(40 0)\lfill
f:0.8 \htext(43 6){\tiny $3$}

\move(50 0)\lvec(60 0)\lvec(60 5)\lvec(50 5)\lvec(50 0)\lfill
f:0.8 \htext(53 6){\tiny $3$}
\move(20 10)\lvec(30 10)\lvec(30 20)\lvec(20 20)\lvec(20
10)\htext(23 13){$_2$}

\move(30 10)\lvec(40 10)\lvec(40 20)\lvec(30 20)\lvec(30
10)\htext(33 13){$_2$}

\move(40 10)\lvec(50 10)\lvec(50 20)\lvec(40 20)\lvec(40
10)\htext(43 13){$_2$}

\move(50 10)\lvec(60 10)\lvec(60 20)\lvec(50 20)\lvec(50
10)\htext(53 13){$_2$}
\move(20 20)\lvec(30 30)\lvec(30 20)\htext(25 21){\tiny $0$}

\move(30 20)\lvec(30 30)\lvec(40 30)\lvec(40 20)\lvec(30
20)\lvec(40 30)\htext(35 21){\tiny $1$}\htext(31 26){\tiny $0$}

\move(40 20)\lvec(40 30)\lvec(50 30)\lvec(50 20)\lvec(40
20)\lvec(50 30)\htext(45 21){\tiny $0$}\htext(41 26){\tiny $1$}

\move(50 20)\lvec(50 30)\lvec(60 30)\lvec(60 20)\lvec(50
20)\lvec(60 30)\htext(55 21){\tiny $1$}\htext(51 26){\tiny $0$}

\move(50 50)\lvec(60 50)\lvec(60 60)\lvec(50 60)\lvec(50
50)\htext(53 53){$_2$}

\move(30 30)\lvec(40 30)\lvec(40 40)\lvec(30 40)\lvec(30
30)\htext(33 33){$_2$}

\move(40 30)\lvec(50 30)\lvec(50 40)\lvec(40 40)\lvec(40
30)\htext(43 33){$_2$}

\move(50 30)\lvec(60 30)\lvec(60 40)\lvec(50 40)\lvec(50
30)\htext(53 33){$_2$}
\move(40 40)\lvec(50 40)\lvec(50 45)\lvec(40 45)\lvec(40
40)\htext(43 41){\tiny $3$}
\move(50 40)\lvec(60 40)\lvec(60 45)\lvec(50 45)\lvec(50
40)\htext(53 41){\tiny $3$}

\move(50 40)\lvec(60 40)\lvec(60 50)\lvec(50 50)\lvec(50
40)\htext(53 46){\tiny $3$}
\end{texdraw}}\ \ \ \ \ $\longleftrightarrow$ \ \ \ \ \
\raisebox{-0.5\height}{\begin{texdraw} \drawdim em \setunitscale
0.16 \linewd 0.5 \fontsize{8}{8}

\htext(10 0){$1$}\htext(20 0){$2$}\htext(30 0){$3$}\htext(40
0){$4$}\htext(50 0){$5$}\htext(60 0){$6$}

\htext(10 -10){$7$}\htext(20 -10){$8$}\htext(30 -10){$9$}\htext(40
-10){$10$}\htext(49 -10){$11$}\htext(59 -10){$12$}

\htext(10 -20){$\vdots$}\htext(20 -20){$\vdots$}\htext(30
-20){$\vdots$}\htext(40 -20){$\vdots$}\htext(50
-20){$\vdots$}\htext(60 -20){$\vdots$}

\move(11 1)\linewd .3 \lcir r:3 \move(51 1)\linewd .3 \lcir r:3
\move(61 1)\linewd .3 \lcir r:3 \move(31 1)\linewd .3 \lcir r:3
\move(21 -9) \linewd .3 \lcir r:3
\end{texdraw}}\vskip 3mm
\end{center}
If we apply $(B_3)$ to \raisebox{-0.25 \height}{\begin{texdraw}
\drawdim em \setunitscale 0.16 \linewd 0.5 \fontsize{8}{8}
\htext(0 0){$1$} \move(1 1)\linewd .3 \lcir r:3 \end{texdraw}} and
\raisebox{-0.25 \height}{\begin{texdraw} \drawdim em \setunitscale
0.16 \linewd 0.5 \fontsize{8}{8} \htext(0 0){$5$} \move(1
1)\linewd .3 \lcir r:3 \end{texdraw}}, then
\begin{center}
$Y'=$\raisebox{-0.4\height}{
\begin{texdraw}
\drawdim em \setunitscale 0.13 \linewd 0.5

\move(30 0)\lvec(40 0)\lvec(40 10)\lvec(30 10)\lvec(30 0)\htext(33
1){\tiny $3$}

\move(40 0)\lvec(50 0)\lvec(50 10)\lvec(40 10)\lvec(40 0)\htext(43
1){\tiny $3$}

\move(50 0)\lvec(60 0)\lvec(60 10)\lvec(50 10)\lvec(50 0)\htext(53
1){\tiny $3$}
\move(30 0)\lvec(40 0)\lvec(40 5)\lvec(30 5)\lvec(30 0)\lfill
f:0.8 \htext(33 6){\tiny $3$}

\move(40 0)\lvec(50 0)\lvec(50 5)\lvec(40 5)\lvec(40 0)\lfill
f:0.8 \htext(43 6){\tiny $3$}

\move(50 0)\lvec(60 0)\lvec(60 5)\lvec(50 5)\lvec(50 0)\lfill
f:0.8 \htext(53 6){\tiny $3$}

\move(30 10)\lvec(40 10)\lvec(40 20)\lvec(30 20)\lvec(30
10)\htext(33 13){$_2$}

\move(40 10)\lvec(50 10)\lvec(50 20)\lvec(40 20)\lvec(40
10)\htext(43 13){$_2$}

\move(50 10)\lvec(60 10)\lvec(60 20)\lvec(50 20)\lvec(50
10)\htext(53 13){$_2$}

\move(30 20)\lvec(30 30)\lvec(40 30)\lvec(30 20)\lvec(40
30)\htext(31 26){\tiny $0$}

\move(40 20)\lvec(40 30)\lvec(50 30)\lvec(50 20)\lvec(40
20)\lvec(50 30)\htext(45 21){\tiny $0$}\htext(41 26){\tiny $1$}

\move(50 20)\lvec(50 30)\lvec(60 30)\lvec(60 20)\lvec(50
20)\lvec(60 30)\htext(55 21){\tiny $1$}\htext(51 26){\tiny $0$}

\move(50 50)\lvec(60 50)\lvec(60 60)\lvec(50 60)\lvec(50
50)\htext(53 53){$_2$}

\move(40 30)\lvec(50 30)\lvec(50 40)\lvec(40 40)\lvec(40
30)\htext(43 33){$_2$}

\move(50 30)\lvec(60 30)\lvec(60 40)\lvec(50 40)\lvec(50
30)\htext(53 33){$_2$}
\move(40 40)\lvec(50 40)\lvec(50 45)\lvec(40 45)\lvec(40
40)\htext(43 41){\tiny $3$}
\move(50 40)\lvec(60 40)\lvec(60 45)\lvec(50 45)\lvec(50
40)\htext(53 41){\tiny $3$}

\move(50 40)\lvec(60 40)\lvec(60 50)\lvec(50 50)\lvec(50
40)\htext(53 46){\tiny $3$}
\end{texdraw}}\ \ \ \ \ $\longleftrightarrow$ \ \ \ \ \
\raisebox{-0.5\height}{\begin{texdraw} \drawdim em \setunitscale
0.16 \linewd 0.5 \fontsize{8}{8}

\htext(10 0){$1$}\htext(20 0){$2$}\htext(30 0){$3$}\htext(40
0){$4$}\htext(50 0){$5$}\htext(60 0){$6$}

\htext(10 -10){$7$}\htext(20 -10){$8$}\htext(30 -10){$9$}\htext(40
-10){$10$}\htext(49 -10){$11$}\htext(59 -10){$12$}

\htext(10 -20){$\vdots$}\htext(20 -20){$\vdots$}\htext(30
-20){$\vdots$}\htext(40 -20){$\vdots$}\htext(50
-20){$\vdots$}\htext(60 -20){$\vdots$}

\move(61 1)\linewd .3  \lcir r:3 \move(31 1)\linewd .3 \fcir f:0.7
r:3 \lcir r:3 \move(21 -9) \linewd .3 \lcir r:3
\end{texdraw}}\ \ \ . \vskip 3mm
\end{center}

(2) Let $Y$ be a proper Young wall in $\Z(\Lambda_0)$ given below.

\begin{center}
$Y=$\raisebox{-0.5\height}{
\begin{texdraw}
\drawdim em \setunitscale 0.13 \linewd 0.5

\move(10 0)\lvec(20 0)\lvec(20 10)\lvec(10 10)\lvec(10 0)\htext(12
6){\tiny $1$}

\move(20 0)\lvec(30 0)\lvec(30 10)\lvec(20 10)\lvec(20 0)\htext(22
6){\tiny $0$}

\move(30 0)\lvec(40 0)\lvec(40 10)\lvec(30 10)\lvec(30 0)\htext(32
6){\tiny $1$}

\move(40 0)\lvec(50 0)\lvec(50 10)\lvec(40 10)\lvec(40 0)\htext(42
6){\tiny $0$}
\move(10 0)\lvec(20 10)\lvec(20 0)\lvec(10 0)\lfill f:0.8
\htext(16 2){\tiny $0$}

\move(20 0)\lvec(30 10)\lvec(30 0)\lvec(20 0)\lfill f:0.8
\htext(26 2){\tiny $1$}

\move(30 0)\lvec(40 10)\lvec(40 0)\lvec(30 0)\lfill f:0.8
\htext(36 2){\tiny $0$}

\move(40 0)\lvec(50 10)\lvec(50 0)\lvec(40 0)\lfill f:0.8
\htext(46 2){\tiny $1$}
\move(20 10)\lvec(30 10)\lvec(30 20)\lvec(20 20)\lvec(20
10)\htext(23 13){$_2$}

\move(30 10)\lvec(40 10)\lvec(40 20)\lvec(30 20)\lvec(30
10)\htext(33 13){$_2$}

\move(40 10)\lvec(50 10)\lvec(50 20)\lvec(40 20)\lvec(40
10)\htext(43 13){$_2$}
\move(20 20)\lvec(30 20)\lvec(30 25)\lvec(20 25)\lvec(20 20)
\htext(23 21){\tiny $3$}

\move(30 20)\lvec(40 20)\lvec(40 25)\lvec(30 25)\lvec(30 20)
\htext(33 21){\tiny $3$}

\move(40 20)\lvec(50 20)\lvec(50 25)\lvec(40 25)\lvec(40 20)
\htext(43 21){\tiny $3$}

\move(40 25)\lvec(50 25)\lvec(50 30)\lvec(40 30)\lvec(40 25)
\htext(43 26){\tiny $3$}
\end{texdraw}} \ \ \ \ \ $\longleftrightarrow$ \ \ \ \ \
\raisebox{-0.5\height}{\begin{texdraw} \drawdim em \setunitscale
0.16 \linewd 0.5 \fontsize{8}{8}

\htext(10 0){$1$}\htext(20 0){$2$}\htext(30 0){$3$}\htext(40
0){$4$}\htext(50 0){$5$}\htext(60 0){$6$}

\htext(10 -10){$7$}\htext(20 -10){$8$}\htext(30 -10){$9$}\htext(40
-10){$10$}\htext(49 -10){$11$}\htext(59 -10){$12$}

\htext(10 -20){$\vdots$}\htext(20 -20){$\vdots$}\htext(30
-20){$\vdots$}\htext(40 -20){$\vdots$}\htext(50
-20){$\vdots$}\htext(60 -20){$\vdots$}

\move(11 1)\linewd .3  \lcir r:3 \move(31 1)\linewd .3 \lcir r:3
\move(41 1) \linewd .3 \lcir r:3 \fontsize{4}{4}\htext(34 4){$2$}
\end{texdraw}}\vskip 3mm
\end{center}
If we apply $(B_5)$ to \raisebox{-0.25 \height}{\begin{texdraw}
\drawdim em \setunitscale 0.16 \linewd 0.5 \fontsize{8}{8}
\htext(0 0){$3$} \move(1 1)\linewd .3 \lcir r:3
\fontsize{4}{4}\move(4 3)\htext{$2$} \end{texdraw}}, then
\begin{center}
$Y'=$\raisebox{-0.5\height}{
\begin{texdraw}
\drawdim em \setunitscale 0.13 \linewd 0.5

\move(30 0)\lvec(40 0)\lvec(40 10)\lvec(30 10)\lvec(30 0)\htext(32
6){\tiny $1$}

\move(40 0)\lvec(50 0)\lvec(50 10)\lvec(40 10)\lvec(40 0)\htext(42
6){\tiny $0$}
\move(30 0)\lvec(40 10)\lvec(40 0)\lvec(30 0)\lfill f:0.8
\htext(36 2){\tiny $0$}

\move(40 0)\lvec(50 10)\lvec(50 0)\lvec(40 0)\lfill f:0.8
\htext(46 2){\tiny $1$}
\move(40 10)\lvec(50 10)\lvec(50 20)\lvec(40 20)\lvec(40
10)\htext(43 13){$_2$}
\move(40 20)\lvec(50 20)\lvec(50 25)\lvec(40 25)\lvec(40 20)
\htext(43 21){\tiny $3$}

\move(40 25)\lvec(50 25)\lvec(50 30)\lvec(40 30)\lvec(40 25)
\htext(43 26){\tiny $3$}
\end{texdraw}} \ \ \ \ \ $\longleftrightarrow$ \ \ \ \ \
\raisebox{-0.5\height}{\begin{texdraw} \drawdim em \setunitscale
0.16 \linewd 0.5 \fontsize{8}{8}

\htext(10 0){$1$}\htext(20 0){$2$}\htext(30 0){$3$}\htext(40
0){$4$}\htext(50 0){$5$}\htext(60 0){$6$}

\htext(10 -10){$7$}\htext(20 -10){$8$}\htext(30 -10){$9$}\htext(40
-10){$10$}\htext(49 -10){$11$}\htext(59 -10){$12$}

\htext(10 -20){$\vdots$}\htext(20 -20){$\vdots$}\htext(30
-20){$\vdots$}\htext(40 -20){$\vdots$}\htext(50
-20){$\vdots$}\htext(60 -20){$\vdots$}

\move(11 1)\linewd .3  \lcir r:3 \move(41 1) \linewd .3 \lcir r:3
\end{texdraw}}\ \ \ . \vskip 3mm
\end{center}}
\end{ex}\vskip 3mm

\subsection{Characterization of $\Z(\Lambda)_{\Lambda-m\delta}$}

For $Y\in \mathcal{Z}(\Lambda)$, let $\widetilde{Y}$ be the proper
Young wall obtained from $Y$ by applying $(B_i)$ ($1\leq i\leq 5$)
until there is no bead movable up or removable. We set
\begin{equation}
\begin{aligned}
\mathcal{Z}(\Lambda)'&=\{\,Y=(y_k)_{k\geq
1}\in\mathcal{Z}(\Lambda)\,|\,|y_k|\equiv r \pmod \ell\,\}, \\
\mathcal{Z}(\Lambda)'_{\lambda}&=\mathcal{Z}(\Lambda)'\cap
\mathcal{Z}(\Lambda)_{\lambda} \ \ \ \text{for $\lambda\leq
\Lambda$},
\end{aligned}
\end{equation}
where  $r=0$ if $\Lambda=\Lambda_0$ or $\Lambda_1$, and $r=n$ if
$\Lambda=\Lambda_n$. Note that for $Y\in\Z(\Lambda)$,
$Y\in\Z(\Lambda)'$ if and only if there exists no bead in a runner
of type I and II in its bead configuration.

\begin{lem}\label{wtB(n)(1)}
Let $Y$ be a proper Young wall in $\Z(\Lambda)$. Then
$Y\in\mathcal{Z}(\Lambda)_{\Lambda-m\delta}$ for some $m\geq 0$ if
and only if $\widetilde{Y}\in
\mathcal{Z}(\Lambda)'_{\Lambda-m'\delta}$ for some $m'\geq 0$.
\end{lem}
\pf We give a proof for the case $\Lambda=\Lambda_0$ or
$\Lambda_1$ (the proof for the case $\Lambda=\Lambda_n$ is
similar). Let $R_k$ be a runner of type I and let $r_k$ be the
number of beads in $R_k$ occurring in $Y\in\mathcal{Z}(\Lambda)$.
Note that $\widetilde{Y}\in \mathcal{Z}(\Lambda)'$ if and only if
(i) $r_k=r_{\ell-k}$ for $1\leq k\leq n-1$, and (ii) there are
even number of beads in $R_{n}$.

Let $Y$ be a proper Young wall in
$\mathcal{Z}(\Lambda)_{\Lambda-m\delta}$ for some $m\geq 0$. By
considering the content of each bead, we see that ${\rm
cont}(Y)=\sum_{i=0}^n c_i\alpha_i + M\delta$ for some $M\geq 0$,
where
\begin{equation}\label{coeff}
\begin{aligned}
& c_0+c_1=\sum_{j=1}^{2n-1}r_j, \\
&c_i=\sum_{j=i}^{2n-i}r_j + \sum_{j=2n-i+1}^{2n-1}2r_j \ \ \
\text{if $2\leq i \leq n-1$}.
\end{aligned}
\end{equation}
Since ${\rm cont}(Y)=m\delta$ and
$\delta=\alpha_0+\alpha_1+\sum_{i=2}^{n}2\alpha_i$, it follows
that $c_0+c_1=c_2=\cdots=c_{n}$, and hence $r_i=r_{2n-i}$ for
$1\leq i\leq n-1$. Thus, $\widetilde{Y}$ has no bead in a runner
of type I. Suppose that there exists at least one bead in $R_n$ in
the bead configuration of $\widetilde{Y}$. Then by ($B_2$) and
($B_5$), there exists only one bead $b$ at $n$ in $R_n$. On the
other hand, if we consider the content of the blocks corresponding
to all the beads in $R_{\ell}$, we see that the coefficient of
$\alpha_n$ is even. Also, the content of the blocks corresponding
to $b$ is $\alpha_{\epsilon}+\alpha_2+\alpha_3+\cdots+\alpha_n$
($\epsilon=0,1$). This contradicts the fact that the coefficient
of $\alpha_n$ in ${\rm cont}(\widetilde{Y})=m'\delta$ ($m'\geq 0$)
is even. Therefore, we have $\widetilde{Y}\in
\mathcal{Z}(\Lambda)'_{\Lambda-m'\delta}$ for some $m'\geq 0$. The
converse is clear from Lemma \ref{abacusB(n)(1)}. \qed\vskip 5mm

Fix $m\geq 0$. Let $Y$ be a proper Young wall in
$\mathcal{Z}(\Lambda)_{\Lambda-m\delta}$. Consider its bead
configuration. Let $r_k$ be the number of beads in $R_k$ of type
I. By Lemma \ref{wtB(n)(1)}, $r_k=r_{\ell-k}$ for $1\leq k\leq
n-1$. Hence, we can associate a unique partition $\lambda^{(k)}$
from the beads in $R_k$ and $R_{\ell-k}$ ($1\leq k\leq n-1$) using
Frobenius notation. Next, we define a partition $\lambda^{(0)}$ as
follows: If $\Lambda=\Lambda_n$, we set
$\lambda^{(0)}=(1^{m_1},2^{m_2},\cdots)$, where $m_k$ is the
number of beads at $k\ell$ in $R_{\ell}$. If $\Lambda=\Lambda_0$
or $\Lambda_1$, then we set
$\lambda^{(0)}=(1^{m_1},3^{m_3},5^{m_5}\cdots)$, where $m_{2k-1}$
is the number of beads at $(2k-1)n$ in $R_{n}$. In the latter
case, $\ell(\lambda^{(0)})$ is even since the number of beads in
$R_n$ is even by Lemma \ref{wtB(n)(1)}. Also, we see that
$\ell(\lambda^{(0)})$ is even if and only if $|\lambda^{(0)}|$ is
even since $\lambda^{(0)}$ is a partition with each part odd.

We define
\begin{equation}
\pi_0(Y)=(\lambda^{(0)},\cdots,\lambda^{(n-1)}).
\end{equation}

\begin{prop}\label{psi-2} For $Y\in \Z(\Lambda)$, we define
$\psi(Y)=(\pi_0(Y),\widetilde{Y})$.
\begin{itemize}
\item[(1)] If $\Lambda=\Lambda_0$ or $\Lambda_1$, then for $m\geq
0$, the map
\begin{equation*}
\psi : \Z(\Lambda)_{\Lambda-m\delta} \longrightarrow
\bigsqcup_{\sum_{i=1}^3m_i=m}\mathscr{OP}(2m_1)\times
\cP^{(n-1)}(m_2) \times \mathcal{Z}(\Lambda)'_{\Lambda-m_3\delta}
\end{equation*}
is a bijection, where $\mathscr{OP}(k)$ is the set of all
partitions of $k$ with odd parts.

\item[(2)] If $\Lambda=\Lambda_n$, then for $m\geq 0$, the map
\begin{equation*}
\psi : \Z(\Lambda)_{\Lambda-m\delta} \longrightarrow
\bigsqcup_{m_1+m_{2}=m}\cP^{(n)}(m_1)\times
\mathcal{Z}(\Lambda)'_{\Lambda-m_2\delta}
\end{equation*}
is a bijection.
\end{itemize}
\qed
\end{prop}\vskip 3mm

It remains to characterize $\Z(\Lambda)'_{\Lambda-m\delta}$
$(m\geq 0)$. First, suppose that $\Lambda=\Lambda_0$ or
$\Lambda_1$. For each $Y\in\Z(\Lambda)'$, let $\lambda_Y$ be the
two-colored partition in $\cP'$, where the multiplicity of $k$
(resp. $\underline{k}$) is the number of white (resp. gray) beads
at position $k\ell$ in $Y$. We observe that the map $Y\mapsto
\lambda_Y$ is a bijection between $\Z(\Lambda)'$ and $\cP'$. We
define
\begin{equation}
\pi_1(Y)=((\lambda_Y)^{0},(\lambda_Y)^{1}),
\end{equation}
(see Definition \ref{lambda01}), where we view $(\lambda_Y)^{1}$
as an ordinary $2$-reduced partition.

\begin{prop}\label{pi1'}
Suppose that $\Lambda=\Lambda_0$ or $\Lambda_1$. Then for $m\geq
0$, the map
\begin{equation*}
\pi_1: \mathcal{Z}(\Lambda)'_{\Lambda-m\delta} \longrightarrow
\bigsqcup_{m_1+2m_2=m}\mathscr{P}(m_1)\times \mathscr{DP}_0(m_2)
\end{equation*}
is a bijection.
\end{prop}
\pf The proof is almost the same as in the case of
$A_{2n-1}^{(2)}$ (see Proposition \ref{pi1} (1)). \qed\vskip 3mm

Next, suppose that $\Lambda=\Lambda_n$. For each
$Y\in\Z(\Lambda)'$, let $\lambda_Y=\lambda$ be the partition in
$\cP'$, where the multiplicity of $k$ (resp. $\underline{k}$) is
the number of white (resp. gray) beads at the position $(2k-1)n$
in the bead configuration of $Y$. Then the map $Y\mapsto
\lambda_Y$ is a bijection between $\Z(\Lambda_n)'$ and $\cP'$. Fix
$Y\in\Z(\Lambda)'_{\Lambda-m\delta}$ ($m\geq 0$). Let $\mu$ be a
unique partition in $\cP'$ satisfying (i) $\lambda_Y=\mu+\nu$ for
some $\nu\in\cP$, and (ii) the proper Young wall in $\Z(\Lambda)'$
corresponding to $\mu$ is reduced. Note that
$\ell(\lambda)=\ell(\mu)$ and $\mu$ is 2-reduced which has one of
the following forms $(0)$,
$(\underline{1}^{m_{\underline{1}}},2^{m_2},\underline{3}^{m_{\underline{3}}}\cdots)$
($m_{\underline{1}}\neq 0$),
$(1^{m_1},\underline{2}^{m_{\underline{2}}},3^{m_3},\cdots)$
($m_1\neq 0$). We define
\begin{equation}
\pi_1(Y)=(\mu,\nu,c),
\end{equation}
where we view $\mu$ as an ordinary partition and
\begin{equation}
c=
\begin{cases}
0 & \text{ if $\mu=(0)$ or
$\mu=(\underline{1}^{m_{\underline{1}}},2^{m_2},\underline{3}^{m_{\underline{3}}}\cdots)$
with $m_{\underline{1}}\neq 0$,}
\\
1 & \text{ if
$\mu=(1^{m_1},\underline{2}^{m_{\underline{2}}},3^{m_3},\cdots)$
with $m_1\neq 0$.}
\end{cases}
\end{equation}

\begin{lem}\label{mu01nu} Under the above hypothesis, we have
\begin{itemize}
\item[(1)] $\mu\in\mathscr{DP}_0$ and $\ell(\lambda)=\ell(\mu)$ is
even,

\item[(2)] $|\mu|+|\nu|-\ell(\mu)/2=m$.
\end{itemize}
\end{lem}
\pf (1) Since ${\rm cont}(Y)=m\delta$, the number of $n$-blocks in
$Y$ is even. This implies that $\ell(\lambda)$ is even. So it
suffices to show that the 2-core of $\mu$ is empty. Let $Z$ be the
reduced proper Young wall in $\Z(\Lambda)'$ corresponding to
$\mu$. Then we have
\begin{equation}\label{content}
{\rm cont}(Y)={\rm cont}(Z)+|\nu|\delta.
\end{equation}
By definition of $\mu$, it is 2-reduced and
$\ell(\lambda)=\ell(\mu)$.  Consider the two-colored Young diagram
of $\mu$. From the pattern for $\Z(\Lambda_n)$, we see that ${\rm
cont}(Z)=m'\delta$ for some $m'\geq 0$ if and only if the number
of gray boxes with residue 0 is equal to the number of gray boxes
with residue 1. Hence, by a similar argument as in Proposition
\ref{pi1}, we conclude that the 2-core of $\mu$ is empty.

(2) Suppose that ${\rm cont}(Z)=m'\delta$ for some $m'\geq 0$.
Note that $2m'$ is the number of $n$-blocks in $Z$ (except the
ones in $Y_{\Lambda}$). For each part $k$ (or $\underline{k}$) in
$\mu$, the number of $n$-blocks in the corresponding column of $Z$
is $2k-1$. Hence,
\begin{equation}
2m'=\sum(2k-1)(m_k+m_{\underline{k}})=2|\mu|-\ell(\mu),
\end{equation}
where $m_k$ (resp. $m_{\underline{k}}$) is the multiplicity of $k$
(resp. $\underline{k}$) in $\mu$. Since $m=m'+|\nu|$ by
\eqref{content}, we get (2).
 \qed

For $m\geq 0$, let $\mathscr{Q}(m)$ be the set of triples
$(\mu,\nu,c)$ such that
\begin{itemize}
\item[(i)] $\mu\in\mathscr{DP}_0$ with $\ell(\mu)$ even,

\item[(ii)] $\nu\in\cP$ with $\ell(\nu)\leq \ell(\mu)$ and
$|\mu|+|\nu|-\ell(\mu)/2=m$,

\item[(iii)] $c=0$ or $1$, and if $\mu=(0)$, then $c=0$.
\end{itemize}

\begin{lem}\label{q(m)} For each $m\geq 0$, the number of elements in
$\mathscr{Q}(m)$ is the coefficients of $q^m$ in
$(q^2)_{\infty}/(q)^2_{\infty}$.
\end{lem}
\pf Let us fix some notations. For variables $a$ and $q$, set
$(a:q)_k=(1-a)(1-aq)\cdots(1-aq^{k-1})$ ($k\geq 1$), and
$(a:q)_{\infty}=\prod_{k\geq 1}(1-aq^{k-1})$. In particular, set
$(q)_k=(q:q)_k$ and $(q)_{\infty}=(q:q)_{\infty}$.

Let $(\mu,\nu,c)$ be an element in  $\mathscr{Q}(m)$  for some
$m\geq 1$. Let $\mu'$ be the unique 2-reduced partition such that
$\mu=\mu'+(1^{2k})$ where $\ell(\mu)=2k$. Then
$|\mu|-\ell(\mu)/2=|\mu'|+k$, and $\ell(\mu')< 2k$ since $\mu$ is
2-reduced. Note that $\mu$ is uniquely determined by $\mu'$ and
$k$. By considering the bead configurations of the partitions with
empty 2-core, the generating functions for the 2-reduced
partitions with empty 2-core and length less than $2k$ is
\begin{equation}
\frac{(q^2)_{2k-1}}{(q^2)_k(q^2)_{k-1}}
\end{equation}
(we leave the proof to the readers as an exercise. see Section 3).
Also the generating functions for the partitions with length less
than or equal to $2k$ is $1/(q)_{2k}$. Therefore, the generating
functions for $\mathscr{Q}(m)$ ($m\geq 0$) is given by
\begin{equation}
1+2\sum_{k\geq
1}\frac{(q^2)_{2k-1}}{(q^2)_k(q^2)_{k-1}(q)_{2k}}q^{k}.
\end{equation}
Note that 2 appears in the above sum because of $c=0,1$ in
$(\mu,\nu,c)$ ($\mu\neq (0)$). By elementary computation, we have
\begin{equation}
\begin{split}
\frac{(q^2)_{2k-1}}{(q^2)_k(q^2)_{k-1}(q)_{2k}}&=
\frac{\prod_{i=1}^{2k-1}(1+q^i)}{\prod_{j=1}^{k}(1-q^{2j})^2} \\
&=\frac{(-q:q)_k(-q^2:q^2)_{k-2}}{(q^2)^2_k} \\
&=\frac{(-q:q)_k(-1:q^2)_{k}}{2(q^2)^2_k}.
\end{split}
\end{equation}
Therefore, the generating function for $\mathscr{Q}(m)$ is
\begin{equation}
\sum_{k\geq 0}\frac{(-q:q)_k(-1:q^2)_{k}}{(q^2)^2_k}.
\end{equation}
If we apply a Heine's identity (see Corollary 2.4 in \cite{And})
by replacing $a$, $b$, $c$, $q$ by $-q$, $-1$, $q^2$, $q^2$,
respectively, we obtain
\begin{equation}
\begin{split}
\sum_{k\geq 0}\frac{(-q:q)_k(-1:q^2)_{k}}{(q^2)^2_k}&=
\frac{(-q:q)_{\infty}(-q^2:q^2)_{\infty}}{(q:q^2)_{\infty}(q^2)_{\infty}}
\\
&=\frac{(-q:q)_{\infty}}{(q)_{\infty}}=\frac{(q^2)_{\infty}}{(q)_{\infty}^2},
\end{split}
\end{equation}
which proves our claim. \qed

\begin{prop}\label{pi1''}
Suppose that $\Lambda=\Lambda_n$. For $m\geq 0$, the map
\begin{equation*}
\pi_1: \mathcal{Z}(\Lambda)'_{\Lambda-m\delta} \longrightarrow
\mathscr{Q}(m)
\end{equation*}
is a bijection.
\end{prop}
\pf By Lemma \ref{mu01nu}, $\pi_1$ is well-defined. Note that $Y$
is uniquely determined by $\pi_1(Y)$. Suppose that
$(\mu,\nu,c)\in\mathscr{Q}(m)$ is given. If $c=0$, then we view
$\mu$ as a two-colored partition where the color of an even (resp.
odd) part is white (resp. gray). If $c=1$, then we view $\mu$ as a
two-colored partition where the color of an even (resp. odd) part
is gray (resp. white). Set $\lambda=\mu+\nu\in\cP'$. Then there
exists a unique $Y\in\Z(\Lambda)'$ corresponding to $\lambda$. By
the same argument in Lemma \ref{mu01nu}, we can check that
$Y\in\Z(\Lambda)'_{\Lambda-m\delta}$ with
$|\mu|+|\nu|-\ell(\mu)/2=m$, and $\pi_1(Y)=(\mu,\nu,c)$, which
defines the inverse map of $\pi_1$. \qed\vskip 3mm

\begin{ex}{\rm
Let $Y$ be a proper Young wall in $\Z(\Lambda_n)'$ such that the
corresponding partition is
$\lambda_Y=(1^2,\underline{2},\underline{3}^2,\underline{4},6,7)$.
If we take
\begin{equation}
\mu=(1^2,\underline{2}^4,3^2), \ \ \ \nu=(1^2,2,3,4),
\end{equation}
then $\lambda_Y=\mu+\nu$ where $\mu$ is the unique $2$-colored
partition whose corresponding proper Young wall in
$\Z(\Lambda_n)'$ is reduced. Also $\mu$ and $\nu$ satisfy the
conditions in Lemma \ref{mu01nu}. Since
$|\mu|+|\nu|-\ell(\mu)/2=23$, we have $Y\in
\Z(\Lambda_n)'_{\Lambda_n-23\delta}$ and
\begin{equation}
\pi_1(Y)=((1^2,2^4,3^2),(1^2,2,3,4),1).
\end{equation}
}
\end{ex}\vskip 5mm

Now, for $m\geq 0$ and $Y\in
\mathcal{Z}(\Lambda)_{\Lambda-m\delta}$
($\Lambda=\Lambda_0,\Lambda_1,\Lambda_n$), we define
\begin{equation}\label{pi}
\pi(Y)=(\pi_0(Y),\pi_1(\widetilde{Y})).
\end{equation}
Then by Proposition \ref{psi-2}, \ref{pi1'} and \ref{pi1''}, we
obtain
\begin{thm}\mbox{}
\begin{itemize}
\item [(1)] If $\Lambda=\Lambda_0$ or $\Lambda_1$, then for $m\geq
0$, the map
\begin{equation*}
\pi : \Z(\Lambda)_{\Lambda-m\delta} \longrightarrow
\bigsqcup_{m_1+m_2+2m_3=m}\mathscr{OP}(2m_1)\times \cP^{(n)}(m_2)
\times \mathscr{DP}_0(m_3)
\end{equation*}
is a bijection.

\item[(2)] If $\Lambda=\Lambda_n$, then for $m\geq 0$, the map
\begin{equation*}
\pi : \Z(\Lambda)_{\Lambda-m\delta} \longrightarrow
\bigsqcup_{m_1+m_{2}=m}\cP^{(n)}(m_1)\times \mathscr{Q}(m_2)
\end{equation*}
is a bijection.
\end{itemize}
\qed
\end{thm}

\begin{ex}{\rm (1) Suppose that $\Lambda=\Lambda_0$.\vskip 5mm

\begin{center}
\hskip -1cm If \ \  $Y=$\ \
\raisebox{-0.5\height}{\begin{texdraw} \drawdim em \setunitscale
0.16 \linewd 0.5 \fontsize{8}{8}

\htext(10 0){$1$}\htext(20 0){$2$}\htext(30 0){$3$}\htext(40
0){$4$}\htext(50 0){$5$} \htext(60 0){$6$}

\htext(10 -10){$7$}\htext(20 -10){$8$}\htext(30 -10){$9$}\htext(40
-10){$10$}\htext(49 -10){$11$}  \htext(59 -10){$12$}

\htext(10 -20){$13$}\htext(20 -20){$14$}\htext(30
-20){$15$}\htext(40 -20){$16$}\htext(49 -20){$17$} \htext(59
-20){$18$}

\htext(10 -30){$19$}\htext(20 -30){$20$}\htext(30
-30){$21$}\htext(40 -30){$22$}\htext(49 -30){$23$} \htext(59
-30){$24$}

\htext(10 -40){$\vdots$}\htext(20 -40){$\vdots$}\htext(30
-40){$\vdots$}\htext(40 -40){$\vdots$}\htext(50
-40){$\vdots$}\htext(60 -40){$\vdots$}

\move(12 -29)\linewd .3 \lcir r:3.5

\move(22 1)\linewd .3 \lcir r:3.5

\move(32 1)\linewd .3 \lcir r:3.5

\move(33 -19)\linewd .3 \lcir r:3.5

\move(33 -29)\linewd .3 \lcir r:3.5

\move(43 -19)\linewd .3 \lcir r:3.5

\move(52 -9)\linewd .3  \lcir r:3.5

\move(62 1)\linewd .3 \fcir f:0.7 r:3.5 \lcir r:3.5

\move(62 -9)\linewd .3 \fcir f:0.7 r:3.5 \lcir r:3.5

\move(62 -29)\linewd .3 \lcir r:3.5

\fontsize{4}{4}\move(34 4)\htext{$2$} \move(34 -16)\htext{$3$}
\move(64 -6)\htext{$4$} \move(64 -26)\htext{$2$}
\end{texdraw}} \ \ , \ \ then  \ \ \ $\widetilde{Y}=$ \ \ \
\raisebox{-0.5\height}{\begin{texdraw} \drawdim em \setunitscale
0.16 \linewd 0.5 \fontsize{8}{8}

\htext(50 0){$6$}

\htext(49 -10){$12$}

\htext(49 -20){$18$}

\htext(49 -30){$24$}

\htext(50 -40){$\vdots$}

\move(52 1)\linewd .3 \lcir r:3.5

\move(52 -9)\linewd .3 \fcir f:0.7 r:3.5 \lcir r:3.5


\move(52 -29)\linewd .3  \lcir r:3.5

\fontsize{4}{4}\move(54 -6)\htext{$4$} \move(54 -26)\htext{$2$}
\end{texdraw}}\end{center}\vskip 5mm
and  $\pi_0(Y)=(\lambda^{(0)},\lambda^{(1)},\lambda^{(2)})$, where
\begin{equation}
\begin{split}
\lambda^{(0)}&=(1^2,5^3,7), \\
\lambda^{(1)}&=((3)|(1))=(1^3,2), \\
\lambda^{(2)}&=((0)|(2))=(3).
\end{split}
\end{equation}
Note that $\lambda_{\widetilde{Y}}=(1,\underline{2}^4,4^2)$ with
$(\lambda_{\widetilde{Y}})^0=(1^5,2^2)$ and
$(\lambda_{\widetilde{Y}})^1=(\underline{1}^4,2^2)$. Since the
$2$-core of $(\lambda_{\widetilde{Y}})^1$ is empty and
$|\lambda_{\widetilde{Y}}|=17$, we have
$\widetilde{Y}\in\Z(\Lambda)'_{\Lambda-17\delta}$ by Proposition
\ref{pi1'} (1), and $\pi_1(\widetilde{Y})=((1^5,2^2),(1^4,2^2))$.
Hence, we have $Y\in\Z(\Lambda)_{\Lambda-37\delta}$ and
\begin{equation}
\pi(Y)=((1^2,5^3,7),(1^3,2),(3),(1^5,2^2),(1^4,2^2)).
\end{equation}

(2) Suppose that $\Lambda=\Lambda_n$.\vskip 5mm

\begin{center}
\hskip -1cm If \ \  $Y=$\ \ \raisebox{-0.5\height}{\begin{texdraw}
\drawdim em \setunitscale 0.16 \linewd 0.5 \fontsize{8}{8}

\htext(10 0){$1$}\htext(20 0){$2$}\htext(30 0){$3$}\htext(40
0){$4$}\htext(50 0){$5$} \htext(60 0){$6$}

\htext(10 -10){$7$}\htext(20 -10){$8$}\htext(30 -10){$9$}\htext(40
-10){$10$}\htext(49 -10){$11$}  \htext(59 -10){$12$}

\htext(10 -20){$13$}\htext(20 -20){$14$}\htext(30
-20){$15$}\htext(40 -20){$16$}\htext(49 -20){$17$} \htext(59
-20){$18$}

\htext(10 -30){$19$}\htext(20 -30){$20$}\htext(30
-30){$21$}\htext(40 -30){$22$}\htext(49 -30){$23$} \htext(59
-30){$24$}

\htext(10 -40){$\vdots$}\htext(20 -40){$\vdots$}\htext(30
-40){$\vdots$}\htext(40 -40){$\vdots$}\htext(50
-40){$\vdots$}\htext(60 -40){$\vdots$}

\move(12 -19)\linewd .3 \lcir r:3.5

\move(22 1)\linewd .3 \lcir r:3.5

\move(32 1)\linewd .3  \lcir r:3.5

\move(33 -19)\linewd .3  \lcir r:3.5

\move(33 -29)\linewd .3 \fcir f:0.7 r:3.5 \lcir r:3.5

\move(32 -9)\linewd .3  \lcir r:3.5

\move(43 -29)\linewd .3 \lcir r:3.5

\move(52 -9)\linewd .3  \lcir r:3.5

\move(62 -19)\linewd .3 \lcir r:3.5

\move(62 -9)\linewd .3 \lcir r:3.5

\move(62 -29)\linewd .3 \lcir r:3.5

\fontsize{4}{4}\move(34 4)\htext{$2$} \move(34 -6)\htext{$2$}
\move(34 -16)\htext{$3$}

\fontsize{4}{4} \move(64 -16)\htext{$3$}
\end{texdraw}} \ \ , \ \ then  \ \ \ $\widetilde{Y}=$ \ \ \
\raisebox{-0.5\height}{\begin{texdraw} \drawdim em \setunitscale
0.16 \linewd 0.5 \fontsize{8}{8} \htext(50 0){$3$}

\htext(50 -10){$9$}

\htext(49 -20){$15$}

\htext(49 -30){$21$}

\htext(50 -40){$\vdots$}

\move(52 1)\linewd .3  \lcir r:3.5

\move(52 -9)\linewd .3  \lcir r:3.5

\move(52 -19)\linewd .3 \fcir f:0.7 r:3.5 \lcir r:3.5

\move(52 -29)\linewd .3 \fcir f:0.7 r:3.5 \lcir r:3.5

\fontsize{4}{4}\move(54 -6)\htext{$2$} \move(54 -16)\htext{$3$}
\move(54 6)\htext{$2$}
\end{texdraw}}\end{center}\vskip 5mm
and $\pi_0(Y)=(\lambda^{(0)},\lambda^{(1)},\lambda^{(2)})$, where
\begin{equation}
\begin{split}
\lambda^{(0)}&=(2,3^3,4), \\
\lambda^{(1)}&=((2)|(1))=(1^2,2), \\
\lambda^{(2)}&=((0)|(3))=(4).
\end{split}
\end{equation}

Also,
$\lambda_{\widetilde{Y}}=(1^2,2^2,\underline{3}^3,\underline{4})$.
If we take $\mu=(1^4,\underline{2}^4)$ and $\nu=(1^5,2)$, then
$\lambda_{\widetilde{Y}}=\mu+\nu$ and the proper Young wall in
$\Z(\Lambda_n)'$ corresponding to $\mu$ is reduced. Since
$((1^4,2^4),(1^5,2),1)\in\mathscr{Q}(15)$, we have
$\widetilde{Y}\in \Z(\Lambda_n)'_{\Lambda_n-15\delta}$ by
Proposition \ref{pi1''}, and
$\pi_1(\widetilde{Y})=((1^4,2^4),(1^5,2),1)$. Hence, we have
$Y\in\Z(\Lambda)_{\Lambda-38\delta}$ and
\begin{equation}
\pi(Y)=((2,3^3,4),(1^2,2),(4),(1^4,2^4),(1^5,2),1).
\end{equation}}
\end{ex}

Therefore, we obtain a new combinatorial proof of the formulas in
\cite{Kac90}.
\begin{cor}\label{Bn1cor1}\mbox{}

\begin{itemize}
\item [(1)] If $\Lambda=\Lambda_0$ or $\Lambda_1$, then we have
\begin{equation*}
\Sigma^{\Lambda}_{\Lambda}(q)=\dfrac{1}{2}\left(
\dfrac{1}{(q^{\frac{1}{2}})_{\infty}(q)^{n-1}_{\infty}(q^2)_{\infty}}+
\dfrac{(q^{\frac{1}{2}})_{\infty}}{(q)^{n+2}_{\infty}} \right).
\end{equation*}

\item[(2)] If $\Lambda=\Lambda_n$, then we have
$\Sigma^{\Lambda}_{\Lambda}(q)=\dfrac{(q^2)_{\infty}}{(q)^{n+2}_{\infty}}$.
\end{itemize}
\end{cor}
\pf (1) Let $\mathscr{O}(m)$ be the number of the partitions of
$m$ with odd parts. Then we have $O(q)=\sum_{m\geq
0}|\mathscr{O}(m)|q^m=(-q:q)_{\infty}=1/(q:q^2)_{\infty}$
($|\mathscr{O}(0)|=1$). Hence the number of the partitions in
$\mathscr{O}(2m)$ is the coefficient of $q^m$ in
\begin{equation}
\frac{1}{2}\left ( O(q^{\frac{1}{2}})+O(-q^{\frac{1}{2}})
\right)=\dfrac{1}{2}\left(
\dfrac{(q)_{\infty}}{(q^{\frac{1}{2}})_{\infty}}+
\dfrac{(q^{\frac{1}{2}})_{\infty}(q^2)_{\infty}}{(q)^{2}_{\infty}}
\right).
\end{equation} Hence, we obtain $\Sigma^{\Lambda}_{\Lambda}(q)$ by multiplying
$\frac{1}{(q)_{\infty}^n(q^2)_{\infty}}$.

(2) This follows from Lemma \ref{q(m)}. \qed

\subsection{Characterization of $\Z(\Lambda_0)_{\Lambda_1-m\delta}$}

\begin{lem}\label{abacusB(n)(1)'}
For $Y\in \Z(\Lambda_0)$, let $r_k$ be the number of beads in
$R_k$ in the bead configuration of $Y$ {\rm (}$1\leq k < \ell,
k\neq n${\rm )}. If $Y\in \Z(\Lambda_0)_{\Lambda_1-m\delta}$ for
some $m\geq 0$, then $r_k=r_{\ell-k}$ for $1\leq k\leq n-1$.
\end{lem}
\pf We see that ${\rm cont}(Y)=\sum_{i=0}^n c_i\alpha_i + M\delta$
for some $M\geq 0$ where $c_i$ is given by \eqref{coeff}. On the
other hand, ${\rm cont}(Y)=\gamma +
m\delta=(m+1)\alpha_0+m\alpha_1+(2m+1)\sum_{i=2}^n\alpha_i$, where
$\gamma=\alpha_0+\sum_{i=2}^n\alpha_i$. This implies that
$c_0+c_1=c_2=\cdots=c_{n}$, and hence $r_k=r_{\ell-k}$ for $1\leq
k\leq n-1$. \qed\vskip 3mm

Suppose that $Y\in \Z(\Lambda_0)_{\Lambda_1-m\delta}$ is given. By
Lemma \ref{abacusB(n)(1)'}, we can associate a unique partition
$\lambda^{(k)}$ ($1\leq k\leq n-1$) from the beads in $R_k$ and
$R_{\ell-k}$ using Frobenius notation.

Let $Y'$ be the proper Young wall obtained by applying $(B_1)$ and
$(B_3)$ to $Y$ until there is no bead in the runners of type I.
Note that $Y'\in \Z(\Lambda_0)_{\Lambda_1-m'\delta}$, where
$m'=m-\sum_{i=1}^{n-1}|\lambda^{(i)}|$. Set
$\lambda^{(0)}=(1^{m_1},3^{m_3},5^{m_5},\cdots)$, where $m_{2k-1}$
is the number of the beads at $(2k-1)n$. Since ${\rm
cont}(Y')=\gamma+m'\delta$, the number of $n$-blocks in $Y'$
(except the ones in $Y_{\Lambda}$) is odd and hence the number of
beads in $R_n$ is odd (or $|\lambda^{(0)}|$ is odd).  We define
\begin{equation}
\pi_0(Y)=(\lambda^{(0)},\cdots,\lambda^{(n-1)}).
\end{equation}

Next, consider $\widetilde{Y'}=\widetilde{Y}$. Note that
$\widetilde{Y}\in\Z(\Lambda_0)_{\Lambda_1-m''\delta}$, where
$m''=m'-\frac{|\lambda^{(0)}|-1}{2}$. From the beads of
$\widetilde{Y}$ in $R_{\ell}$, we define $\mu$ to be the partition
in $\cP'$, where the multiplicity of $k$ (resp. $\underline{k}$)
is given by
\begin{equation}\label{mu}
\begin{split}
m_{k}&=
\begin{cases}
\text{the number of white beads at $k\ell$} & \text{ if
$\ell(|\widetilde{Y}|)$ is odd},\\
\text{the number of gray beads at $k\ell$} & \text{ if
$\ell(|\widetilde{Y}|)$ is even},
\end{cases}
\\
m_{\underline{k}}&=
\begin{cases} \text{the number of gray beads at
$k\ell$} & \text{ if
$\ell(|\widetilde{Y}|)$ is odd},\\
\text{the number of white beads at $k\ell$} & \text{ if
$\ell(|\widetilde{Y}|)$ is even},
\end{cases}
\end{split}
\end{equation}

We define
\begin{equation}
\pi_1(\widetilde{Y})=(\mu^0,\mu^1),
\end{equation}
where we view $\mu^1$ as an ordinary $2$-reduced partition.

\begin{lem}\label{mu01} Under the above hypothesis, we have
\begin{itemize}
\item[(1)] $\mu^1\in \mathscr{DP}_0$,

\item[(2)] $|\mu^0|+|\mu^1|=m''$, where
$\widetilde{Y}\in\Z(\Lambda_0)_{\Lambda_1-m''\delta}$.
\end{itemize}
\end{lem}
\pf First, note that in $\widetilde{Y}$, there exists exactly one
bead $b$ in $R_n$ and it is placed at $n$ since the number of the
beads in $R_n$ in the bead configuration of $Y'$ is odd. Since
there is an one-to-one correspondence between $\cP'$ and
$\Z(\Lambda)'$ ($\Lambda=\Lambda_0,\Lambda_1$), there exist unique
 proper Young walls $Z\in \Z(\Lambda_r)'$
corresponding to $\mu$ where $r\equiv \ell(|\widetilde{Y}|)+1
\pmod 2$. In fact, $Z$ can be obtained in the following way:
\begin{itemize}
\item[(i)] if $\ell(|\widetilde{Y}|)$ is odd, then remove the
left-most column of $\widetilde{Y}$ whose content is $\gamma$. The
resulting proper Young wall is $Z\in\Z(\Lambda_0)'$.

\item[(ii)] if $\ell(|\widetilde{Y}|)$ is even, then shift all the
blocks to the left column by one position following the pattern,
except the first $n$ blocks (except the one in $Y_{\Lambda}$) from
the bottom. There are $n$ blocks left in the first column whose
content is $\gamma$. If we cut out this first column, then the
resulting proper Young wall is $Z\in\Z(\Lambda_1)'$.
\end{itemize}
From the above facts, we have
\begin{equation}
{\rm cont}(\widetilde{Y})={\rm cont}(Z)+ \gamma,
\end{equation}
where $\gamma=\alpha_{0}+\sum_{i=2}^n\alpha_i$. Since ${\rm
cont}(\widetilde{Y})=m''\delta+\gamma$, we have ${\rm
cont}(Z)=m''\delta$, i.e. $Z\in
\Z(\Lambda_r)'_{\Lambda_r-m''\delta}$. Hence from Proposition
\ref{pi1'}, $\mu$ satisfies the conditions in (1) and (2).\qed
\vskip 3mm

Now, for each $m\geq 0$ and $Y\in
\mathcal{Z}(\Lambda_0)_{\Lambda_1-m\delta}$, we define
\begin{equation}\label{pi'}
\pi(Y)=(\pi_0(Y),\pi_1(\widetilde{Y})).
\end{equation}

Then, we obtain
\begin{thm}
For each $m\geq 0$, the map
\begin{equation*}
\pi : \Z(\Lambda_0)_{\Lambda_1-m\delta} \longrightarrow
\bigsqcup_{m_1+m_2+2m_3=m}\mathscr{OP}(2m_1+1)\times\cP^{(n)}(m_2)
\times \mathscr{DP}_0(m_3)
\end{equation*}
is a bijection.
\end{thm}
\pf By Lemma \ref{mu01}, $\pi$ is well-defined. The inverse map
can be defined by reversing the construction of $\pi$ naturally.
\qed

\begin{ex}{\rm Consider the following proper Young wall in
$\Z(\Lambda_0)$.\vskip 5mm

\begin{center}
  \ \  $Y=$\ \  \raisebox{-0.5\height}{\begin{texdraw}
\drawdim em \setunitscale 0.16 \linewd 0.5 \fontsize{8}{8}

\htext(10 0){$1$}\htext(20 0){$2$}\htext(30 0){$3$}\htext(40
0){$4$}\htext(50 0){$5$} \htext(60 0){$6$}

\htext(10 -10){$7$}\htext(20 -10){$8$}\htext(30 -10){$9$}\htext(40
-10){$10$}\htext(49 -10){$11$}  \htext(59 -10){$12$}

\htext(10 -20){$13$}\htext(20 -20){$14$}\htext(30
-20){$15$}\htext(40 -20){$16$}\htext(49 -20){$17$} \htext(59
-20){$18$}

\htext(10 -30){$19$}\htext(20 -30){$20$}\htext(30
-30){$21$}\htext(40 -30){$22$}\htext(49 -30){$23$} \htext(59
-30){$24$}

\htext(10 -40){$\vdots$}\htext(20 -40){$\vdots$}\htext(30
-40){$\vdots$}\htext(40 -40){$\vdots$}\htext(50
-40){$\vdots$}\htext(60 -40){$\vdots$}

\move(12 -29)\linewd .3 \lcir r:3.5

\move(22 1)\linewd .3 \lcir r:3.5

\move(32 1)\linewd .3 \lcir r:3.5

\move(33 -19)\linewd .3 \lcir r:3.5

\move(33 -29)\linewd .3 \lcir r:3.5

\move(43 -19)\linewd .3 \lcir r:3.5

\move(52 -9)\linewd .3  \lcir r:3.5

\move(62 1)\linewd .3  \lcir r:3.5

\move(62 -9)\linewd .3 \lcir r:3.5

\move(62 -29)\linewd .3 \lcir r:3.5

\fontsize{4}{4}\move(34 4)\htext{$2$} \move(34 -16)\htext{$2$}
\move(64 -6)\htext{$4$} \move(64 -26)\htext{$3$}
\end{texdraw}}\end{center}\vskip 5mm

Then we have  \ \ \ $\widetilde{Y}=$ \ \ \
\raisebox{-0.5\height}{\begin{texdraw} \drawdim em \setunitscale
0.16 \linewd 0.5 \fontsize{8}{8}

\htext(39 0){$3$}

\htext(39 -10){$9$}

\htext(39 -20){$15$}

\htext(39 -30){$21$}

\htext(40 -40){$\vdots$}

\htext(50 0){$6$}

\htext(49 -10){$12$}

\htext(49 -20){$18$}

\htext(49 -30){$24$}

\htext(50 -40){$\vdots$}

\move(41 1)\linewd .3  \lcir r:3.5

\move(52 1)\linewd .3  \lcir r:3.5

\move(52 -9)\linewd .3 \fcir f:0.7 r:3.5 \lcir r:3.5


\move(52 -29)\linewd .3  \lcir r:3.5

\fontsize{4}{4}\move(54 -6)\htext{$4$} \move(54 -26)\htext{$3$}
\end{texdraw}} \ \ \
and \ \ \  $\pi_0(Y)=(\lambda^{(0)},\lambda^{(1)},\lambda^{(2)})$,
where \vskip 5mm
\begin{equation}
\begin{split}
\lambda^{(0)}&=(1^2,5^2,7), \\
\lambda^{(1)}&=((3)|(1))=(1^3,2), \\
\lambda^{(2)}&=((0)|(2))=(3).
\end{split}
\end{equation}
Note that $\ell(|\widetilde{Y}|)$ is odd. Following the rule in
\eqref{mu}, we get $\mu=(1,\underline{2}^4,4^3)$, where
$\mu^0=(1^5,2^3)$ and $\mu^1=(\underline{1}^4,2^3)$. Since $\mu^1$
has an empty $2$-core, we have
$\widetilde{Y}\in\Z(\Lambda_0)_{\Lambda_1-21\delta}$ and
$\pi_1(\widetilde{Y})=((1^5,2^3),(1^4,2^3))$. Therefore,
$Y\in\Z(\Lambda_0)_{\Lambda_1-38\delta}$ and
\begin{equation}
\pi(Y)=((1^2,5^2,7),(1^3,2),(3),(1^5,2^3),(1^4,2^3)).
\end{equation}
}
\end{ex}

Therefore, we also obtain another proof of the formula in
\cite{Kac90}.
\begin{cor}
\begin{equation*}
\Sigma^{\Lambda_0}_{\Lambda_1}(q)=\dfrac{1}{2q^{\frac{1}{2}}}\left(
\dfrac{1}{(q^{\frac{1}{2}})_{\infty}(q)^{n-1}_{\infty}(q^2)_{\infty}}-
\dfrac{(q^{\frac{1}{2}})_{\infty}}{(q)^{n+2}_{\infty}} \right).
\end{equation*}
\end{cor}
\pf As in Corollary \ref{Bn1cor1} (1), the number of partitions in
$\mathscr{O}(2m+1)$ ($m\geq 0$) is the coefficient of $q^m$ in
\begin{equation}
\frac{1}{2q^{\frac{1}{2}}}\left(
O(q^{\frac{1}{2}})-O(-q^{\frac{1}{2}})
\right)=\dfrac{1}{2q^{\frac{1}{2}}}\left(
\dfrac{(q)_{\infty}}{(q^{\frac{1}{2}})_{\infty}}-
\dfrac{(q^{\frac{1}{2}})_{\infty}(q^2)_{\infty}}{(q)^{2}_{\infty}}
\right).
\end{equation} Hence, we obtain the result.\qed\end{document}